
\magnification=\magstep1
\baselineskip=18pt
\parindent 18pt

\hsize=5.9truein
\hoffset=0.5truein
 

\def\authors{b. j. c. baxter}
\def\chptitle{}



\font\sc=cmcsc10

\font\headfont=cmsl9

\def\pf{{\noindent{\it Proof.\ }}}
\def\square{\vrule height5pt width5pt}
\def\qed{{\hfill$\square$\vskip 10pt}}

\def\R1{{\cal R}}
\def\Rd{{\R1^d}}
\def\Rn{{\R1^n}}
\def\ZZ{{{\cal Z} }}
\def\Z0{{\sum_{j=1}^n y_j = 0}}
\def\Zd{{\ZZ^d}}
\def\T1{{[0,2\pi]}}
\def\Td{{[0,2\pi]^d}}
\def\CC{{\cal C}}

\def\half{{1\over2}}
\def\OO{{\cal O}}
\def\phi{{\varphi}}
\def\phihat{{\hat\varphi}}

\outer\def\beginsection#1\par{\vskip0pt plus.3\vsize\penalty-150
 \vskip0pt plus-.3\vsize\bigskip\vskip\parskip
 \message{#1}\centerline{\bf#1}\nobreak\smallskip\noindent}

\def\proclaim #1. #2\par{\medbreak
\noindent{\bf#1.\enspace}{\it#2}\par\medbreak}

\def\proclaimdefn #1. #2\par{\medbreak
\noindent{\bf#1.\enspace}#2\par\medbreak}



\def\sect#1{\goodbreak\bigskip\smallskip\centerline{\bf #1}
\medskip\noindent\ignorespaces}


\def\rightheadline{\ifnum\pageno=1 \hfill%
  \else\sl\hfil\chptitle\hfil\fi}

\def\leftheadline{\ifnum\pageno=1 \hfill%
  \else\hfil\sc\authors\hfil\fi}  


\def\pts{{(x_j)_{j=1}^n}}

\def\sjk{{\sum_{j,k=1}^n y_j y_k\,}}
\def\djk{{\| x_j - x_k \|}}
\def\xjk{{x_j-x_k}}
\def\yjk{{|\sum_{j=1}^n y_j e^{ix_j t}|^2}}
\def\trig{{|\sum_{j=1}^n y_j e^{ijt}|^2}}

\def\skzd{{\sum_{k\in\Zd}}}

\def\rint{{\int_{-\infty}^{\infty}}}
\def\hint{{\int_0^\infty}}

\def\Sdm1{{S^{d-1}}}

\def\ysqr{{\Vert y \Vert^2}}

\def\ainv{{A_n^{-1}}}
\def\aninv{{\Vert A_n^{-1} \Vert_2}}

\def\zsum{{\sum_{k=-\infty}^\infty}}

\def\hi{{(2\pi)^{-1}}}
\def\tint{{\hi\int_0^{2\pi}}}
\def\tintd{{(2\pi)^{-d}\int_{\Td}}}
\def\fhat{{\hat f}}
\def\ajxi{{\Bigl|\sum_{j\in\Zd} a_j \exp(ij\xi)\Bigl|^2}}
\def\Rdm0{{\Rd\setminus\{0\}}}
\def\ghat{{\hat g}}
\def\ytrigo{{\Bigl|\sum_{j\in\Zd} y_j \exp(ix_j \xi)\Bigl|^2}}
\def\ytrig2{{\Bigl|\sum_{j\in\Zd} y_j \exp(ij \xi)\Bigl|^2}}
\def\ynj{{y_j^{(n)}}}

\def\phichat{{\phihat_c}}
\def\chichat{{{\hat \chi}_c}}

\def\min{{\rm min}}

\def\Rn{{\R1^n}}
\def\Rnn{{\R1^{n \times n}}}

\input psfig.sty


\pageno=-1
\centerline{\bf THE INTERPOLATION THEORY}
\centerline{\bf OF RADIAL BASIS FUNCTIONS}

\medskip
\centerline{by}
\medskip
\centerline{\sc Bradley John Charles Baxter}
\centerline{of}
\centerline  {\sc Trinity College}
\bigskip
\bigskip
\centerline{A dissertation presented in fulfilment of the requirements for}
\centerline{the degree of Doctor of Philosophy, Cambridge University}
\medskip
\centerline{\sl August 1992}

\vfill \eject
\centerline{\bf THE INTERPOLATION THEORY}
\centerline{\bf OF RADIAL BASIS FUNCTIONS}

\medskip
\centerline{\sc B. J. C. Baxter}
\medskip
\beginsection Summary

The problem of interpolating functions of $d$ real
variables ($d > 1$) occurs naturally in many areas of
applied mathematics and the sciences. Radial basis
function methods can provide interpolants to function values
given at irregularly positioned points for any value of
$d$. Further, these interpolants are often excellent
approximations to the underlying function, even when
the number of interpolation points is small.

In this dissertation we begin with the existence
theory of radial basis function interpolants. It is first shown that, when the
radial basis function is a $p$-norm and $1 < p < 2$,
interpolation is always possible when the 
points are all different and there are at least two of
them. Our approach extends the analysis of the case
$p=2$ devised in the 1930s by Schoenberg. We then show
that interpolation is not always possible when $p > 2$.
Specifically, for every $p > 2$, we construct a set of
different points in some $\Rd$ for which the
interpolation matrix is singular. This construction seems
to have no precursor in the literature.

The greater part of this work investigates the sensitivity of radial basis
function interpolants to changes in the function values at the interpolation
points. This study was motivated by the observation that large condition
numbers occur in some practical calculations. Our early results show that it is
possible to recast the work of Ball, Narcowich and Ward in the language of
distributional Fourier transforms in an elegant way. We then use this language
to study the interpolation matrices generated by subsets of regular grids. In
particular, we are able to extend the classical theory of Toeplitz operators to
calculate sharp bounds on the spectra of such matrices. Moreover, we also
describe some joint work with Charles Micchelli in which we use the theory of
P\'olya frequency functions to continue this work, as well as shedding new
light on some of our earlier results.

Applying our understanding of these spectra, we construct
preconditioners for the conjugate gradient solution of the interpolation
equations. The preconditioned conjugate gradient algorithm was first suggested
for this problem by Dyn, Levin and Rippa in 1983, who were motivated by the
variational theory of the thin plate spline. In contrast, our approach is
intimately connected to the theory of Toeplitz forms. Our main result is that the
number of steps required to achieve solution of the linear system to within a
required tolerance can be independent of the number of interpolation points. In
other words, the number of floating point operations needed for a regular grid is
proportional to the cost of a matrix-vector multiplication. The Toeplitz
structure allows us to use fast Fourier transform techniques, which implies that
the total number of operations is a multiple of $n \log n$, where $n$ is the
number of interpolation points.  

Finally, we use some of our methods to study the behaviour of the multiquadric
when its shape parameter increases to infinity. We find a surprising link with
the {\it sinus cardinalis} or {\it sinc} function of Whittaker. Consequently, it
can be highly useful to use a large shape parameter when approximating
band-limited functions.

\vfill \eject
\beginsection Declaration


In this dissertation, all of the work is my own with the exception of
Chapter~5, which contains some results of my collaboration with Dr
Charles Micchelli of the IBM Research Center, Yorktown Heights, New
York, USA. This collaboration was approved by the Board of Graduate
Studies.

No part of this thesis has been submitted for a degree elsewhere.
However, the contents of several chapters have appeared, or are to appear, in
journals. In particular, we refer the reader to Baxter (1991a, b, c)
and Baxter (1992a, b). Furthermore, Chapter 2 formed a Smith's Prize essay
in an earlier incarnation.

\vfill \eject
\beginsection Preface

It is a pleasure to acknowledge the support I have received during my
doctoral research.

First, I must record my gratitude to Professor Michael Powell for his
support, patience and understanding whilst supervising my studies.
His enthusiasm, insight, precision, and distrust for gratuitous
abstraction have enormously influenced my development as a
mathematician. In spite of his many commitments he has always been
generous with his time. In particular, I am certain that the great
care he has exhibited when reading my work will be of lasting benefit;
there could be no better training for the preparation and refereeing
of technical papers. 

The Numerical Analysis Group of the University of Cambridge has
provided an excellent milieu for research, but I am especially
grateful to Arieh Iserles, whose encouragement and breadth of
mathematical knowledge have been of great help to me. In particular, it was
Arieh who introduced me to the beautiful theory of Toeplitz
operators.

Several institutions have supported me financially.
The Science and Engineering Research Council
and A.E.R.E. Harwell provided me with a CASE Research Studentship
during my first three years. At this point, I must thank Nick Gould
for his help at Harwell.
Subsequently I have been aided by
Barrodale Computing Limited, the Amoco Research Company, Trinity
College, Cambridge, and the endowment of the John Humphrey Plummer
Chair in Applied Numerical Analysis, for which I am once more indebted
to Professor Powell. Furthermore, these institutions and the
Department of Applied Mathematics and Theoretical Physics, have
enabled me to attend conferences and enjoy the opportunity to work
with colleagues abroad. I would also like to thank David Broomhead,
Alfred Cavaretta, Nira Dyn, David Levin, Charles Micchelli, John
Scales and Joe Ward, who have invited and sponsored my visits, and
have invariably provided hospitality and kindness.

There are many unmentioned people to whom I owe thanks. Certainly this
work would not have been possible without the help of my friends and
family. In particular, I thank my partner, Glennis Starling, and my
father. I am unable to thank my mother, who died during the last weeks
of this work, and I have felt this loss keenly. This dissertation is
dedicated to the memories of both my mother and my grandfather,
Charles S.~Wilkins. 

\vfill \eject
\parindent=30pt
\beginsection Table of Contents

\bigskip
\noindent Chapter 1 : {\bf Introduction}
\item{1.1} Polynomial interpolation \dotfill page 3 
\item{1.2} Tensor product methods \dotfill page 4
\item{1.3} Multivariate splines \dotfill page 4
\item{1.4} Finite element methods \dotfill page 5
\item{1.5} Radial basis functions \dotfill page 6
\item{1.6} Contents of the thesis \dotfill page 10
\item{1.7} Notation \dotfill page 12
\medskip

\bigskip
\noindent Chapter 2 : {\bf Conditionally positive definite functions
and $p$-norm distance matrices}
\item{2.1} Introduction \dotfill page 14 
\item{2.2} Almost negative matrices \dotfill page 15
\item{2.3} Applications \dotfill page 18
\item{2.4} The case $p > 2$ \dotfill page 25
\medskip

\bigskip
\noindent Chapter 3 : {\bf Norm estimates for distance matrices}
\item{3.1} Introduction \dotfill page 31
\item{3.2} The univariate case for the Euclidean norm \dotfill page 31
\item{3.3} The multivariate case for the Euclidean norm \dotfill page 35
\item{3.4} Fourier transforms and Bessel transforms \dotfill page 37
\item{3.5} The least upper bound for subsets of a grid \dotfill page
40
\medskip

\bigskip
\noindent Chapter 4 : {\bf Norm estimates for Toeplitz distance
matrices I}
\item{4.1} Introduction \dotfill page 42
\item{4.2} Toeplitz forms and Theta functions \dotfill page 44
\item{4.3} Conditionally negative definite functions of order 1 \dotfill page 50
\item{4.4} Applications \dotfill page 58
\item{4.5} A stability estimate \dotfill page 62
\item{4.6} Scaling the infinite grid \dotfill page 65
\item{}    Appendix \dotfill page 68
\medskip

\bigskip
\noindent Chapter 5 : {\bf Norm estimates for Toeplitz distance
matrices II}
\item{5.1} Introduction \dotfill page 70
\item{5.2} Preliminary facts \dotfill page 71
\item{5.3} P\'olya frequency functions \dotfill page 80
\item{5.4} Lower bounds on eigenvalues \dotfill page 88
\item{5.5} Total positivity and the Gaussian cardinal function \dotfill page 92
\medskip

\bigskip
\noindent Chapter 6 : {\bf Norm estimates and preconditioned conjugate gradients}
\item{6.1} Introduction \dotfill page 95
\item{6.2} The Gaussian \dotfill page 96
\item{6.3} The multiquadric \dotfill page 102
\medskip

\bigskip
\noindent Chapter 7 : {\bf On the asymptotic cardinal function for the
multiquadric}
\item{7.1} Introduction \dotfill page 118
\item{7.2} Some properties of the multiquadric \dotfill page 119
\item{7.3} Multiquadrics and entire functions of exponential type $\pi$ \dotfill page 121
\item{7.4} Discussion \dotfill page 125
\medskip

\bigskip
\noindent Chapter 8 : {\bf Conclusions} \dotfill page 126
\medskip

\bigskip
\noindent References \dotfill page 128

\vfill \eject
\pageno=1

\def\chapno{1.}
\headline{\ifnum\pageno=1\hfil\else{\headfont\hfil\chptitle\hfil}\fi}
\def\chptitle{Introduction}


\bigskip
\centerline{\bf 1 : Introduction}
\medskip

\noindent
The multivariate interpolation problem occurs frequently
in many branches of science and engineering. Typically, we are given a
discrete set $I$ in $\Rd$, where $d$ is greater than one, and real  numbers
$\{f_i\}_{i\in I}$. Our task is to construct a continuous
or sufficiently differentiable function $s \colon \Rd \to
\R1$ such that
$$ s(i) = f_i, \qquad i \in I, \eqno{(1.1)}$$
and we say that $s$ interpolates the data $\{(i, f_i): i \in
I \}$. Interpolants can be highly useful. For example, we
may need to approximate a function whose values are 
known only at the interpolation points, that is we are
ignorant of its behaviour outside $I$.
Alternatively, the underlying function might be far
too expensive to evaluate at a large number of points, in which case the aim
is to choose an interpolant which is cheap to compute. We
can then use our interpolant in other algorithms in order
to, for example, calculate approximations to extremal values of
the original function. Another application is
data-compression, where the size of our initial data
$\{(i,f_i): i \in {\hat I}\}$ exceeds the storage capacity of 
available computer hardware. In this case, we can choose
a subset $I$ of ${\hat I}$ and use the corresponding data to construct
an interpolant with which we estimate the remaining
values. It is important to note that in general $I$ will consist of scattered
points, that is its elements can be irregularly positioned. Thus
algorithms that apply to arbitrary distributions of points are necessary.
Such algorithms exist and are well understood in the univariate case (see, for
instance, Powell (1981)), but many difficulties intrude when $d$ is bigger than
one.

There are many applications of multivariate interpolation, but we prefer to treat a
particular application in some detail rather than provide a list. 
Therefore we consider the following
interesting example of Barrodale {\it et al} (1991).

When a time-dependent system is under observation, it is often necessary to
relate pictures of the system taken at different times. For example, when
measuring the growth of a tumour in a patient, we must expect many changes to
occur between successive X-ray photographs,
such as the position of the patient or the amount of fluid in the body's
tissues. If we can identify corresponding points on the two photographs, such as
parts of the bone structure or intersections of particular veins, then these
pairs of points can be viewed as the data for two interpolation problems.
Specifically, let $(x_j,y_j)_{j=1}^n$ be the coordinates of the points in one
picture, and let the corresponding points in the second picture be $(\xi_j,
\eta_j)_{j=1}^n$. We need functions $s_x \colon \R1^2 \to \R1$ and $s_y \colon
\R1^2 \to \R1$ such that 
$$ s_x(x_j,y_j) = \xi_j \hbox{ and } s_y(x_j,y_j) =
\eta_j \hbox{ for } j = 1, \ldots, n. \eqno{(1.2)}$$ 
Therefore we see that the scattered data
interpolation problem arises quite naturally as an attempt to approximate the
non-linear coordinate transformation mapping one picture into the next.

It is important to understand that interpolation is not always desirable.
For example, our data may be corrupted by measurement errors, in which case there is
no good reason to choose an approximation which satisfies the interpolation
equations, but we do want to construct an approximation which is close to the
function values in some sense. One option is to choose
our function $s \colon \Rd \to \R1$ from some family (usually a linear space) of
functions so as to minimize a certain functional $G$, such as 
$$ G(s-f) =
\sum_{i\in I} [f_i - s(i)]^2, \eqno{(1.3)}$$ 
which is the familiar least-squares
fitting problem. Of course this can require the solution of a nonlinearly 
constrained optimization problem, depending on the family of functions and the
functional $G$.
Another alternative to interpolation takes $s$ to be the sum 
of decaying functions, each centred at a point in $I$ and taking the
function value at that point. Such an approximation is usually called a
{\it quasi-interpolant}, reflecting the requirement that it should resemble the
interpolant in some suitable way. These methods are of both practical and
theoretical importance, but we emphasize that this dissertation is
restricted to interpolation, specifically interpolation using radial basis
functions, for which we refer the reader to Section 1.5 and the later
chapters of the dissertation.

We now briefly describe some other multivariate approximation schemes. Of course,
our treatment does not provide a thorough overview of the field, for which we
refer the reader to de Boor (1987), Franke
(1987) or Hayes (1987).
However, it is interesting to contrast radial basis functions with some of the
other methods. In fact, the memoir of Franke (1982) is dedicated to this
purpose; it contains careful numerical experiments using some thirty methods, including
radial basis functions, and provides an excellent reason for their theoretical
study: they obtain excellent  accuracy when interpolating scattered data.
Indeed, Franke found them to excel in this sense when compared to the other
tested methods, thus providing an excellent reason for their
theoretical study.

\sect {1.1 Polynomial interpolation}
Let $P$ be a linear space of polynomials in $d$ real variables spanned by
$(p_i)_{i\in I}$, where $I$ is the discrete subset of $\Rd$ discussed at the
beginning of the introduction. Then an interpolant $s \colon \Rd \to \R1$ of
the form 
$$ s(x) = \sum_{i\in I} c_i p_i(x), \qquad x \in \Rd, \eqno{(1.4)}$$
exists if and only if the matrix $(p_i(j))_{i,j \in I}$ is invertible. We see
that this property depends on the geometry of the centres when $d > 1$, which
is a significant difficulty. One solution is to choose a particular geometry. As
an example we describe the tensor product approach on a ``tartan grid''.
Specifically, let $I = \{ (x_j,y_k): 1 \le j \le l, 1 \le k \le m \}$, where
$x_1 < \cdots < x_l$ and $y_1 < \cdots < y_m$ are given real numbers, and let
$\{f_{(x_j,y_k)}: 1 \le j \le l, 1\le k \le m\}$ be the function values at
these centres. We let $(L^1_j)_{j=1}^l$ and $(L^2_k)_{k=1}^m$ be the usual
univariate Lagrange interpolating polynomials associated with the numbers
$(x_j)_1^l$ and $(y_k)_1^m$ respectively and define our interpolant $s \colon
\R1^2 \to \R1$ by the equation
$$ s(x,y) = \sum_{j=1}^l \sum_{k=1}^m f_{(x_j,y_k)} L^1_j(x) L^2_k(y), \qquad
(x,y) \in \R1^2. \eqno{(1.5)}$$
Clearly this approach extends to any number of dimensions $d$.

\sect {1.2 Tensor product methods}
The tensor product scheme for tartan grids described in the previous section is
not restricted to polynomials. Using the same notation as before, we replace
$(L_j^1)_{j=1}^l$ and $(L_k^2)_{k=1}^m$ by univariate functions $(P_j)_{j=1}^l$
and $(Q_k)_{k=1}^m$ respectively. Our interpolant takes the form
$$ s(x,y) = \sum_{j=1}^l \sum_{k=1}^m y_{jk} P_j(x) Q_k(y),
\qquad (x,y) \in \R1^2, \eqno{(1.7)}$$
from which we obtain the coefficients $(y_{jk})$. By adding points
outside the interval $[x_1, x_l]$ and $[y_1, y_m]$ we can choose
$(P_j)$ and $(Q_k)$ to be
univariate B-splines. In this case the linear systems involved are
invertible and banded, so that the number of operations and the
storage required are both multiples of the total number of points in
the tartan grid. Such methods are extremely
important for the subtabulation of functions on regular grids, and clearly the
scheme exists for any number of dimensions $d$. A useful survey is the book of
Light and Cheney (1986)

\sect {1.3 Multivariate Splines}
Generalizing some of the properties of univariate splines to a multivariate
setting has been an {\it id\'ee fixe} of approximation theory. Thus the name 
``spline'' is overused, being applied to almost any extension of univariate
spline theory. In this section we briefly consider box splines.
These are compactly supported piecewise polynomial functions which extend
Schoenberg's characterization of the $B$-spline $B(\cdot; t_0, \ldots, t_k)$
with arbitrary knots $t_0, \ldots, t_k$ as the ``shadow'' of a $k$-dimensional
simplex (Schoenberg (1973), Theorem 1, Lecture 1). Specifically, the box spline $B(\cdot;A)$ associated with the $d \times
n$ matrix $A$ is the distibution defined by  $$ B(\cdot; A) : C_0^\infty(\Rd)
\to \R1 : \phi \mapsto \int_{[-1/2,1/2]^n} \phi(Ax)\,dx, $$ where
$C_0^\infty(\Rd)$ is the vector subspace of $C^\infty(\Rd)$ whose elements vanish
at infinity. If we let $a_1, \ldots, a_n \in \Rd$ be the columns of $A$, then
the Fourier transform of the box spline is given by $$ {\hat B}(\xi; A) =
\prod_{j=1}^n \hbox{sinc}\ \xi^T a_j, \qquad \xi \in \Rd, $$ 
where $\hbox{sinc}(x) = \sin(x/2)/(x/2)$. We see that a simple example of a
box spline is a tensor product of univariate B-splines. It can be shown that
there exist box splines with smaller supports than tensor product B-splines.

A large body of mathematics now exists, and a suitable comprehensive
review is the long paper of Dahmen and Micchelli (1983). Further, this
theory is also yielding useful results in the study of wavelets (see
Chui (1992)). However, there are many computational difficulties. At present,
box spline software is not available from the main providers of scientific
computation packages.

\sect {1.4 Finite element methods}
Finite element methods can provide extremely flexible
piecewise polynomial spaces for approximation and scattered data interpolation.
When $d=2$ we first choose a {\it triangulation} of the points. Then a
polynomial is constructed on each triangle, possibly using function values and
partial derivative values at other points in addition to the vertices of the
triangulation. This is a non-trivial problem, since we usually require some
global differentiability properties, that is the polynomials must fit together
in a suitably smooth way. Further, the partial derivatives are frequently
unknown, and these methods can be highly sensitive to the accuracy of
their estimates (Franke (1982)). 

Much recent research has been directed
towards the choice of triangulation. The Delaunay triangulation (Lawson (1977))
is often recommended, but some work of Dyn, Levin and Rippa (1986) indicates that
greater accuracy can be achieved  using {\it data-dependent triangulations},
that is triangulations whose component triangles reflect the geometry of the
function in some way. Finally, the complexity of constructing triangulations in
higher dimensions effectively limits these methods to two and three dimensional
problems.

\sect {1.5 Radial basis functions}
A radial basis function approximation takes the form
$$ s(x) = \sum_{i \in I} y_i \phi(\|x-i\|), \qquad x \in \Rd, \eqno{(1.8)}$$
where $\phi \colon [0,\infty) \to \R1$ is a fixed univariate function and the
coefficients $(y_i)_{i\in I}$ are real numbers. We do not place any
restriction on the norm $\| \cdot \|$ at this
point, although we note that the Euclidean norm is the most common
choice. Therefore our approximation $s$ is a linear combination of
translates of a fixed function $x \mapsto \phi(\|x\|)$ which is ``radially
symmetric'' with respect to the given norm, in the sense that it clearly
possesses the symmetries of the unit ball. We shall often say that
the points $(x_j)_{j=1}^n$ are the {\it centres} of the radial basis
function interpolant. Moreover, it is usual to refer to $\phi$ as the
radial basis function, if the norm is understood.

If $I$ is a finite set, say $I = (x_j)_{j=1}^n$, the interpolation
conditions provide the linear system 
$$ A y = f, \eqno{(1.9)}$$
where
$$ A = \Bigl( \phi(\|x_j - x_k\|) \Bigr)_{j,k=1}^n, \eqno{(1.10)}$$
$y = (y_j)_{j=1}^n$ and $f = (f_j)_{j=1}^n$. 

One of the most attractive features of radial basis function methods
is the fact that a unique interpolant is often guaranteed under
rather mild conditions on the centres. In several important cases, the only
restrictions are that there are at least two centres and they are all
distinct, which are as simple as one could wish. However, one
important exception to this statement is the {\it thin plate spline}
introduced by Duchon (1975, 1976), where we choose $\phi(r) = r^2 \log
r$. It is easy to see that the interpolation matrix $A$ given by
(1.10) can be singular for non-trivial sets of distinct centres. For
example, choosing $x_2, \ldots, x_n$ to be any different points on
the sphere of unit radius whose centre is $x_1$, we conclude that the first
row and column of $A$ consist entirely of zeros. Of course,
such examples exist for any function $\phi$ with more than one zero.
Fortunately, it can be shown that it is suitable to add a polynomial of
degree $m \ge 1$ to the definition of $s$ if the centres are {\it
unisolvent}, which means that the zero polynomial is the only polynomial of
degree $m$ which vanishes at every centre (see, for instance, Powell
(1992)). 
The extra degrees of freedom are usually taken up by moment conditions on the coefficients
$(y_j)_{j=1}^n$. Specifically, we have the equations
$$\eqalign{
 \sum_{k=1}^n y_k \phi(\|x_j - x_k\|) + P(x_j) &= f_j, \qquad j = 1, 2,
\ldots, n, \cr
 \sum_{k=1}^n y_k p(x_k) &= 0 \qquad \hbox{ for every } p \in \Pi_m(\Rd), 
}\eqno{(1.11)}$$
where $\Pi_m(\Rd)$ denotes the vector space of polynomials in $d$ real
variables of total degree $m$, and the theory guarantees the existence
of a unique vector $(y_j)_{j=1}^n$ and a unique polynomial $P \in
\Pi_m(\Rd)$ satisfying (1.11). Moreover, because (1.8) does not
reproduce polynomials when $I$ is a finite set, it is sometimes useful
to augment $s$ in this way. 

In fact Duchon derived (1.11) as the solution to a
variational problem when $d = 2$: he proved that the function $s$ given by (1.11)
minimizes the integral
$$ \int_{\R1^2} [s_{x_1 x_1}]^2 + 2 [s_{x_1 x_2}]^2 + [s_{x_2 x_2}]^2 \,
dx, $$
where $m=1$ and $s$ satisfies some differentiability conditions.
Duchon's treatment is somewhat abstract, using sophisticated distribution
theory techniques, but a detailed alternative may be found in
Powell (1992). We do not study the thin plate spline in this
dissertation, although many of our results are highly relevant to its
behaviour.

In his comparison of multivariate approximation methods,
Franke (1982) considered several radial basis functions including the
thin plate spline. Therefore we briefly consider some of these functions.

\medskip
\noindent{\bf The multiquadric}

\noindent Here we choose $\phi(r) = (r^2 + c^2)^{1/2}$, where $c$ is a
real constant. The interpolation matrix $A$ is invertible provided
only that the points are all different and there are at least two of
them. Further, this matrix has an important spectral property: it is
{\it almost negative definite}; we refer the reader to Section 2 for
details.

Franke found that this radial basis function provided the most
accurate interpolation surfaces of all the methods tried for
interpolation in two dimensions. His centres were mildly irregular in
the sense that the range of distances between centres was not so large
that the average distance became useless. He found that the method
worked best when $c$ was chosen to be close to this average distance.
It is still true to say that we do not know how to choose $c$ for a
general function. Buhmann and Dyn (1991) derived error estimates which
indicated that a large value of $c$ should provide excellent accuracy.
This was borne out by some calculations and an analysis of Powell
(1991) in the case when the centres formed a regular grid in
one dimension. Specifically, he found that the uniform norm of the
error in interpolating $f(x) = x^2$ on the integer grid decreased by a
factor of $10^3$ when $c$ increased by one; see Table 6 of Powell
(1991) for these stunning results. In Chapter 7 of this thesis we are
able to show that the interpolants converge uniformly as $c \to
\infty$ if the underlying function is square-integrable and {\it band-limited}, that is its
Fourier transform is supported by the interval $[-\pi,\pi]^d$. Thus,
for many functions, it would seem to be useful to choose a large value of $c$.
Unfortunately, if the centres form a finite regular grid, then
we find that the smallest eigenvalue of the distance decreases
exponentially to zero as $c$ tends to infinity. Indeed, the reader is
encouraged to consider Table 4.1, where we find that the smallest
eigenvalue decreases by a factor of about 20 when $c$ is increased by
one and the spacing of the regular grid is unity.

We do not consider the polynomial reproduction properties of the
multiquadric discovered by Buhmann (1990) in this dissertation, but we
do make use of some of his work, in particular his formula for the
cardinal function's Fourier transform (see Chapter 7). However, we cannot
resist mentioning one of the brilliant results of Buhmann, in particular the
beautiful and surprising result that the degree of polynomials reproduced by
interpolation on an infinite regular grid actually increases with the
dimension. The work of Jackson (1988) is also highly relevant here.

\medskip
\noindent{\bf The Gaussian}

\noindent
There are many reasons to advise users to avoid the Gaussian $\phi(r)
= \exp(-cr^2)$. Franke
(1982) found that it is very sensitive to the choice of parameter $c$,
as we might expect. Further, it cannot even reproduce constants when
interpolating function values given on an infinite regular grid (see
Buhmann (1990)). Thus its potential for practical computer
calculations seems to be small. However, it possesses many properties which
continue to win admirers in spite of these problems. In particular, it
seems that users are seduced by its smoothness and rapid decay.
Moreover the Gaussian interpolation matrix (1.10) is positive definite
if the centres are distinct, as well as being suited to iterative
techniques. I suspect that this state of affairs will continue until
good software is made available for radial basis functions such as the
multiquadric. Therefore I wish to emphasize that this thesis addresses
some properties of the Gaussian because of its theoretical importance
rather than for any use in applications. 

In a sense it is true to say that the Gaussian generates all of the
radial basis functions considered in this thesis. Here we are thinking
of the Schoenberg characterization theorems for conditionally negative
definite functions of order zero and order one. These theorems and
related results occur many times in this dissertation. 

\medskip
\noindent{\bf The inverse multiquadric}

\noindent
Here we choose $\phi(r) = (r^2 + c^2)^{-1/2}$. Again , Franke (1982) found that
this radial basis function can provide excellent approximations, even when the
number of centres is small. As for the multiquadric, there is no good choice of
$c$ known at present. However, the work presented in Chapter 7 does extend to
this function (although this analysis is not presented here), so that sometimes
a large value of $c$ can be useful.

\medskip
\noindent{\bf The thin plate spline}

\noindent
We have hardly touched on this highly important function, even though the works
of Franke (1982) and Buhmann (1990) indicate its importance is two dimensions
(and, more generally, in even dimensional spaces). 
However, we aim to generalize the
norm estimate material of Chapters 3--5 to this function in future. There is no
numerical evidence to indicate that this ambition is unfounded, and the
preconditioning technique of Chapter 6 works equally well when applied to this
function. Therefore we are optimistic that these properties will be understood
more thoroughly in the near future.

\sect{1.6 Contents of the thesis}
Like Gaul, this thesis falls naturally into three parts, namely Chapter 2,
Chapters 3--6, and Chapter 7. In Chapter 2 we study and extend the work of
Schoenberg and Micchelli on the nonsingularity of interpolation
matrices. One of
our main discoveries is that it is sometimes possible to prove nonsingularity
when the norm is non-Euclidean. Specifically, we prove that the interpolation
matrix is non-singular if we choose a $p$-norm for $1 < p < 2$ and if the
centres are different and there are at least two of them. This complements the
work of Dyn, Light and Cheney (1991) which investigates the case when $p = 1$.
They find that a necessary and sufficient condition for nonsingularity
when $d = 2$ 
is that the points should not form the vertices of a {\it closed path},
which is a closed polygonal curve consisting of alternately horizontal and
vertical arcs. For example, the $1$-norm interpolation matrix generated by the
vertices of any  rectangle is singular. Therefore it may be useful that we can
avoid these difficulties by using a $p$-norm for some $p \in (1,2)$. However,
the situation is rather different when $p > 2$. This is probably the most
original contribution of this section, since it makes use of a device that seems
to have no precursor in the literature and is wholly independent of the
Schoenberg-Micchelli corpus. We find that, if both $p$ and the dimension $d$
exceed two, then it is possible to construct sets of distinct points which
generate a singular interpolation matrix. It is interesting to relate that these
sets were suggested by numerical experiment, and the author is grateful to
M.~J.~D.~Powell for the use of his TOLMIN optimization software.

The second part of this dissertation is dedicated to the study of the spectra of
interpolation matrices. Thus, having studied the nonsingularity (or otherwise)
of certain interpolation matrices, we begin to quantify .
This study was initiated by the beautiful papers of Ball (1989), and Narcowich
and Ward (1990, 1991), which provided some spectral bounds for several functions,
including the multiquadric. Our main findings are that it is possible to use
Fourier transform methods to address these questions, and that, if the centres
form a subset of a regular grid, then it is possible to provide a sharp upper
bound on the norm of the inverse of the interpolation matrix. Further,
we are able to understand the distribution of {\it all} the
eigenvalues using some work of Grenander and Szeg\H{o} (1984). This work
comprises Chapters 3 and 4. In the latter section, it turns out that everything
depends on an infinite product expansion for a Theta function of Jacobi type. This
connection with classical complex analysis still excites the author, and this
excitement was shared by Charles Micchelli. Our collaboration, which forms
Chapter 5, explores a property of P\'olya frequency functions which generalizes
the product formula mentioned above. Furthermore, Chapter 5 contains several
results which attack the norm estimate problem of Chapter 4 using a
slightly different
approach. We find that we can remove some of the assumptions required
at the expense of a little 
more abstraction. This work is still in progress, and we cannot yet say
anything about the approximation properties of our suggested class of
functions. We have included this work because we think it is interesting and,
perhaps more importantly, new mathematics is frequently open-ended.

Chapters 6 and 7 apply the work of previous chapters. In Chapter 6 we use our
study of Toeplitz forms in Chapter 4 to suggest a preconditioner for the
conjugate gradient solution of the interpolation equations, and the results are
excellent, although they only apply to finite regular grids. Of course it is
our hope to extend this work to arbitrary point sets in future. We remark that
our approach is rather different from the variational heuristic of 
Dyn, Levin and Rippa (1986), which concentrated  on preconditioners for thin
plate splines in two dimensions. Probably our most important practical finding
is that the number of iterations required to attain a solution to within a
particular tolerance seems to be independent of the number of centres.

Next, Chapter 7 is unique in that it is the
only chapter of this thesis which concerns itself with the approximation power
of radial basis function spaces. Specifically, we investigate the behaviour of
interpolation on an infinite regular grid using a multiquadric $\phi(r) = (r^2
+ c^2)^{1/2}$ when the parameter $c$ tends to infinity. We find an interesting
connection with the classical theory of the Whittaker cardinal spline: the
Fourier transform of the cardinal (or fundamental) function of interpolation
converges (in the $L^2$ norm) to the characteristic function of the cube
$[-\pi,\pi]^d$. This enables us to show that the interpolants to certain
band-limited functions converge uniformly to the underlying function when $c$
tends to infinity.

\medskip
\noindent{\bf An aside}
\noindent
\def\Sdash{S^\prime}
Finally, we cannot resist the following excursion into the theory of
conic sections, whose only purpose
is to lure the casual reader. Let $S$ and $\Sdash$ be different points
in $\R1^2$ and let $f \colon \R1^2 \to \R1$ be the function defined by
$$ f(x) = \|x-S\| + \|x-S^\prime\|, \qquad x \in \R1^2, $$
where $\|\cdot\|$ is the Euclidean norm. Thus the contours of $f$
constitute the set of all ellipses whose focal points are $S$ and
$\Sdash$. By direct calculation we obtain the expression
$$ \nabla f(x) = \Bigl({x-S \over \|x-S\|}\Bigr) 
+ \Bigl({x-\Sdash \over \|x - \Sdash\|}\Bigr)$$
which implies the relations
$$ \Bigl({x-S \over \|x-S\|}\Bigr)^T \nabla f(x) 
 = 1 + \Bigl({x-S \over \|x-S\|}\Bigr)^T \Bigl({x-\Sdash \over
\|x-\Sdash\|}\Bigr)
 = \Bigl({x-\Sdash \over \|x-\Sdash\|}\Bigr)^T \nabla f(x), $$
whose geometric interpretation is the reflector property of the
ellipse. A similar derivation exists for the hyperbola.

\sect {1.7 Notation}
We have tried to use standard notation throughout this thesis with a few
exceptions. Usually we denote a finite sequence of points in
$d$-dimensional real space $\Rd$ by  subscripted
variables, for example $(x_j)_{j=1}^n$. However we have avoided this usage when
coordinates of points occur. Thus Chapters 2 and 5 use superscripted variables,
such as $(x^j)_{j=1}^n$, and coordinates are then indicated by
subscripts. For example, $x^j_k$ denotes the $k$th coordinate of the
$j$th vector of a sequence of vectors $(x^j)_{j=1}^n$.     
The inner product of two vectors $x$ and $y$ is
denoted $xy$ in the context of a Fourier transform, but we have used the more
traditional linear algebra form $x^T y$ in Chapter 6 and in a few
other places. We have used no special notation for vectors, and we
hope that no ambiguity arises thereby.

Given any absolutely integrable function $f
\colon \Rd \to \R1$, we define its Fourier transform by the equation
$$ \fhat(\xi) = \int_\Rd f(x) \exp(-ix\xi)\,dx, \qquad \xi \in \Rd.$$
We also use this normalization when discussing distributional Fourier
transforms.  Thus, if it is permissible to invert the Fourier
transform, then the integral takes the form
$$ f(x) = (2\pi)^{-d} \int_\Rd \fhat(\xi) \exp(ix\xi)\,d\xi, \qquad
x\in \Rd. $$

The norm symbol ($\|\cdot\|$) will usually denote the Euclidean norm,
but this is not so in Chapter 1. Here the Euclidean norm is denoted by
$|\cdot|$ to distinguish it from other norm symbols. 

Finally, the reader will find that the term ``radial basis
function'' can often mean the univariate function $\phi \colon
[0,\infty) \to \R1$ {\it and} the multivariate function $\Rd \ni x
\mapsto \phi(\|x\|)$. This abuse of notation was inherited from the
literature and seems to have become quite standard. However, such
potential for ambiguity is bad. It is perhaps unusual for the author
of a dissertation to deride his own notation, but it is hoped that the
reader will not perpetuate this terminology.

\vfill \eject
\def\chapno{2.}
\headline{\ifnum\pageno=14\hfil\else{\headfont\hfil\chptitle\hfil}\fi}
\def\chptitle{Conditionally positive functions and $p$-norm distance
matrices}


\bigskip
\centerline{\bf 2 : Conditionally positive functions and}
\smallskip
\centerline{\bf $\bf p$-norm distance matrices}
\medskip

\sect {\chapno1. Introduction}
The real multivariate interpolation problem is as follows. Given distinct points
$x^1, \ldots , x^n \in\Rd$ and real scalars $f_1, \ldots ,
f_n$, we wish to construct a continuous function $s\colon\Rd
\to\R1$ for which 
$$s(x^i) = f_i, \quad\hbox{ for } i = 1,\ldots ,n .$$ 
The radial basis function approach is to choose a
function $\phi\colon[0, \infty ) \to [0, \infty )$ and a norm
$\|\cdot \| $ on $\Rd$ and then let $s$ take the form 
$$s(x) = \sum_{i = 1}^n \lambda_i \, \phi (\| x - x^i \|) .$$
Thus $s$ is chosen to be an element of the vector space spanned by the
functions $\xi \mapsto \phi (\| \xi - x^i \| )$, for $i = 1,
\ldots , n $.  The interpolation conditions then define a linear
system $ A \lambda = f$, where $A \in\Rnn$ is given
by 
$$ A_{ij} = \phi ( \| x^i - x^j \| ), \quad\hbox{ for } 1 \le i, j \le n,$$ 
and where $\lambda = ( \lambda_1, ..., \lambda_n )$ and
$f = (f_1, ...,f_n)$. In this thesis, a matrix such as $A$ will be
called a distance matrix.

Usually $\|\cdot\|$ is chosen to be the Euclidean norm,
and in this case Micchelli (1986) has shown the distance matrix generated
by distinct points to be invertible for several useful choices of
$\phi $.  In this chapter, we investigate the invertibility of the
distance matrix when $\|\cdot\| $ is a $p$-norm for $1 < p <
\infty$, $p
\ne 2$, and $\phi (t) = t$, the identity. We find that $p$-norms do
indeed provide invertible distance matrices given distinct points, for
$1 < p \le 2$. Of course, $p = 2$ is the Euclidean case mentioned
above and is not included here.  Now Dyn, Light and Cheney (1991) have
shown that the $1-$norm distance matrix may be singular on quite
innocuous sets of distinct points, so that it might be useful to
approximate $\|\cdot\|_1$ by $\|\cdot\|_p$ for some $p \in (1,
2]$.  This work comprises section \chapno3. The framework of the proof is
very much that of Micchelli (1986).

For every $p > 2$, we find that distance matrices can be singular
on certain sets of distinct points, which we construct. We find that
the higher the dimension of the underlying vector space for the points
$x^1, \ldots , x^n$, the smaller the least $p$ for which there exists
a singular $p$-norm.

\sect {\chapno2. Almost negative matrices}
Almost every matrix considered in this section will induce a
non-positive form on a certain hyperplane in $\Rn$.
Accordingly, we first define this ubiquitous subspace and fix
notation.

\proclaimdefn {Definition \chapno2.1}. For any positive integer $n$, let
$$ Z_n = \{ \ y \in\Rn : \sum_{i=1}^n y_i = 0 \ \}. $$

\noindent Thus $Z_n$ is a hyperplane in $\Rn $. We note that $Z_1 = \{ 0
\}$.

\proclaimdefn {Definition \chapno2.2}. We shall call $A \in\Rnn $ 
{\it almost negative definite} (AND) 
if $A$ is symmetric and $$ y^T Ay \le 0, \quad\hbox{ whenever } y \in
Z_n.$$ Furthermore, if this inequality is strict for all non-zero $y
\in Z_n$, then we shall call $A$ {\it strictly AND}.
\vskip 18pt

\proclaim {Proposition \chapno2.3}. Let $A \in\Rnn $ be
strictly AND with non-negative trace. Then $$ (-1)^{n-1} \det A > 0.$$

\pf~We remark that there are no strictly AND $1 \times 1$
matrices, and hence $n \ge 2$. Thus $A$ is a symmetric matrix inducing
a negative-definite form on a subspace of dimension $n-1 > 0$, so that
$A$ has at least $n-1$ negative eigenvalues. But trace $A \ge 0$, and
the remaining eigenvalue must therefore be positive.~\qed

\noindent Micchelli (1986) has shown that both $A_{ij} = |x^i-x^j|$ and $A_{ij} =
(1+|x^i-x^j|^2)^{1\over 2}$ are AND, where here and subsequently $|\cdot|$
denotes the Euclidean norm.  In fact, if the points $x^1, \ldots, x^n$
are distinct and $n \ge 2$, then these matrices are strictly AND.
Thus the Euclidean and multiquadric interpolation matrices generated
by distinct points satisfy the conditions for proposition \chapno2.3.

Much of the work of this chapter rests on the following characterization of
AND matrices with all diagonal entries zero. This theorem is stated
and used to good effect by Micchelli (1986), who omits much of the proof
and refers us to Schoenberg (1935). Because of its extensive use 
we include a proof for the convenience of the reader.  The
derivation follows the same lines as that of Schoenberg (1935).

\proclaim {Theorem \chapno2.4}. Let $A \in\Rnn $ have all diagonal entries zero. Then $A$ is
AND if and only if there exist $n$ vectors $ y^1, \ldots, y^n \in
\Rn $ for which $$ A_{ij} = |y^i - y^j|^2. $$

\pf~Suppose $ A_{ij} = |y^i - y^j|^2 $ for vectors $ y^1,
\ldots, y^n \in\Rn$. Then $A$ is symmetric and the
following calculation completes the proof that $A$ is AND. Given any $z
\in Z_n$, we have 
$$\eqalign { 
 z^T A z &= \sum_{i,j=1}^n z_i z_j |y^i -y^j|^2 \cr 
         &= \sum_{i,j=1}^n z_i z_j ( |y^i|^2 + |y^j|^2 -
2(y^i)^T(y^j) )\cr 
         &= -2 \sum_{i,j=1}^n z_i z_j (y^i)^T(y^j) \quad
\hbox{ since the coordinates of $z$ sum to zero, } \cr 
         &= -2 \Bigl|\sum_{i=1}^n z_i y^i \Bigr|^2 \le 0 .
}$$

\noindent This part of the proof is given in Micchelli (1986).  The converse
requires two lemmata.

\proclaim {Lemma \chapno2.5}. Let $B \in\R1^{k \times k} $ be a symmetric 
non-negative definite matrix. 
Then we can find $\xi^1,
\ldots, \xi^k \in\R1^k $ such that $$ B_{ij} = |\xi^i|^2 +
|\xi^j|^2 - |\xi^i - \xi^j|^2 . $$

\pf~Since $B$ is symmetric and non-negative definite, we have
$ B = P^T P$, for some $P\in\R1^{k \times k} $.
Let $p^1, \ldots, p^k$ be the columns of $P$. Thus $$ B_{ij} = (p^i)^T
(p^j).$$
Now $$ |p^i - p^j|^2 = |p^i|^2 + |p^j|^2 - 2(p^i)^T(p^j).$$
Hence $$ B_{ij} = {1 \over 2}( |p^i|^2 + |p^j|^2 - |p^i - p^j|^2 ).$$
All that remains is to define $ \xi^i = p^i / \surd 2$ , for $i = 1,
\ldots, k$.~\qed

\proclaim {Lemma \chapno2.6}. Let $A \in\Rnn $. Let $e^1, \ldots, e^n$ denote
the standard basis for $\Rn $, and define 
$$\eqalign { 
   f^i &= e^n - e^i ,\hbox{ for } i = 1, \ldots, n-1, \cr 
   f^n &= e^n.
}$$
Finally, let $F \in\Rnn$ be the matrix with columns
$ f^1, \ldots, f^n.$ Then 
$$\eqalign { 
  (-F^T AF)_{ij} &= A_{in} + A_{nj} - A_{ij} - A_{nn}, 
                       \quad\hbox{ for } 1 \le i, j \le n-1, \cr
  (-F^T AF)_{in} &= A_{in} - A_{nn}, \cr 
  (-F^T AF)_{ni} &= A_{ni} - A_{nn}, 
                       \quad\hbox{ for } 1 \le i \le n-1,\cr 
  (-F^T AF)_{nn} &= -A_{nn}.
}$$

\pf~We simply calculate $(-F^T AF)_{ij} \equiv -(f^i)^T A
(f^j)$.~\qed

We now return to the proof of Theorem \chapno2.4: Let $A \in\Rnn$ 
be AND with all diagonal entries zero.  Lemma
\chapno2.6 provides a convenient basis from which to view the action of $A$.
Indeed, if we set $ B = -F^T AF$, as in Lemma \chapno2.6, we see that the
principal submatrix of order $n-1$ is non-negative definite, since
$f^1, \ldots, f^{n-1}$ form a basis for $Z_n$. Now we appeal to Lemma
\chapno2.5, obtaining $ \xi^1,\ldots, \xi^{n-1} \in\R1^{n-1}$ such that
$$ B_{ij} = |\xi^i|^2 + |\xi^j|^2 - |\xi^i -
\xi^j|^2 \ , \hbox{ for } 1 \le i, j \le n-1, $$
while Lemma \chapno2.6 gives $$ B_{ij} = A_{in} + A_{jn} - A_{ij}.$$
Setting $i = j$ and recalling that $A_{ii} = 0$, we find 
$$ A_{in} = |\xi^i|^2 , \quad\hbox{ for} 1 \le i\le n-1$$ 
and thus we obtain 
$$A_{ij} = |\xi^i - \xi^j|^2, \quad\hbox{ for } 1 \le i, j \le n-1 .$$

Now define $\xi^n = 0$. Thus $ A_{ij} = |\xi^i - \xi^j|^2$, for $1 \le
i, j \le n$, where $ \xi^1,\ldots, \xi^n \in\R1^{n-1}$. We may
of course embed $\R1^{n-1}$ in $\Rn$. More formally, let
$\iota \colon\R1^{n-1} \hookrightarrow\Rn $ be the map
$\iota \colon(x_1, \ldots, x_{n-1}) \mapsto (x_1, \ldots, x_{n-1}, 0)$, and,
for $i = 1, \ldots, n$, define $y^i = \iota ( \xi^i )$. Thus $y^1,
\ldots, y^n
\in\Rn$ and $$ A_{ij} = |y^i - y^j|^2.$$
The proof is complete.~\qed

Of course, the fact that $ y^n = 0 $ by this
construction is of no import; we may take any translate of the $n$
vectors $y^1,\ldots, y^n$ if we wish.

\sect {\chapno3. Applications}
In this section we introduce a class of functions inducing AND
matrices and then use our characterization Theorem \chapno2.4 to prove a
simple, but rather useful, theorem on composition within this class.
We illustrate these ideas in examples \chapno3.3-\chapno3.5. The remainder of the
section then uses Theorems \chapno2.4 and \chapno3.2 to deduce results concerning
powers of the Euclidean norm. This enables us to derive the promised
$p$-norm result in Theorem \chapno3.11.

\proclaimdefn {Definition \chapno3.1}. We shall call $f\colon[0,\infty) \to [0,\infty) $ 
a {\it conditionally
negative definite function of {order 1} (CND1)} if, for any positive
integers $n$ and $d$, and for any points $ x^1, \ldots, x^n
\in\Rd$, the matrix $A \in\Rnn $
defined by $$ A_{ij} = f(|x^i - x^j|^2), \hbox{ for } 1 \le i, j \le
n, $$ is AND. Furthermore, we shall call $f$ {\it strictly CND1} if
the matrix $A$ is strictly AND whenever $n \ge 2$ and the points $x^1,
\ldots, x^n$ are distinct.

This terminology follows that of Micchelli (1986), Definition \chapno3.1 .  We
see that the matrix $A$ of the previous definition satisfies the
conditions of proposition \chapno2.3 if $f$ is strictly CND1, $n \ge 2$ and
the points $x^1, \ldots, x^n$ are distinct.

\proclaim{Theorem \chapno3.2}.
\item{(1)} Suppose that $f$ and $g$ are CND1 functions and that $f(0) =
0 $. Then $g \circ f$ is also a CND1 function. Indeed, if $g$ is
strictly CND1 and f vanishes only at $0$, then $g \circ f$ is strictly
CND1. 
\item{(2)} Let A be an AND matrix with all diagonal entries zero.  Let
$g$ be a CND1 function.  Then the matrix defined by $$ B_{ij} =
g(A_{ij}), \hbox{ for } 1 \le i, j \le n,$$ is AND. Moreover, if $n
\ge 2$ and no off-diagonal elements of $A$ vanish, then $B$ is
strictly AND whenever $g$ is strictly AN. 

\pf~
\item{(1)} The matrix $ A_{ij} = f(|x^i - x^j|^2)$ is an AND matrix with all
diagonal entries zero.  Hence, by Theorem \chapno2.4, we can find $n$ vectors
$ y^1, \ldots, y^n \in\Rn $ such that $$ f(|x^i - x^j|^2) =
|y^i - y^j|^2. $$
But g is a CND1 function, and so the matrix $B \in\Rnn$ 
defined by $$ B_{ij} = g(|y^i - y^j|^2) = g \circ f(|x^i -
x^j|^2),$$ is also an AND matrix. Thus $g \circ f$ is a CND1 function.
The condition that $f$ vanishes only at $0$ allows us to deduce that
$y^i \ne y^j$, whenever $i \ne j$. Thus $B$ is strictly AND if $g$ is
strictly CND1.
\item{(2)} We observe that $A$ satisfies the hypotheses of Theorem \chapno2.4. We
may therefore write $A_{ij} = |y^i - y^j|^2$, and thus $B$ is AND
because $g$ is CND1. Now, if $A_{ij} \ne 0$ if $i \ne j$, then the
vectors $y^1, ...  , y^n$ are distinct, so that $B$ is strictly AND if
$g$ is strictly CND1.~\qed

For the next two examples only, we shall need the following concepts.
Let us call a function $g\colon[0,\infty)
\to [0,\infty) $ {\it positive definite} if, for any positive integers
$n$ and $d$, and for any points $ x^1, \ldots, x^n \in\Rd$,
the matrix $A \in\Rnn $ defined by $$ A_{ij} =
g(|x^i - x^j|^2), \hbox{ for } 1 \le i, j \le n, $$ is non-negative
definite.  Furthermore, we shall call $g$ {\it strictly positive
definite} if the matrix $A$ is positive definite whenever the points $
x^1, \ldots, x^n $ are distinct.  We reiterate that these last two
definitions are needed only for examples \chapno3.3 and \chapno3.4.

\proclaimdefn{Example \chapno3.3}. A Euclidean distance matrix $A$ is AND, indeed
strictly so given distinct points. This was proved by Schoenberg (1938)
and rediscovered by Micchelli (1986). Schoenberg also proved the stronger
result that the matrix
$$ A_{ij} = |x^i - x^j|^\alpha, \hbox{ for } 1 \le i,j \le n,$$
is strictly AND given distinct points $x^1, \ldots, x^n \in\Rd$, $n \ge 2$
and $0 < \alpha < 2$. We shall derive this fact using Micchelli's methods in
Corollary \chapno3.7 below, but we shall use the result here to
illustrate Theorem \chapno3.2. 
We see that, by Theorem \chapno2.4, there exist $n$ vectors $ y^1, \ldots, y^n
\in\Rn $ such that $$ A_{ij} \equiv |x^i - x^j|^\alpha = |y^i -
y^j|^2. $$
The vectors $ y^1, \ldots, y^n$ must be distinct whenever the points $
x^1, \ldots, x^n \in\Rd$ are distinct, since $A_{ij} \ne 0$
whenever $i \ne j$.

Now let $g$ denote any strictly positive definite function. Define $B \in \Rnn$ 
by $$ B_{ij} \equiv g(A_{ij}). $$ Thus $$
g(|x^i - x^j|^\alpha) = g(|y^i - y^j|^2) .$$ Since we have shown that the
vectors $ y^1, \ldots, y^n$ are distinct, the matrix $B$ is therefore
positive definite.

For example, the function $g(t) = \exp(-t)$ is a strictly positive
definite function.  For an elementary proof of this fact, see Micchelli
(1986), p.15 .  Thus the
matrix whose elements are $$ B_{ij} = \exp(\ -|x^i - x^j|^\alpha), 1 \le i,j
\le n, $$ is always (i) non-negative definite, and (ii) positive
definite whenever the points $ x^1, \ldots, x^n $ are distinct
\qed

\proclaimdefn{Example \chapno3.4}. This will be our first example using a $p$-norm
with $p \ne 2$. Suppose we are given distinct points $ x^1, \ldots,
x^n
\in\Rd$. Let us define $A \in\Rnn
$ by $$ A_{ij} = \| x^i - x^j \|_1 . $$
Furthermore, for $k = 1, \ldots, d$, let $A^{(k)} \in\Rnn$ 
be given by $$ A_{ij}^{(k)} = | x^i_k - x^j_k |,$$
recalling that $ x^i_k$ denotes the $k^{th}$ coordinate of the point
$x^i$.

We now remark that $ A = \sum_{i=1}^d A^{(k)} $. But every $A^{(k)}$
is a Euclidean distance matrix, and so every $A^{(k)}$ is AND.
Consequently $A$, being the sum of AND matrices, is itself AND.  Now
$A$ has all diagonal entries zero. Thus, by Theorem \chapno2.4, we can
construct $n$ vectors $ y^1, \ldots, y^n \in\Rn $ such that $$
A_{ij} \equiv \| x^i - x^j \|_1 = |y^i - y^j|^2. $$
As in the preceding example, whenever the points $ x^1, \ldots, x^n $
are distinct, so too are the vectors $ y^1, \ldots, y^n$.

\noindent This does not mean that $A$ is non-singular.  Indeed, Dyn, Light and
Cheney (1991) observe that the 1-norm distance matrix is singular for the
distinct points $\{ (0,0), (1,0), (1,1), (0,1) \} $.
\vskip 12pt
\hskip 12pt Now let $g$ be any strictly positive definite function. Define $B \in \Rnn$
by $$ B_{ij} = g(A_{ij}) = g( \| x^i - x^j
\|_1 ) = g(|y^i - y^j|^2) .$$
Thus $B$ is positive definite.

For example, we see that the matrix $ B_{ij} = \exp(\ -\| x^i - x^j
\|_1 )$ is positive definite whenever the points $ x^1, \ldots, x^n $
are distinct.~\qed

\proclaimdefn{Example \chapno3.5}. As in the last example, let $ A_{ij} = \| x^i -
x^j \|_1 $, where $n \ge 2$ and the points $ x^1, \ldots, x^n $ are
distinct.  Now the function $ f(t) = (1+t)^{1 \over 2}$ is strictly
CND1 (\ Micchelli (1986)\ ).  This is the CND1 function generating the
multiquadric interpolation matrix.  We shall show the matrix $B \in
\Rnn $ defined by $$ B_{ij} = f(A_{ij}) = (1+\|
x^i - x^j \|_1)^{1 \over 2}$$ to be strictly AND.

Firstly, since the points $ x^1, \ldots, x^n $ are distinct, the previous
example shows that we may write $$ A_{ij} = \| x^i - x^j \|_1 =
|y^i - y^j|^2 ,$$ where the vectors $y^1, \ldots, y^n$ are distinct.
Thus, since $f$ is strictly CND1, we deduce from Definition \chapno3.1 that
$B$ is a strictly AND matrix.~\qed

We now return to the main theme of this chapter.
Recall that a function $f$ is completely monotonic provided
that
$$ (-1)^k f^{(k)}(x) \ge 0 , \hbox{ for every }k = 0, 1, 2, \ldots
\hbox{ and for } 0 < x < \infty .$$
We now require a theorem of Micchelli (1986), restated in our notation.

\proclaim{Theorem \chapno3.6}. Let $f\colon[0,\infty) \to [0,\infty)$ have a completely monotonic
derivative. Then $f$ is a CND1 function. Further, if $f^\prime$ is non-constant, then
$f$ is strictly CND1.

\def\micchelli{}
\pf~This is Theorem \micchelli2.3 of Micchelli (1986).~\qed

\proclaim{Corollary \chapno3.7}. The function $g(t) = t^\tau$ is strictly CND1 for every $\tau \in (0,1)$.

\pf~The conditions of the previous theorem are satisfied by $g$.~\qed

We see now that we may use this choice of $g$ in Theorem \chapno3.2, as in the
following corollary.

\proclaim{Corollary \chapno3.8}. For every $\tau \in (0,1)$ and for every positive integer $k \in [1, d]$, define
$A^{(k)} \in\Rnn $ by $$ A_{ij}^{(k)} = |x^i_k -
x^j_k|^{2\tau }, \hbox{ for } 1 \le i,j \le n.$$ Then every $A^{(k)}$
is AND.

\pf~For each $k$, the matrix $(|x^i_k-x^j_k|)_{i,j=1}^n $
is a Euclidean distance matrix. 
Using the function $g(t) = t^\tau$, we
now apply Theorem \chapno3.2 (2) to deduce that $A^{(k)} = 
g( |x^i-x^j|^2)$ is AND.~\qed

We shall still use the notation $\| . \|_p$ when $p
\in (0,1)$, although of course these functions are not norms .

\proclaim{Lemma \chapno3.9}. For every $p \in (0,2)$, the matrix $A
\in\Rnn $ defined by $$ A_{ij} = \| x^i -
x^j \|_p^p, \hbox{ for } 1 \le i, j \le n, $$ is AND.  If $n \ge 2
$ and the points $x^1, \ldots, x^n$ are distinct, then we can find
distinct $y^1, \ldots, y^n \in\Rn$ such that $$ \| x^i -
x^j \|_p^p = |y^i - y^j|^2 . $$ 

\pf~If we set $p = 2\tau $, then we see that $\tau \in (0,1)$
and $A = \sum_{k=1}^d A^{(k)}$, where the $A^{(k)}$ are those matrices
defined in Corollary \chapno3.8. Hence so that each $A^{(k)}$ is AND, and
hence so is their sum.  Thus, by Theorem \chapno2.4, we may write $$A_{ij} =
\| x^i - x^j
\|_p^p = |y^i - y^j|^2 . $$ Furthermore, if $n \ge 2 $ and the points
$x^1, \ldots, x^n$ are distinct, then $A_{ij} \ne 0$ whenever $i \ne
j$, so that the vectors $y^1, \ldots, y^n$ are distinct.~\qed

\proclaim {Corollary \chapno3.10}. For any $p \in (0,2)$ and for any $\sigma \in (0,1)$, define 
$B \in\Rnn $ by $$B_{ij} = ( \| x^i - x^j
\|_p^p )^\sigma .$$ Then $B$ is AND. As before, if $n \ge 2 $ and
the points $x^1, \ldots, x^n$ are distinct, then $B$ is strictly AND.

\pf~Let $A$ be the matrix of the previous lemma and 
let $g(t) = t^\tau$. We now apply Theorem \chapno3.2 (2)
\qed

\proclaim {Theorem \chapno3.11}. For every $p \ \in (1,2)$, the $p$-norm distance matrix $B \in\Rnn $,
that is: $$B_{ij} = \| x^i - x^j \|_p , \hbox{ for } 1 \le i, j
\le n,$$ is AND. Moreover, it is strictly AND if $n
\ge 2 $ and the points $x^1, \ldots, x^n$ are distinct, in which case $$ (-1)^{n-1}\det B > 0 .$$

{\it Proof. } If $p \in (1,2)$, then $\sigma \equiv 1/p \ \in (0,1)$.
Thus we may apply Corollary \chapno3.12. The final inequality follows from
the statement of proposition \chapno2.3.~\qed
\vskip 12pt
\hskip 12pt We may also apply Theorem \chapno3.2 to the $p-$norm distance matrix, 
for $p \in (1,2]$, or indeed to the $p^{th}$ power
of the $p-$norm distance matrix, for $p \in (0,2)$. Of course, we do not have a norm
for $0<p<1$, but we define the function in the obvious way.
We need only note
that, in these cases, both classes satisfy the conditions of Theorem
\chapno3.2 (2). We now state this formally for the $p-$norm distance matrix

\proclaim {Corollary \chapno3.12}. Suppose the matrix $B$ is the $p-$norm 
distance matrix defined in Theorem \chapno3.13. Then,
if $g$ is a CND1 function, the matrix $g(B)$ defined by $$ g(B)_{ij}
= g(B_{ij}), \hbox{ for } 1 \le i,j \le n,$$ is AND. Further, if $n
\ge 2$ and the points $x^1, \ldots, x^n$ are distinct, then $g(B)$ is
strictly AND whenever $g$ is strictly AN.

\pf~This is immediate from Theorem \chapno3.11 and the statement of
Theorem \chapno3.2 (2).~\qed

\sect {\chapno4. The case $p > 2$}
We are unable to use the ideas developed in the previous section to
understand this case. However, numerical experiment suggested the
geometry described below, which proved surprisingly fruitful.  We
shall view $\R1^{m+n}$ as two orthogonal slices $\R1^m
\oplus\Rn$. Given any $p > 2$, we take the vertices $\Gamma_m$
of $[-m^{-1/p}, m^{-1/p}]^m \subset\R1^m$ and embed this in
$\R1^{m+n}$. Similarly, we take the vertices $\Gamma_n$ of
$[-n^{-1/p}, n^{-1/p}]^n \subset\Rn$ and embed this too in
$\R1^{m+n}$. We see that we have constructed two orthogonal cubes
lying in the $p$-norm unit sphere.

\proclaimdefn{Example}. If $m = 2$ and $n = 3$, then $\Gamma_m = \{ (\pm \alpha,
\pm \alpha , 0, 0, 0) \}$ and $\Gamma_n = \{ (0, 0, \pm \beta, \pm \beta , \pm \beta ) \}
$, where $\alpha = 2^{-1/p}$ and $\beta = 3^{-1/p}$.

Of course, given $m$ and $n$, we are interested in values of $p$ for which the 
$p-$norm distance matrix generated by $\Gamma_m \cup \Gamma_n$ is
singular. Thus we ask whether there exist scalars $\{\lambda_y\}_{\{y
\in \Gamma_m\}}$ and $\{\mu_z\}_{\{z \in \Gamma_n\}}$, not all zero,
such that the function $$ s(x) = \sum_{y \in \Gamma_m} \lambda_y \|
x-y \|_p + \sum_{z \in \Gamma_n} \mu_z \| x-z \|_p $$
vanishes at every interpolation point. In fact, we shall show that
there exist scalars $\lambda$ and $\mu$, not both zero, for which the
function $$ s(x) = \lambda \sum_{y \in \Gamma_m} \| x-y \|_p +
\mu \sum_{z \in \Gamma_n} \| x-z \|_p $$ vanishes at every
interpolation point.

We notice that
\item{(i)}For every $y \in \Gamma_m$ and $z \in \Gamma_n$,
we have $\| y - z \|_p = 2^{1/p}$.
\item{(ii)} The sum $\sum_{y \in \Gamma_m} \| \tilde y - y \|_p$
takes the same value for every vertex $\tilde y \in \Gamma_m$, and
similarly, {\it mutatis mutandis}, for $\Gamma_n$.

Thus our interpolation equations reduce to two in number: $$ \lambda
\sum_{y \in \Gamma_m} \| \tilde y - y \|_p
\ + \ 2^{n+1/p}\mu \ = \ 0,$$
and $$ 2^{m+1/p}\lambda \ + \ \mu \sum_{z \in \Gamma_n} \|
\tilde z - z \|_p \ = \ 0 , $$ where by (ii) above, we see that
$\tilde y $ and $\tilde z$ may be any vertices of $\Gamma_m, \Gamma_n$
respectively.

We now simplify the (1,1) and (2,2) elements of our reduced system by
use of the following lemma.

\proclaim{Lemma \chapno4.1}. Let $\Gamma $ denote the vertices of $[0, 1]^k
$.  Then  
$$ \sum_{x \in
\Gamma } \| x \|_p = \sum_{l=0}^k {k \choose l} l^{1/p} .$$

\pf~Every vertex of $\Gamma $ has coordinates taking the
values $0$ or $1$. Thus the distinct $p$-norms occur when exactly $l$
of the coordinates take the value $1$, for $l = 0, \ldots, k$; each of
these occurs with frequency $k \choose l$.~\qed

\proclaim{Corollary \chapno4.2}. $$ \sum_{y \in \Gamma_m} \| \tilde y - y
\|_p = 2 \sum_{k=0}^m {m \choose k} (k/m)^{1/p}, \hbox{ for every }
\tilde y \in \Gamma_m, \hbox{ and }$$
$$ \sum_{z \in \Gamma_n} \| \tilde z - z \|_p = 2 \sum_{l = 0}^n
{n \choose l} (l/n)^{1/p}, \hbox{ for every }
\tilde z \in \Gamma_n.$$ 

\pf~We simply scale the result of the previous lemma by
$2m^{-1/p}$ and $2n^{-1/p}$ respectively.~\qed

With this simplification, the matrix of our system
becomes $$ \left(\matrix{ 2\sum_{k=0}^m {m \choose k} (k/m)^{1/p} & \
&2^n.2^{1/p}\cr
\ & \ & \ \cr          
2^m.2^{1/p}& \ & 2\sum_{l=0}^n {n \choose l} (l/n)^{1/p}\cr}
\right).$$

We now recall that
$$ B_i(f_p, 1/2) = 2^{-i} \sum_{j=0}^i {i \choose j} (j/i)^{1/p}$$ is
the Bernstein polynomial approximation of order $i$ to the function
$f_p(t) = t^{1/p}$ at $t = 1/2$. Our reference for properties for
Bernstein polynomial approximation will be Davis (1975), sections 6.2 and
6.3. Hence, scaling the determinant of our matrix by $2^{-(m+n)}$, we
obtain the function $$\phi_{m,n}(p) = 4\/B_m(f_p, 1/2)\/B_n(f_p,
1/2) - 2^{2/p} .$$ We observe that our task reduces to investigation
of the zeros of $\phi_{m,n}$.

We first deal with the case $m = n$, noting the factorization:
$$ \phi_{n,n}(p) = \{ 2B_n(f_p,1/2) + 2^{1/p}\} \{ 2B_n(f_p,1/2) -
2^{1/p}\} .$$ Since $f_p(t) \ge 0$, for $t \ge 0$ we deduce from the
monotonicity of the Bernstein approximation operator that $B_n(f_p,
1/2) \ge 0$. Thus the zeros of $\phi_{n,n}$ are those of the factor
$$\psi_n(p) = 2B_n(f_p, 1/2) - 2^{1/p}.$$

\proclaim{Proposition \chapno4.3}. $\psi_n$ enjoys the following properties.
\item{(1)} $\psi_n(p) \to \psi (p)$, where $\psi (p) =
2^{1 - 1/p} - 2^{1/p}$, as $n \to \infty$.
\item{(2)} For every $p > 1$, $\psi_n (p) < \psi_{n+1}(p)$, for every
positive integer $n$.
\item{(3)} For each $n$, $\psi_n$ is strictly increasing for $p \in 
[1, \infty)$.
\item{(4)} For every positive integer $n$, $\lim_{p \to \infty} \psi_n
(p) = 1 - 2^{1-n}$. 

\pf
\item{(1)} This is a consequence of the convergence of Bernstein polynomial
approximation.
\item{(2)} It suffices to show that $B_n(f_p, 1/2) < B_{n+1}(f_p, 1/2)$, for
$p > 1$ and $n$ a positive integer. We shall use Davis (1975), Theorem
6.3.4: If $g$ is a convex function on $[0,1]$, then $B_n(g,x) \ge
B_{n+1}(g,x)$, for every $x \in [0,1]$. Further, if $g$ is non-linear
in each of the intervals $[{ {j-1}\over n}, {j \over n}]$, for $j = 1,
\ldots, n$, then the inequality is strict.
Every function $f_p$ is concave and non-linear on $[0,1]$
for $p > 1$, so that this inequality is strict and reversed.
\item{(3)} We recall that 
$$ \psi_n(p) = 2B_n(f_p, 1/2) - 2^{1/p} = 2^{1-n} \sum_{k=0}^n {n
\choose k} (k/n)^{1/p} - 2^{1/p}.$$ Now, for $p_2 > p_1 \ge 1$, we
note that $t^{1/p_2} > t^{1/p_1}$, for $t \in (0,1)$, and also that
$2^{1/p_2} < 2^{1/p_1}$. Thus $(k/n)^{1/p_2} > (k/n)^{1/p_1}$, for $k
= 1, \ldots, n-1$ and so $\psi_n(p_2) > \psi_n(p_1)$.
\item{(4)} We observe that, as $p \to \infty$,
$$ \psi_n(p) \to 2^{1-n} \sum_{k=1}^n {n \choose k} - 1 = 2( 1
- 2^{-n}) - 1 = 1 - 2^{1-n}.$$
\qed

\proclaim{Corollary \chapno4.4}. For every integer $n > 1$, each $\psi_n$ has
a unique root $p_n \in (2, \infty )$.  Further, $p_n \to 2$
strictly monotonically as $n \to \infty$.

\pf~We first note that $\psi (2) = 0$, and that this is the
only root of $\psi$.  By proposition \chapno4.3 (1) and (2), we see that
$$\lim_{n \to \infty}{\psi_n(2)} = \psi (2) = 0 \hbox{ and }
\psi_n(2) < \psi_{n+1}(2) < \psi (2) = 0.$$ By proposition \chapno4.3 (4), we
know that, for $n > 1$, $\psi_n$ is positive for all sufficiently
large $p$.  Since every $\psi_n$ is strictly increasing by proposition
\chapno4.3 (3), we deduce that each $\psi_n$ has a unique root $p_n \in (2,
\infty)$ and that $\psi_n(p) < (>) 0 $ for $p < (>) p_n$.

We now observe that $\psi_{n+1}(p_n) > \psi_n(p_n) = 0$, by proposition \chapno4.3 (2), whence
$2 < p_{n+1} < p_n$.  Thus $(p_n)$ is a monotonic decreasing sequence
bounded below by $2$.  Therefore it is convergent with limit in $[2,
\infty)$.  Let $p^*$ denote this limit.  To prove that $p^* = 2$, it
suffices to show that $\psi(p^*) = 0$, since $2$ is the unique root of
$\psi$.  
Now suppose that $\psi(p^*) \ne 0$. By continuity, $\psi$ is bounded away from zero
in some compact neighbourhood $N$ of $p^*$. We now recall the following theorem
of Dini: If we have a monotonic increasing sequence of continuous real-valued functions
on a compact metric space with continuous limit function, then the convergence is uniform.
A proof of this result may be found in many texts, for example Hille (1962), p. 78.
Thus $\psi_n \to \psi$ uniformly in $N$. Hence there is an integer $n_0$ 
such that $\psi_n$ is bounded away from zero for every $n \ge n_0$. But 
$p^* = \lim p_n$ and $\psi_n(p_n) = 0$ for each $n$, so that we have reached
a contradiction. Therefore $\psi(p^*) = 0$ as required.~\qed

Returning to our original scaled determinant $\phi_{n,n}$, we see that
$\Gamma_n \cup \Gamma_n$ generates a singular $p_n$-norm distance
matrix and $p_n \searrow 2$ as $n \to \infty$.  Furthermore $$
\phi_{m,m}(p) < \phi_{m,n}(p) < \phi_{n,n}(p) , \hbox{ for }
1 < m < n, $$ using the same method of proof as in proposition \chapno4.3
(2).  Thus $\phi_{m,n}$ has a unique root $p_{m,n}$ lying in the
interval $(p_n, p_m).$ We have therefore proved the following theorem.

\proclaim{Theorem \chapno4.5}. For any positive integers $m$ and $n$, both
greater than $1$, there is a $p_{m,n} > 2$ such that the $\Gamma_m
\cup
\Gamma_n$-generated $p_{m,n}$-norm distance matrix is singular.
Furthermore, if $1 < m < n$, then $$ p_m \equiv p_{m,m} > p_{m,n} >
p_{n,n} \equiv p_n, $$ and $p_n \searrow 2 \hbox{ as } n \to
\infty$. 

Finally, we deal with the ``gaps'' in the sequence $(p_n)$ as follows.
Given a positive integer $n$, we take the configuration $\Gamma_n \cup
\Gamma_n (\theta )$, where $\Gamma_n (\theta )$ denotes the vertices
of the scaled cube $[-\theta n^{-1/p}, \theta n^{-1/p}]^n$ and
$\theta > 0$. The $2 \times 2$ matrix deduced from corollary \chapno4.2 on
page 8 becomes $$ \left(\matrix{ 2 \sum_{k=0}^n {n \choose k}
(k/n)^{1/p} & \ &2^n (1+\theta^p)^{1/p}\cr
\ & \ & \ \cr
2^n (1+\theta^p)^{1/p}& \ & 2 \theta \sum_{k=0}^n {n \choose k}
(k/n)^{1/p}\cr}
\right).$$
Thus, instead of the function $\phi_{n,n}$ discussed above, we now
consider its analogue: $$ \phi_{n,n,\theta}(p) = 4\theta
B_n^2(f_p, 1/2) - (1 + \theta^p)^{2/p}. $$ If $p > p_n$, the unique
zero of our original function $\phi_{n,n}$, we see that
$\phi_{n,n,1}(p)
\equiv \phi_{n,n} (p) > 0 $, because every $\phi_{n,n}$ is strictly
increasing, by proposition \chapno4.3 (3). However, we notice that
$\lim_{\theta \to 0} \phi_{n,n,\theta}(p) = -1$, so that
$\phi_{n,n,\theta}(p) < 0$ for all sufficiently small $\theta
> 0$.  Thus there exists a $\theta^* > 0 \hbox{ such that }
\phi_{n,n,\theta^* }(p) = 0$. Since this is true for every $p >
p_n$, we have strengthened the previous theorem. We now state this
formally.

\proclaim{Theorem \chapno4.6}. For every $p > 2$, there is a configuration of
distinct points generating a singular $p$-norm distance matrix. 

It is interesting to investigate how rapidly the sequence of zeros
$(p_n)$ converges to $2$. We shall use Davis (1975), Theorem 6.3.6, which
states that, for any bounded function $f$ on $[0,1]$,
$$\lim_{n\to\infty}{n(B_n(f,x) - f(x)) } = {1 \over 2} x(1-x)f^{\prime
\prime}(x),
\hbox{ whenever }f^{\prime \prime }(x) \hbox{ exists}. $$Applying this to 
$$\psi_n(p) = 2 B_n(f_p, 1/2) - 2^{1/p}, $$ we shall derive the
following bound.

\proclaim{Proposition \chapno4.7}. $p_n = 2 + O(n^{-1}).$

\pf~We simply note that 
$$\eqalign { 
  0 &= \psi_n(p_n) \cr 
    &= \psi (p_n) + O(n^{-1}), \hbox{ by Davis (1975) 6.3.6,} \cr 
    &= \psi (2) + (p_n-2) \psi^\prime (2) + {\it o}(p_n-2) + O(n^{-1}) 
.}$$
Since $\psi^\prime (2) \ne 0$, we have $p_n - 2 = O(n^{-1})$.~\qed

\vfill \eject
\def\chapno{3.}
\headline{\ifnum\pageno=31\hfil\else{\headfont\hfil\chptitle\hfil}\fi}
\def\chptitle{Norm estimates for distance matrices}
\bigskip
\centerline{\bf 3 : Norm estimates for distance matrices}
\medskip

\def\pts{{(x_j)_{j=1}^n}}

\def\sjk{{\sum_{j,k=1}^n y_j y_k\,}}
\def\djk{{\| x_j - x_k \|}}
\def\xjk{{x_j-x_k}}
\def\yjk{{|\sum_{j=1}^n y_j e^{ix_j t}|^2}}
\def\trig{{|\sum_{j=1}^n y_j e^{ijt}|^2}}

\def\skzd{{\sum_{k\in\Zd}}}

\def\rint{{\int_{-\infty}^{\infty}}}
\def\hint{{\int_0^\infty}}

\def\Sdm1{{S^{d-1}}}

\def\ysqr{{\Vert y \Vert^2}}

\def\ainv{{A_n^{-1}}}
\def\aninv{{\Vert A_n^{-1} \Vert_2}}

\def\zsum{{\sum_{k=-\infty}^\infty}}

\def\hi{{(2\pi)^{-1}}}
\def\tint{{\hi\int_0^{2\pi}}}
\def\tintd{{(2\pi)^{-d}\int_{\Td}}}
\def\fhat{{\hat f}}
\def\ajxi{{\Bigl|\sum_{j\in\Zd} a_j \exp(ij\xi)\Bigl|^2}}
\def\Rdm0{{\Rd\setminus\{0\}}}
\def\ghat{{\hat g}}
\def\ytrigo{{\Bigl|\sum_{j\in\Zd} y_j \exp(ix_j \xi)\Bigl|^2}}
\def\ytrig2{{\Bigl|\sum_{j\in\Zd} y_j \exp(ij \xi)\Bigl|^2}}
\def\ynj{{y_j^{(n)}}}

\def\phichat{{\phihat_c}}
\def\chichat{{{\hat \chi}_c}}

\def\min{{\rm min}}

\def\ainv{{A_n^{-1}}}
\def\aninv{{\Vert A_n^{-1} \Vert_2}}
\def\djk{{\Vert x_j - x_k \Vert}}
\def\ghat{{\hat g}}
\def\half{{\textstyle 1\over2}}
\def\hint{{\int_0^\infty}}

\def\omd{{\omega_{d-1}}}
\def\Omd{{\Omega_d}}
\def\OO{{\cal O}}
\def\points{{(x_j)_1^n}}
\def\pts{{(x_j)_{j=1}^n}}
\def\phi{{\varphi}}
\def\phihat{{\hat\varphi}}
\def\Rd{{\R1^d}}
\def\Rdm0{{\Rd\setminus\{0\}}}
\def\rint{{\int_{-\infty}^{\infty}}}
\def\Rn{{\R1^n}}
\def\Rnn{{\R1^{n\times n}}}
\def\Sdm1{{S^{d-1}}}
\def\sjk{{\sum_{j,k=1}^n y_j y_k\,}}

\def\xjk{{x_j-x_k}}

\def\yjk{{|\sum_{j=1}^n y_j e^{ix_j t}|^2}}
\def\ysqr{{\Vert y \Vert^2}}

\def\skzd{{\sum_{k\in\Zd}}}
\def\Z0{{\sum_{j=1}^n y_j = 0}}
\def\Zd{{\ZZ^d}}

\sect {\chapno1. Introduction}
In this chapter we use Fourier transform techniques to derive inequalities of the form
$$ y^T\,Ay \le -\mu\,y^T y, \qquad y \in \Rn, \eqno{(\chapno1)}$$
where $\mu$ is a positive constant and $\sum_{j=1}^n y_j = 0$. Here we are
using the notation of the abstract. It can be shown that equation (\chapno1) implies
the bound $\|A^{-1}\|_2 \le 1/\mu$ (see Chapter 4). Such estimates
have been derived in Ball (1989), Narcowich and Ward (1990, 1991) and
Sun (1990), using a different technique. The author
submits that the derivation presented here for the Euclidean norm is more perspicuous. 
Further, we relate the generalized Fourier transform to the measure that occurs in
an important characterization theorem for those functions $\phi$ considered here.
This is useful because tables of generalized Fourier transforms are widely available,
thus avoiding several of the technical calculations of Narcowich
(1990, 1991). Finally, we mention
some recent work of the author that provides the {\sl least upper bound} on $\|A^{-1}\|_2$
when the points $(x_j)_{j\in\Zd}$ form a subset of $\Zd$.

The norm $\|\cdot \|$ will always be the Euclidean norm in this section. We shall denote
the inner product of two vectors $x$ and $y$ by $xy$.

\sect {\chapno2. The Univariate Case for the Euclidean Norm}
Let $n \ge 2$ and let $\points$ be points in $\R1$ satisfying the condition
$\djk \ge 1$ for $j \ne k$. We shall prove that
$$\Bigl| \sjk \djk \Bigr| \ge \half\ysqr, $$
whenever $\Z0$.

We shall use the fact that
the generalized Fourier transform of $\phi(x) = |x|$ is $\phihat(t)= -2/t^2 $ in the univariate case. 
A proof of this may be found in Jones (1982), Theorem 
7.32.

\proclaim {Proposition \chapno2.1}. If \ $\Z0$, then
$$\eqalign{
 \sjk\djk
    &= (2\pi)^{-1} \rint (-2/t^2) \sjk \exp(i(x_j-x_k)t)\,dt \cr
    &= -\pi^{-1} \rint  \yjk t^{-2}\,dt  .
} \eqno{(\chapno2)}$$

\pf The two expressions on the righthand side above are equal
because of the useful identity
$$ \sjk \exp(i(\xjk)t) = \yjk.$$ This identity will be used several times below. We now let
$${\hat g}(t) =  (-2t^{-2})\ \yjk , \ \ \hbox{ for } t \in \R1.$$ 
The condition $\Z0$ implies that $\hat g$ is uniformly bounded. Further, since
$\hat g(t) = \OO(t^{-2})$ for large $|t|$, we see that $\hat g$ is absolutely
integrable.
Thus we have the equation
$$ g(x) = (2\pi)^{-1} \rint \hat g(t) \exp(ixt)\,dt.$$ A standard result of the theory of 
generalized Fourier transforms (cf. Jones (1982), Theorem 7.14, pages 224ff) provides 
the expression
$$ \sjk \| x + \xjk \| = (2\pi)^{-1} \rint (-2t^{-2})\yjk \exp(ixt) \,dt,$$ where we have used the identity stated at the beginning of this proof. 
We now need only set $x = 0$ in this final equation.~\qed

\proclaim {Proposition \chapno2.2}. Let $B:\R1 \rightarrow \R1$ be a continuous
function such that $\hbox{ supp }(B)$ is contained
in the interval $ [-1,1]$ and $0 \le {\hat B}(t) \le t^{-2}$. 
If $n \ge 2$, $\djk \ge 1$ for $j\ne k$, and $\Z0$, then
$$ \sjk\djk \le -2B(0) \ysqr .$$ 

\pf By Proposition \chapno2.1 and properties of Fourier transforms, 
$$\eqalign{
   \sjk\djk &\le (2\pi)^{-1} \rint (-2{\hat B}(t)) \sjk \exp(i(x_j-x_k)t)dt \cr
            &= -2 \sjk B(x_j-x_k) \cr
            &= -2 B(0) \ysqr , 
}$$
where the first inequality follows from the condition ${\hat B}(t) \le t^{-2}$.
The last line is a consequence of $\hbox{supp} (B) \subset [-1,1]$.~\qed

\proclaim {Corollary \chapno2.3}. Let 
$$B(x) = \cases{ (1-|x|)/4, &if $|x| \le 1$\cr 0, &otherwise.\cr}  $$ Then $B$ satisfies the conditions of 
Proposition \chapno2.2 and $B(0) = 1/4$.

\pf  By direct calculation, we find that
$$ {\hat B}(t) = {{\sin^2 (t/2)}\over{t^2}} \le {1\over{t^2}}.  $$ It is clear that the other conditions of
Proposition \chapno2.2 are satisfied.~\qed

\noindent
We have therefore shown the following theorem to be true.

\proclaim {Theorem \chapno2.4}. Let $\points$ be points in $\R1$ such that $n \ge 2$ and $\djk \ge 1$
when $j \ne k$. If \ $\Z0$, then
$$ \sjk\djk \le -{1\over2} \ \ysqr .$$ 

\noindent We see that a consequence of this result is the non-singularity
of the Euclidean distance matrix when the points $\points$ are distinct
and $n \ge 2$. It is important to realise that the homogeneity of the
Euclidean norm allows us to replace the condition ``$\djk \ge 1$ if $j \ne
k$''
by ``$\djk \ge \epsilon$ if $j \ne k$''. We restate Theorem \chapno2.4 in this
form for the convenience of the reader:

\proclaim {Theorem \chapno2.4b}. Choose any $\epsilon > 0$ and let $\points$ be
points in $\R1$ such that $n \ge 2$ and $\djk \ge \epsilon$
when $j \ne k$. If \ $\Z0$, then

$$ \sjk\djk \le -{1\over2} \ \epsilon \ysqr .$$ 
We shall now show that this bound is optimal. Without loss of
generality, we return to the case $\epsilon = 1$.
We take our points
to be the integers $0, 1, \ldots, n$, so that the Euclidean distance matrix,
$A_n$ say, is given by
$$ A_n = \pmatrix{ 0 & 1 & 2 & \ldots & n\cr
                   1 & 0 & 1 & \ldots & n-1\cr
                   \vdots & \vdots &\ddots & \vdots\cr
                   n & n-1 & n-2 & \ldots & 0\cr} .$$
It is straightforward to calculate the inverse of $A_n$:
$$ \ainv = \pmatrix{ (1-n)/2n&1/2& & & & 1/2n\cr
                      1/2&-1&1/2& & & \cr
                       &1/2&-1& & & \cr
                       & & & \ddots& &\cr
                       & & & &-1&1/2\cr
                      1/2n & & & &1/2&(1-n)/2n }
 .$$                      

\proclaim{Proposition \chapno2.5}. We have the inequality $2 - (\pi^2/2n^2) \le \aninv \le 2$.

\pf We observe that $\aninv\le\| \ainv \|_1 = 2$, establishing the upper bound.
For the lower bound, we focus attention on the $(n-1)\times(n-1)$ symmetric
tridiagonal minor of $\ainv$ formed by deleting its first and last rows and
columns, which we shall denote by $T_n$. Thus we have
$$
   T_n = \pmatrix{ -1&1/2& & & & \cr
                 1/2&-1&1/2& & & \cr
                  &1/2&-1& & & \cr
                  & & &\ddots& & \cr
                  & & & &-1&1/2\cr
                  & & & &1/2&-1 }
.$$
Now
$$\eqalign{
\| T_n \|_2
   &= \max\{ y^T\ainv y : y^T y=1 \hbox{ and } y_1 = y_{n+1} = 0 \} \cr
   &\le \max\{ y^T\ainv y : y^T y=1 \} \cr
   &= \aninv, \cr
}$$
so that $\| T_n \|_2 \le \aninv \le 2$.
But the eigenvalues of $T_n$ are given by
$$ \lambda_k = -1 + \cos(k\pi/n), \hbox{ for } k = 1, 2, \ldots, n-1.$$
Thus $\| T_n \|_2 = 1 - \cos(\pi-\pi/n) \ge 2 - \pi^2/2n^2,$ where we have
used an elementary inequality based on the Taylor series for the cosine function.
The proposition is proved.~\qed

\sect {\chapno3. The Multivariate Case for the Euclidean Norm}
We first prove the multivariate versions of Propositions \chapno2.1 and \chapno2.2, which
generalize in a very straightforward way. We shall require the fact
that the generalized Fourier transform of $\phi(x) = \| x \|$ in
$\Rd$
is given by
$$ \phihat(t) = -c_d \| t \|^{-d-1},$$ where $$ c_d = 2^d\pi^{(d-1)/2} \Gamma((d+1)/2) .$$ 
This may be found in Jones (1982), Theorem 7.32. We now deal with
the analogue of Proposition \chapno2.1.

\proclaim {Proposition \chapno3.1}. If \ $\Z0$, then
$$
 \sjk \djk = -c_d(2\pi )^{-d} \int_\Rd \yjk \|t\|^{-d-1}\,dt . \eqno{(\chapno3)}
$$

\pf We define  
$$ {\hat g}(t) = -c_d \| t \|^{-d-1} \yjk. $$ 
The condition $\Z0$ implies this function is uniformly bounded and
the decay for large argument is sufficient to ensure absolute
integrability. The argument now follows the proof of Proposition \chapno2.1, with 
obvious minor changes.~\qed

\proclaim {Proposition \chapno3.2}. Let $B: \Rd \rightarrow \R1$ be a
continuous
function such that $\hbox{supp}\!(B)$ is contained in the ball  $\{ x \in \Rd : \| x \| \le 1
\}$,
$0 \le {\hat B}(t) \le \| t \|^{-d-1}$ and $B(0) > 0$.
If $n \ge 2$, $\djk \ge 1$ for $j \ne k$, and $\Z0$, then
$$ \sjk \djk \le -c_d B(0) \ysqr .$$ 

\pf The proof of Proposition \chapno2.2 clearly generalizes to this
case.~\qed

However, to exhibit a function $B$ satisfying the conditions of
Proposition \chapno3.2 is harder than in the univariate case. We modify a construction
of Narcowich and Ward (1990) and Sun (1990). Let $$B_0(x) = \cases{ 1, &if $\| x \| \le 1/2$ \cr 0, 
&otherwise.} $$ Then, using Narcowich and Ward (1990), equation 1.10 or [9],
Lemma \chapno3.1, we find that
$${\hat B_0}(t) = (2\| t \|)^{-d/2} J_{d\over 2}(\|
t\|/2),$$ where $J_k$ denotes the $k^{th}$-order Bessel function of the
first kind. Further, ${\hat B_0}$ is a radially symmetric function
since $B_0$ is radially symmetric. We now define $$B = B_0 * B_0,$$ so that, by the convolution theorem, 
$$\eqalign{ \hat B(t) &= (\hat B_0)^2 (t) \cr
       &= (2\| t \|)^{-d}  J_{d\over 2}^2(\| t\|/2) , 
}$$
and the behaviour of $J_0$ for large argument provides the inequality
$$ \hat B(t) \le \mu_d \| t \|^{-d-1} , $$ for some constant $\mu_d$. Since the conditions of
Proposition \chapno3.2 are now easy to verify when $B$ is scaled by
$\mu_d^{-1}$, we see that we are done .

\sect{\chapno4. Fourier Transforms and Bessel Transforms}
Here we relate our technique to the work of Ball (1989) and Narcowich
and Ward (1990, 1991).

\proclaimdefn{Definition \chapno4.1}. A real sequence $(y_j)_{j\in\Zd}$ is said to be
{\it zero-summing} if it is finitely supported and $\sum_{j\in\Zd} y_j = 0$.

\proclaimdefn{Definition \chapno4.2}. A function $\phi\colon[0,\infty) \to \R1$
will be said to be  {\it conditionally
negative definite of order 1 on} $\Rd$, hereafter shortened to 
CND1($d$), if it is continuous and, for any points $(x_j)_{j\in\Zd}$ in $\Rd$
and any zero-summing sequence $(y_j)_{j\in\Zd}$, we have
$$ 
  \sum_{j, k\in\Zd} y_j y_k \phi(\djk) \le 0 .
$$

\noindent Such functions were characterized by von Neumann and Schoenberg (1941). 
For every positive integer $d$, let $\Omd \colon [0,\infty) \to \R1$ be defined by
$$ \Omd(r) = \omd^{-1} \int_\Sdm1 \cos(ryu)\,dy , $$
where $u$ may be any unit vector in $\Rd$, $\Sdm1$ denotes the unit sphere
in $\Rd$, and $\omd$ its $(d-1)$-dimensional Lebesgue measure. Thus $\Omd$
is essentially the Fourier transform of the normalized rotation invariant
measure on the unit sphere. 

\proclaim {Theorem \chapno4.3}. Let $\phi\colon[0,\infty)\to \R1$ 
be a continuous function. A necessary and sufficient condition that
$\phi$ be a CND1($d$) function is that
it have the form
$$ \phi(r) = \phi(0) + \hint \left( 1 - \Omd(rt) \right) t^{-2} d\beta(t) ,
$$for every $r \ge 0$,
where $\beta\colon[0,\infty)\to\R1$ is a non-decreasing function
such that $\int_1^\infty t^{-2}\,d\beta(t) < \infty$ and $\beta(0) = 0$.
Furthermore,  $\beta$ is uniquely determined by $\phi$.

\pf  The first part of this result is Theorem 7 of von Neumann and
Schoenberg (1941), restated in our terminology. 
The uniqueness of $\beta$ is a consequence of Lemma 2 of that paper.~\qed

\noindent It is a consequence of this theorem that there exist constants $A$
and $B$ such that $\phi(r) \le Ar^2 + B$. For we have
$$
 \left|\int_1^\infty \left( 1 - \Omd(rt) \right) t^{-2}\,d\beta(t) \right|
 \le 2 \int_1^\infty t^{-2}\,d\beta(t) < \infty ,
$$
using the fact that $|\Omd(r)| \le 1$ for every $r \ge 0$.
Further, we see that
$$
 0 \le 1- \Omd(\rho) = 2 \omd^{-1} \int_\Sdm1 \sin^2(\rho ut/2)\,dt \le \rho^2/2 ,
$$  
which provides the bound
$$
 \Bigl|\int_0^1 \left(1 - \Omd(rt)\right) t^{-2}\,d\beta(t)\Bigr|
    \le r^2 \int_0^1 \half\,d\beta(t) =  \half r^2 \beta(1).
$$  
Thus $A = \beta(1)/2$ and $B = \phi(0) + 2 \int_1^\infty t^{-2}\,d\beta(t)$ suffice.
Therefore the function $\{ \phi(\| x \|) : x \in \Rd\}$
is a tempered distribution in the sense of Schwartz (1966) and possesses
a generalized Fourier transform $\{\phihat(\| \xi \|): \xi \in\Rd\}$. There is
a rather simple relation between the generalized Fourier transform and the
nondecreasing function of Theorem \chapno4.3 for a certain class of functions. This is
our next topic.

\proclaim{Definition \chapno4.4}. A function $\phi\colon[0,\infty) \to \R1$ will
be termed {\it admissible} if it is a continuous function of algebraic growth
which satisfies the following conditions:
\item{1. } $\phihat$ is a continuous function on $\Rdm0$.
\item{2. } The limit $\lim_{\,\|\xi\| \to 0} \|\xi\|^{d+1}\phihat(\|\xi\|)$ exists.
\item{3. } The integral $\int_{\{\|\xi\| \ge 1\}} |\phihat(\|\xi\|)|\, d\xi$ exists.

\vskip 12pt

It is straightforward to prove the analogue of Propositions \chapno2.1 and \chapno3.1 for an admissible
function.

\proclaim {Proposition \chapno4.5}. Let $\phi\colon[0,\infty)\to\R1$ be an admissible
function and let $(y_j)_{j\in\Zd}$ be a zero-summing sequence. Then for any choice
of points $(x_j)_{j\in\Zd}$ in $\Rd$ we have the identity
$$ \sum_{j,k\in\Zd} y_j y_k \phi(\|x_j-x_k\|)
      = (2\pi)^{-d} 
  \int_\Rd \Bigl|\sum_{j\in\Zd} y_j \exp(ix_j \xi)\Bigl|^2 \phihat(\|\xi\|)\,d\xi .
\eqno{(\chapno4)}$$
         
 \pf Let $\ghat\colon \Rd \to \R1$ be the function defined by
 $$ \ghat(\xi) = \Bigl|\sum_{j\in\Zd} y_j \exp(ix_j \xi)\Bigl|^2 \phihat(\|\xi\|).$$
 Then $\ghat$ is an absolutely integrable function on $\Rd$, because of the conditions
 on $\phi$ and because $(y_j)_{j\in\Zd}$ is a zero-summing sequence. 
Thus $\ghat$ is the generalized transform of $\sum_{j,k} y_j y_k \phi(\|\cdot+x_j-x_k\|)$,
and by standard properties of 
generalized Fourier transforms we deduce that
 $$ \sum_{j,k} y_j y_k \phi(\|x + x_j - x_k\|)
    = (2\pi)^{-d} \int_\Rd 
        \Bigl|\sum_{j\in\Zd} y_j \exp(ix_j \xi)\Bigl|^2 \phihat(\|\xi\|)
\exp(ix\xi)\,d\xi .$$
 The proof is completed by setting $x = 0$.~\qed

\proclaim {Proposition \chapno4.6}. Let $\phi\colon[0,\infty)\to\R1$ be 
an admissible CND1($d$) function. Then
$$\,d\beta(t) = -(2\pi)^{-d}\omd\phihat(\|tu\|)t^{d+1}\,dt,
$$where $u$ may be any unit vector in $\Rd$.
 
\pf 
Let $\mu$ and $\nu$ be different integers and let $(y_j)_{j\in\Zd}$ be a sequence
with only two nonzero elements, namely $y_\mu = -y_\nu = 2^{-1/2}$. Choose any
point $\zeta \in \Rd$ and set $x_\mu = 0$, $x_\nu = \zeta$,
so that equation (\chapno4) provides the expression
$$ \phi(0) - \phi(\|\zeta\|) =
  (2\pi)^{-d} \int_\Rd \left( 1 - \cos(\zeta\xi)\right) \phihat(\|\xi\|)\,d\xi .$$
Employing spherical polar coordinates, this integral takes the form
 $$ \phi(0) - \phi(\|\zeta\|) = 
2\pi)^{-d} \omd \hint  \left(1-\Omd(t\|\zeta\|)\right) \phihat(\|t u\|) 
t^{d-1}\,dt ,$$
where $u$ may be any unit vector in $\Rd$.
Setting $r = \|\zeta\|$, we have 
$$\phi(r) = \phi(0) + \hint   \left(1-\Omd(rt)\right) \gamma(t) t^{-2}\,dt ,$$
where $\gamma(t) = -(2\pi)^{-d}\omd\phihat(\|tu\|)t^{d+1}$. 
Now Theorem 4.2.6 of the following chapter implies that $\phihat$ is a nonpositive function. Thus there exists
a nondecreasing function $\tilde\beta\colon[0,\infty) \to \R1$ such that
$\gamma(t)\,dt = d\tilde\beta(t)$, and 
$\int_1^\infty t^{-2}\,d\tilde\beta(t)$ is finite and $\tilde\beta(0) = 0$. But
the uniqueness of the representation of Theorem \chapno4.3 implies that $\beta =
\tilde\beta$, that is
$$\,d\beta(t) = -(2\pi)^{-d}\omd\phihat(\|tu\|)t^{d+1}\,dt, $$and the proof is complete.~\qed

This proposition is useful if we want to calculate $\beta$ for a particular
function $\phi$, since tables of generalized Fourier transforms are readily
available.

\proclaim{Example \chapno4.7}. Let $\phi(r) = (r^2 + 1)^{1/2}$.
This is a non-negative CND1($d$) function for all $d$ (see Micchelli (1986)).
When $d=3$, the generalized Fourier transform is $\phihat(r) = -4\pi r^{-2} K_2(r)$.
Here $K_2$ is a modified
Bessel function which is positive and smooth in $\R1^+$, has a pole at the
origin, and decays exponentially (See Abramowitz and Stegun (1970)). Consequently
$\phihat$ is a non-negative admissible function. Applying Theorem \chapno4.7 gives the
equation
$$\eqalign{
 d\beta(r) &=   (2\pi)^{-3} (4\pi) r^4 (4\pi r^{-2} K_2(r) )\cr   
           &=   (2r^2/\pi) K_2(r) dr,
}$$
agreeing with Narcowich and Ward (1991), equation 3.12.

\sect{\chapno5. The Least Upper Bound for Subsets of the Integer Grid}
In the next chapter we use extensions of the
technique provided here to derive the the following result.

\proclaim{Theorem \chapno5.1}. Let $\phi\colon[0,\infty)\to\R1$ be an admissible function 
that is not identically zero, let $\phi(0)\ge 0$, and let $\phi$ be CND1($d$) for 
every positive integer $d$. Further, let $(x_j)_{j\in\Zd}$ be any elements of $\Zd$
and let $A=(\phi(\|x_j-x_k\|))_{j,k\in{\cal N}}$, where ${\cal N}$ can
be any finite subset of $\Zd$. Then we have the inequality
$$ \|A^{-1}\| \le \Bigl( \sum_{k\in\Zd} |\phihat(\|\pi e + 2\pi k\|)| \Bigr)^{-1}, $$
where $e = [1, \ldots, 1]^T \in \Rd$ and $\phihat$ is the generalized
Fourier transform of $\phi$. Moreover, this is the least upper bound
valid for all finite subsets of $\Zd$.

\pf See Section 4.4 of the thesis.~\qed

\vfill \eject
\def\chapno{4.}

\headline{\ifnum\pageno=42\hfil\else{\headfont\hfil\chptitle\hfil}\fi}
\def\chptitle{Norm estimates for Toeplitz distance matrices I}
\bigskip
\centerline{\bf 4 : Norm estimates for Toeplitz distance matrices I}
\medskip


\def\pts{{(x_j)_{j=1}^n}}

\def\sjk{{\sum_{j,k=1}^n y_j y_k\,}}
\def\djk{{\| x_j - x_k \|}}
\def\xjk{{x_j-x_k}}
\def\yjk{{|\sum_{j=1}^n y_j e^{ix_j t}|^2}}
\def\trig{{|\sum_{j=1}^n y_j e^{ijt}|^2}}

\def\skzd{{\sum_{k\in\Zd}}}

\def\rint{{\int_{-\infty}^{\infty}}}
\def\hint{{\int_0^\infty}}

\def\Sdm1{{S^{d-1}}}

\def\ysqr{{\Vert y \Vert^2}}

\def\ainv{{A_n^{-1}}}
\def\aninv{{\Vert A_n^{-1} \Vert_2}}

\def\zsum{{\sum_{k=-\infty}^\infty}}

\def\hi{{(2\pi)^{-1}}}
\def\tint{{\hi\int_0^{2\pi}}}
\def\tintd{{(2\pi)^{-d}\int_{\Td}}}
\def\fhat{{\hat f}}
\def\ajxi{{\Bigl|\sum_{j\in\Zd} a_j \exp(ij\xi)\Bigl|^2}}
\def\Rdm0{{\Rd\setminus\{0\}}}
\def\ghat{{\hat g}}
\def\ytrigo{{\Bigl|\sum_{j\in\Zd} y_j \exp(ix_j \xi)\Bigl|^2}}
\def\ytrig2{{\Bigl|\sum_{j\in\Zd} y_j \exp(ij \xi)\Bigl|^2}}
\def\ynj{{y_j^{(n)}}}

\def\phichat{{\phihat_c}}
\def\chichat{{{\hat \chi}_c}}

\def\min{{\rm min}}


\sect {\chapno1. Introduction}
The multivariate interpolation problem is as follows: given
points
$\pts$ in $\Rd$ and real numbers $(f_j)_{j=1}^n$, construct
a function $s\colon\Rd\to\R1 $ such that $s(x_k)=f_k$, for $k=1,\ldots,n$.
The radial basis function approach is to choose a univariate function
$\phi\colon[0,\infty) \to \R1$, a norm $\|\,.\, \|$ on $\Rd$,
and to let $s$ take the form
$$ s(x) = \sum_{j=1}^n y_j \phi(\| x - x_j \|) . $$
The norm $\|\,.\,\|$ will be the Euclidean norm throughout this chapter. Thus
the radial basis function interpolation problem has a unique solution for any
given scalars $(f_j)_{j=1}^n$ if and only if the matrix
$(\phi(\djk))_{j,k=1}^n$
is invertible. Such a matrix will, as before, be called a distance matrix.
These functions provide a useful and flexible form for multivariate approximation,
but their approximation power as a space of functions is not addressed here.

A powerful and elegant theory was developed by I.~J.~Schoenberg and others some fifty
years ago which may be used to analyse the singularity of
distance matrices. Indeed, in Schoenberg 
(1938) it was shown that the Euclidean distance matrix, which is the case
$\phi(r)=r$, is invertible if $n\ge 2$ and the points $\pts$ are distinct.
Further, extensions of this work by Micchelli (1986)
proved that the distance matrix is invertible
for several classes of functions, including the Hardy
multiquadric, the only restrictions on the points $\pts$
being that they are distinct and that $n\ge 2$. Thus the singularity of the
distance matrix has been successfully investigated for 
many useful radial basis functions. In this chapter,
we bound the eigenvalue of smallest modulus for certain distance matrices.
Specifically, we provide the greatest lower bound on the moduli of the  eigenvalues
in the case when the points $\pts$ form a subset of the integers
$\Zd$,
our method of analysis applying to a wide class of functions which includes the
multiquadric. More precisely, let $N$ be any finite subset of the integers $\Zd$
and let $\lambda^N_{\min}$ be the smallest eigenvalue in modulus of the 
distance matrix $(\phi(\|j-k\|))_{j,k\in N}$. Then the results of Sections~3
and 4 provide the inequality
$$ |\lambda^N_{\min}| \ge C_\phi, \eqno{(\chapno1.1)}$$
where $C_\phi$ is a positive constant for which an elegant formula is derived.
We also provide a constructive proof that $C_\phi$ cannot be replaced by any
larger number, and it is for this reason that we shall describe inequality (\chapno1.1)
as an {\it optimal} lower bound. Similarly, we shall say that an upper bound is 
optimal if none of the constants appearing in the inequality can be replaced by smaller
numbers.

It is crucial to our analysis that the distance matrix $(\phi(\|j-k\|))_{j,k\in N}$
may be embedded in the bi-infinite matrix $(\phi(\|j-k\|))_{j,k\in \Zd}$.
Such a bi-infinite matrix is called a Toeplitz matrix if $d = 1$. We shall
use this name for all values of $d$, since we use the multivariate form of the
Fourier analysis of Toeplitz forms (see Grenander and Szeg\H{o} (1984)).

Of course, inequality (\chapno1.1) also provides an upper bound on the norm
of the inverse of the distance matrices generated by finite subsets of the
integers $\Zd$. This is not the first paper to address the 
problem of bounding
the norms of inverses of distance matrices and we acknowledge the papers
of Ball (1989) and Narcowich and Ward [1990, 1991], which first interested the
author in such estimates. Their results are not limited to the case when
the data points are a subset of the integers.
Instead, they apply when the points satisfy the condition $\djk \ge \epsilon$
for $j\ne k$, where $\epsilon$ is a positive constant, and they provide
lower
bounds on the smallest modulus of an eigenvalue for several functions $\phi$,
including the multiquadric. We will find that these bounds are not optimal,
except in the special case of the Euclidean norm in the univariate case.
Further, our bounds apply to all the conditionally negative definite functions
of order 1. The definition of this class of functions may be found in Section~\chapno3.

As in the previous section, we make extensive use of the theory of generalized Fourier transforms,
for which our principal reference will still be Jones (1982). 
These transforms
are precisely the Fourier transforms of tempered distributions constructed in 
Schwartz (1966). First, however, Section~2 presents several theorems which
require only the classical theory of the Fourier transform. These results will
be necessary in Section~\chapno3.

\sect{\chapno2. Toeplitz forms and Theta functions}
We require several properties of the Fej\'er kernel, which is defined as follows.
For each positive integer $n$, the $n^{\rm th}$ univariate 
Fej\'er kernel is the  positive trigonometric polynomial 
$$\eqalign{
   K_n(t) &=  \sum_{k=-n}^n \left( 1 - |k|/n \right) \exp(ikt) \cr
          &=  {{\sin^2\,nt/2}\over{n\,\sin^2\,t/2}} . 
}\eqno{(\chapno2.1)}$$ Further, the $n^{\rm th}$ multivariate Fej\'er kernel is defined
by the product
$$ K_n(t_1, \ldots, t_d) = K_n(t_1)K_n(t_2)\cdots K_n(t_d) , \qquad
 t \in \Rd. \eqno{(\chapno2.2)}$$         

\proclaim {Lemma \chapno2.1}. The univariate kernel enjoys the following property:
for any continuous $2\pi$-periodic function
$f\colon\R1\to\R1$ and for all $x \in \R1$ we have
$$ \lim_{n\to\infty} \tint K_n(t-x)f(t) \,dt = f(x). $$
Moreover, we have the equations
$$ (2\pi)^{-1} \int_0^{2\pi} K_n(t) \,dt = 1 \eqno{(\chapno2.3)}$$
and
$$ K_n(t) = \Bigl| n^{-1/2} \sum_{k=0}^{n-1} \exp(ikt) \Bigr|^2 . \eqno{(\chapno2.4)}$$

\pf Most text-books on harmonic analysis contain the first property and
(\chapno2.3). For example,
see pages 89ff, volume I, Zygmund (1979). It is elementary to deduce (\chapno2.4)
from (\chapno2.1).~\qed

\proclaim {Lemma \chapno2.2}. For every continuous $\Td$-periodic function
$f\colon\Rd \to\R1$, the multivariate Fej\'er kernel gives the convergence
property
$$ \lim_{n\to\infty} \tintd K_n(t- x)f(t)\,dt = f(x) $$
for every $x \in \Rd$. Further, 
$K_n$ is the square of the modulus of a trigonometric polynomial with real
coefficients and
$$ (2\pi)^{-d} \int_\Td K_n(t) \,dt = 1 .$$

\pf  The first property is Theorem 1.20 of chapter 17 of Zygmund (1979).
The last part of the lemma is an immediate consequence of (\chapno2.3), (\chapno2.4)
and the definition of the multivariate
Fej\'er kernel.~\qed

All sequences will be real sequences here. Further,
we shall say that a sequence $(a_j)_{\Zd} := \lbrace a_j \rbrace_{j\in\Zd} $ 
is {\it finitely supported}
if it contains only finitely many nonzero terms. The scalar product
of two vectors $x$ and $y$ in $\Rd$ will be denoted by $xy$. 

\proclaim {Proposition \chapno2.3}. Let $f\colon \Rd \to \R1$ be an absolutely integrable
continuous function whose Fourier transform $\fhat$ is also absolutely integrable. 
Then for
any finitely supported sequence $(a_j)_\Zd$, and for any choice of points
$(x_j)_\Zd$ in $\Rd$, we have the identity
$$ \sum_{j,k\in\Zd} a_j a_k f(x_j - x_k) = 
  (2\pi)^{-d} \int_\Rd \Bigl| \sum_{j\in\Zd} a_j \exp(ix_j \xi) \Bigl|^2 \fhat(\xi)\,d\xi .$$
  
\pf The function $\{ \sum_{j,k} a_j a_k f(x + x_j - x_k): x \in \Rd \}$ is absolutely
integrable. Its Fourier transform is given by
$$\eqalign{
   \Bigl[ \sum_{j,k\in\Zd} a_j a_k f(\cdot + x_j - x_k) \Bigr]^\land (\xi)
     &=  \sum_{j,k\in\Zd} a_j a_k \exp(i(x_j-x_k)\xi) \fhat(\xi) \cr
     &=  \Bigl| \sum_{j\in\Zd} a_j \exp(ix_j \xi) \Bigl|^2 \fhat(\xi) , \qquad \xi \in \Rd,
}$$
and is therefore absolutely integrable.
Therefore the Fourier inversion 
theorem states that
$$ 
\sum_{j,k\in\Zd} a_j a_k f(x + x_j - x_k)
 = (2\pi)^{-d} \int_\Rd 
     \Bigl| \sum_{j\in\Zd} a_j \exp(ix_j \xi) \Bigl|^2 \fhat(\xi) 
\exp(ix\xi)\,d\xi .$$ Setting $x=0$ produces the stated equation.~\qed

In this dissertation a key {\it r\^ole} will be played by the {\it symbol function}
$$ \sigma(\xi) = \sum_{k\in\Zd} \fhat(\xi+2\pi k), \qquad \xi \in \Rd .\eqno{(\chapno2.5)}$$
If $\fhat \in L^1(\Rd)$, then $\sigma$ is an absolutely integrable function on
$\Td$ and its defining series is absolutely convergent almost everywhere.
These facts are consequences of the relations
$$ \infty
 > \int_\Rd |\fhat(\xi)|\,d\xi  
 =   \sum_{k\in\Zd} \int_{\Td}  |\fhat(\xi+2\pi k)|\,d\xi 
 =   \int_{\Td} \sum_{k\in\Zd} |\fhat(\xi+2\pi k)|\,d\xi, 
$$
the exchange of integration and summation being a consequence of Fubini's theorem.
If the points $(x_j)_\Zd$ are integers, then we readily deduce the following
bounds on the quadratic form.

\proclaim {Proposition \chapno2.4}. Let $f$ satisfy the conditions of Proposition \chapno2.3
and let  $(a_j)_\Zd$ be a finitely supported sequence. Then we have the identity
$$ \sum_{j,k\in\Zd} a_j a_k f(j-k) = \tintd \Bigl|\sum_{j\in\Zd}
 a_j \exp(ij\xi)\Bigl|^2 \sigma(\xi)\,d\xi. \eqno{(\chapno2.6)} $$
Further, letting $m = \inf \lbrace \sigma(\xi) : \xi \in \Td \rbrace$ and
$M = \sup \lbrace \sigma(\xi) : \xi \in \Td \rbrace$, we have the bounds
$$ m \ \sum_{j\in\Zd} a_j^2 \le \sum_{j,k\in\Zd}
 a_j a_k f(j-k) \le M \ \sum_{j\in\Zd} a_j^2 .$$

\pf Proposition \chapno2.3 implies the equation
$$\eqalign{
   \sum_{j,k\in\Zd} a_j a_k f(j-k)
      &=  \sum_{k\in \Zd} \tintd \ajxi \fhat(\xi+2\pi k)\,d\xi \cr
      &=  \tintd \ajxi \sigma(\xi)\,d\xi,
 }$$
the exchange of integration and summation being justified by Fubini's theorem.
For the upper bound, the Parseval theorem yields the expressions
$$\eqalign{
   \sum_{j,k\in\Zd} a_j a_k f(j-k) 
      &=  \tintd  \ajxi \sigma(\xi)\,d\xi \cr
      &\le  M \sum_{j\in\Zd} a_j^2.
}$$
The lower bound follows similarly and the proof is complete.
\qed

The inequalities of the last proposition enjoy the following optimality
property.

\proclaim {Proposition \chapno2.5}. Let $f$ satisfy the conditions of Proposition \chapno2.3 and
suppose that the symbol function is continuous. Then the inequalities of Proposition
\chapno2.4 are optimal lower and upper bounds.

\pf Let $\xi_M \in \Td$ be a point such that $\sigma(\xi_M) = M$, which exists
by continuity of the symbol function.
We shall construct finitely supported sequences
$\{ (a_j^{(n)})_{j\in\Zd}: n = 1, 2, \ldots\, \}$
such that $\sum_{j\in\Zd} (a_j^{(n)})^2 = 1$, for all $n$, and
$$ \lim_{n \to \infty} \sum_{j,k\in\Zd} a_j^{(n)} a_k^{(n)} f(j-k)
   = M . \eqno{(\chapno2.7)}$$
   
We recall from Lemma \chapno2.2 that the multivariate Fej\'er kernel is the square of
the modulus of a trigonometric polynomial with real coefficients. 
Therefore there exists a finitely supported sequence $(a_j^{(n)})_\Zd$ satisfying the relation
$$ \Bigl|\sum_{j\in\Zd} a_j^{(n)} \exp(ij\xi) \Bigl|^2 = K_n(\xi-\xi_M),
\qquad \xi \in \Rd. \eqno{(\chapno2.8)}$$
Further, the Parseval theorem and Lemma \chapno2.2 provide the equations
$$ \sum_{j\in\Zd} (a_j^{(n)})^2 = (2\pi)^{-d} \int_\Td K_n(\xi - \xi_M)\,d\xi = 1 $$
and
$$
 \lim_{n \to \infty} \tintd K_n(\xi-\xi_M) \sigma(\xi)\,d\xi 
    = \sigma(\xi_M) = M.
$$
It follows from (\chapno2.6) and (\chapno2.8) that the limit (\chapno2.7) holds.
The lower bound of Proposition~\chapno2.4 is dealt with in the same fashion.~\qed

The set of functions satisfying the conditions of Proposition \chapno2.5 is nonvoid. For
example, suppose that we have  $\fhat(\xi) = \OO(\|\xi\|^{-d-\delta})$,
for large $\|\xi\|$,
where $\delta$ is a positive constant. Then the series defining the symbol function
$\sigma$ converges uniformly, by the Weierstrass M-test, and $\sigma$ is continuous, being a
uniformly convergent sum of continuous functions. These remarks apply when $f$ is a Gaussian,
which is the subject of the rest of this section. We shall see
that the analysis of the Gaussian provides the key to many of our results.

\proclaim {Proposition \chapno2.6}. Let $\lambda$ be a positive constant
and let $f(x) = \exp(-\lambda\|x\|^2)$, for $x \in \Rd$. Then $f$ satisfies
the conditions of Proposition \chapno2.5.

\pf The Fourier transform of $f$ is the function $\fhat(\xi) = (\pi/\lambda)^{d/2}
\exp(-\|\xi\|^2 /4\lambda)$, which is a standard calculation of the classical
theory of the Fourier transform. It is clear that $f$ satisfies the conditions
of Proposition \chapno2.3, and that the symbol function is the expression
$$ \sigma(\xi) = (\pi/\lambda)^{d/2}
       \sum_{k\in\Zd} \exp(-\|\xi + 2\pi k\|^2/4\lambda) , \quad \xi \in \Rd. 
\eqno{(\chapno2.9)}$$
Finally, the decay of the Gaussian ensures that $\sigma$ is continuous,
being a uniformly convergent sum of continuous functions.~\qed

This result is of little use unless we know the minimum and maximum values
of the symbol function for the Gaussian.
Therefore we show next that explicit expressions for these numbers may be calculated
from properties of Theta functions. Lemmata \chapno2.7 and \chapno2.8 address the cases when
$d=1$ and $d \ge 1$ respectively.

\proclaim {Lemma \chapno2.7}. Let $\lambda$ be a positive constant and let 
 $E_1\colon\R1 \to \R1$ be the $2\pi$-periodic
function 
$$ E_1(t) = \zsum \exp\left(-\lambda(t+2k\pi)^2\right) . $$
Then $E_1(0) \ge E_1(t) \ge E_1(\pi)$ for all $t \in \R1$.

\pf An application of the Poisson summation formula provides the relation
$$\eqalign{
  E_1(t) &=  (4\pi\lambda)^{-1/2} \zsum e^{-k^2/4\lambda} e^{ikt} \cr
              &=  (4\pi\lambda)^{-1/2} \left( 1 + 2 \sum_{k=1}^\infty e^{-k^2/4\lambda}
               \cos(kt)  \right) .
}$$
This is a Theta function. Indeed, using the notation of Whittaker
and Watson (1927), Section~21.11, it is a Theta function of Jacobi type
$$ \vartheta_3 (z, q ) = 1 + 2 \sum_{k=1}^\infty q^{k^2} \cos(2kz) ,$$
where $q \in \CC$ and $|q| < 1$. Choosing $q = e^{-1/4\lambda}$ we obtain the relation
$$ E_1(t) = (4\pi\lambda)^{-1/2} \vartheta_3( t/2, q ). $$
The useful product formula for $\vartheta_3$:
$$ \vartheta_3( z, q ) = G \prod_{k=1}^\infty (1 + 2q^{2k-1}\cos\,2z + q^{4k-2}), $$
where $G = \prod_{k=1}^\infty (1-q^{2k})$, is given in Whittaker
and Watson (1927), Sections~21.3 and 21.42. Thus
$$ E_1(t) = 
   (4\pi\lambda)^{-1/2} G \prod_{k=1}^\infty (1 + 2q^{2k-1}\cos\,t + q^{4k-2}),
\quad t \in \R1.$$
Now each term of the infinite product is a decreasing function on the interval
$[0,\pi]$, which implies that $E_1$  is a decreasing function on $[0,\pi]$.
Since $E_1$ is an even $2\pi$-periodic function, we
deduce that $E_1$ attains its global minimum at $t=\pi$ and its maximum 
at $t=0$.~\qed

\proclaim {Lemma \chapno2.8}. Let $\lambda$ be a positive constant and let 
$E_d \colon \Rd \to \Rd$ be the
$[0,2\pi]^d$-periodic function given by
$$ E_d(x) = \skzd \exp(-\lambda\| x + 2k\pi \|^2) .$$
Then $E_d(0) \ge E_d(x) \ge E_d(\pi e)$, where $e = [1, 1, \ldots,
1]^T$.

\pf The key observation is the equation
$$ E_d(x) = \prod_{k=1}^d E_1(x_k). $$
Thus $E_d(0) = \prod_{k=1}^d E_1(0) \ge \prod_{k=1}^d
E_1(x_k) = E_d(x) \ge \prod_{k=1}^d E_1(\pi) = E_d(\pi e)$,
using the previous lemma.~\qed       
       
These lemmata imply that in the Gaussian case the maximum and minimum
values of the symbol function occur at $\xi = 0$ and $\xi = \pi e$ 
respectively, where $e = [1, \ldots, 1]^T$. Therefore we deduce from 
formula (\chapno2.9) that the
constants of Proposition \chapno2.4 are the
expressions
$$\eqalign{
 m &=  (\pi/\lambda)^{d/2} \sum_{k\in\Zd} \exp(-\|\pi e + 2\pi k\|^2/4\lambda) \quad \hbox{ and } \cr
 M &=  (\pi/\lambda)^{d/2} \sum_{k\in\Zd} \exp(-\|\pi k\|^2/\lambda).
}\eqno{(\chapno2.10)}$$

\sect{\chapno3. Conditionally negative definite functions of order 1}
In this section we derive the optimal lower bound on the eigenvalue
moduli of the
distance matrices generated by the integers for a class of functions
including the Hardy multiquadric.

\proclaimdefn{Definition \chapno3.1}. A real sequence $(y_j)_\Zd$ is said to be
{\it zero-summing} if it is finitely supported and $\sum_{j\in\Zd} y_j = 0$.
\vskip 12pt

Let $\phi\colon[0,\infty) \to \R1$ be a continuous function of algebraic
growth. Thus it is meaningful to speak of the generalized Fourier
transform of the radially symmetric function $\{\phi(\|x\|): x \in \Rd\}$.
We denote this transform by $\{\phihat(\|\xi\|): \xi \in \Rd\}$, so emphasizing
that it is a radially symmetric distribution, but we note that $\phihat$
depends on $d$. We shall restrict attention
to the collection of functions described below.

\proclaimdefn{Definition \chapno3.2}. A function $\phi\colon[0,\infty) \to \R1$ will
be termed {\it admissible} if it is a continuous function of algebraic growth
which satisfies the following conditions:

\frenchspacing
\item{1. } $\phihat$ is a continuous function on $\Rdm0$.

\item{2. } The limit $\lim_{\,\|\xi\| \to 0} \|\xi\|^{d+1}\phihat(\|\xi\|)$ exists.

\item{3. } The integral $\int_{\{\|\xi\| \ge 1\}} |\phihat(\|\xi\|)|\, d\xi$ exists.

\vskip 12pt

It is straightforward to prove the analogue of Proposition \chapno2.3 for an admissible
function.

\proclaim {Proposition \chapno3.3}. Let $\phi\colon[0,\infty)\to\R1$ be an admissible
function and let $(y_j)_\Zd$ be a zero-summing sequence. Then for any choice
of points $(x_j)_\Zd$ in $\Rd$ we have the identity
$$ \sum_{j,k\in\Zd} y_j y_k \phi(\|x_j-x_k\|)
      = (2\pi)^{-d} 
  \int_\Rd \Bigl|\sum_{j\in\Zd} y_j \exp(ix_j \xi)\Bigl|^2 \phihat(\|\xi\|)\,d\xi .
\eqno{(\chapno3.1)}$$

 \pf Let $\ghat\colon \Rd \to \R1$ be the function defined by
 $$ \ghat(\xi) = \Bigl|\sum_{j\in\Zd} y_j \exp(ix_j \xi)\Bigl|^2 \phihat(\|\xi\|).$$
 Then $\ghat$ is an absolutely integrable function on $\Rd$, because of the conditions
 on $\phi$ and because $(y_j)_\Zd$ is a zero-summing sequence. 
Thus $\ghat$ is the generalized transform of $\sum_{j,k} y_j y_k \phi(\|\cdot+x_j-x_k\|)$,
and by standard properties of 
generalized Fourier transforms we deduce that
 $$ \sum_{j,k} y_j y_k \phi(\|x + x_j - x_k\|)
    = (2\pi)^{-d} \int_\Rd 
        \Bigl|\sum_{j\in\Zd} y_j \exp(ix_j \xi)\Bigl|^2 \phihat(\|\xi\|)
\exp(ix\xi)\,d\xi .$$
 The proof is completed by setting $x = 0$. \qed

 We come now to the subject that is given in the title of this section.

\proclaimdefn{Definition \chapno3.4}. Let $\phi\colon [0,\infty) \to \R1$ be a continuous
function. We shall say that $\phi$ is {\it conditionally negative definite of
order 1 on every} $\Rd$, hereafter shortened to CND1, if 
we have the inequality
$$ \sum_{j,k\in\Zd} y_j y_k \phi(\|x_j - x_k\|) \le 0, $$
for every positive integer $d$,
for every zero-summing sequence $(y_j)_\Zd$ and for any choice of points
$(x_j)_\Zd$ in $\Rd$.
\vskip  12pt   

Such functions were completely characterized by I. J. Schoenberg (1938).

\proclaim {Theorem \chapno3.5}. A continuous function $\phi\colon [0,\infty) \to \R1$
is CND1 if and only if there exists a nondecreasing function 
$\alpha\colon [0,\infty) \to \R1$ such that
$$ \phi(r) = \phi(0) + 
     \int_0^\infty \lbrack 1 - \exp(-tr^2)\rbrack  t^{-1}
d\alpha(t), \quad\hbox{ for } r > 0, $$ and the integral $
\int_1^\infty t^{-1}\,d\alpha(t)$ exists.

\pf This is Theorem 6 of Schoenberg (1938). \qed

\noindent Thus $d\alpha$ is a positive Borel measure such that
$$ \int_0^1 d\alpha(t) < \infty \hbox{ and } \int_1^\infty t^{-1} \,d\alpha(t) < \infty .$$
Further, it is a consequence of this theorem that there exist constants $A$ and $B$ such that
$\phi(r) \le Ar^2 + B$, where $A$ and $B$ are constants. In order to prove
this assertion we note the elementary inequalities
$$ \int_1^\infty \lbrack 1-\exp(-tr^2)\rbrack  t^{-1}\,d\alpha(t)
        \le  \int_1^\infty t^{-1}\,d\alpha(t) < \infty, $$
and
$$ \int_0^1  \lbrack 1-\exp(-tr^2)\rbrack  t^{-1}\,d\alpha(t)       
        \le r^2 \int_0^1\,d\alpha(t) . $$
Thus $A = r^2( \alpha(1)
- \alpha(0))$ and $B = \phi(0) + \int_1^\infty t^{-1}\,d\alpha(t)$ suffice. Therefore we may 
regard a CND1 function as a tempered distribution and it possesses a generalized
Fourier transform. The following relation between the transform and the integral
representation of Theorem \chapno3.5 will be essential to our needs.

\proclaim {Theorem \chapno3.6}. Let $\phi\colon[0,\infty) \to \R1$ be an admissible
CND1 function. For $\xi \in \Rdm0$, we have the formula
$$
   \phihat(\|\xi\|) 
    = -\!\int_0^\infty \exp(-\|\xi\|^2/4t) (\pi/t)^{d/2} t^{-1}\,d\alpha(t).
\eqno{(\chapno3.2)}
$$

Before embarking on the proof of this theorem, we require some groundwork. We
shall say
that a function $f\colon\Rdm0 \to \R1$ is {\it symmetric} if $f(-x) = f(x)$, for every
$x \in \Rdm0$.

\proclaim {Lemma \chapno3.7}. Let $\alpha \colon [0,\infty) \to \R1$ be a
nondecreasing function such that the integral $\int_1^\infty
t^{-1}\,d\alpha(t)$ exists. Then the function  
$$\psi(\xi) =  - \int_0^\infty
\exp(-\|\xi\|^2/4t)\/ (\pi/t)^{d/2} t^{-1}\,d\alpha(t), \quad \xi \in \Rdm0, \eqno{(\chapno3.3)}$$  
is a symmetric smooth function, that is every derivative exists.

\pf For every nonzero $\xi$, the limit
$$\lim_{t \to 0} \exp(-\|\xi\|^2/4t) (\pi/t)^{d/2} t^{-1} = 0$$
implies that the integrand of expression (\chapno3.3) is a continuous function on $[0,\infty)$.
Therefore it follows from the inequality 
$$\int_1^\infty \exp(-\|\xi\|^2/4t) (\pi/t)^{d/2} t^{-1}\,d\alpha(t)
\le \pi^{d/2} \int_1^\infty t^{-1}\,d\alpha(t) < \infty $$
that the integral is well-defined.
Further, a similar argument for nonzero $\xi$ shows that every derivative of the integrand
with respect to $\xi$ is
also absolutely integrable for $t \in [0,\infty)$, which implies that every derivative
of $\psi$ exists. The proof is complete, the symmetry of $\psi$ being obvious.
\qed

\def\aa{{\pmatrix{
\fhat(0) & \fhat(\xi) \cr
\fhat(\xi) & \fhat(0) \cr } }}

\def\vv{{\pmatrix{
a_\alpha \cr
a_\beta \cr } }}

\proclaim {Lemma \chapno3.8}. Let $f\colon\Rd \to \R1$ be a symmetric absolutely
integrable function such that
$$ \int_\Rd \Bigl|\sum_{j\in\Zd} a_j \exp(ix_j t)\Bigl|^2 f(t) \,dt = 0 ,$$
for every finitely supported sequence $(a_j)_\Zd$ and for any choice of points
$(x_j)_\Zd$. Then $f$ must vanish almost everywhere.

\pf The given conditions on $f$ imply that the Fourier transform 
$\fhat$ is a symmetric function that satisfies the equation
$$ \sum_{j,k\in\Zd} a_j a_k \fhat(x_j - x_k) = 0 , $$
for every finitely supported sequence $(a_j)_\Zd$ and for all points $(x_j)_\Zd$ 
in $\Rd$. Let $\alpha$ and $\beta$ be different integers and let $a_\alpha$ and
$a_\beta$ be the only nonzero elements of $(a_j)_\Zd$.
We now choose any point $\xi \in \Rdm0$ and set
$x_\alpha = 0$, $x_\beta = \xi$, which provides the equation
$$ \vv^T \aa \vv = 0, \qquad \hbox{ for all } a_\alpha, a_\beta \in \R1    .$$
Therefore $\fhat(0)= \fhat(\xi) = 0$, and since $\xi$ was arbitrary, 
$\fhat$ can only be the zero function. Consequently $f$ must
vanish almost everywhere.~\qed

\proclaim {Corollary \chapno3.9}. Let $g\colon\Rdm0 \to \R1$ be a symmetric continuous 
function such that
$$
 \int_\Rd \ytrigo |g(\xi)|\,d\xi < \infty \eqno{(\chapno3.4)}
$$
and
$$
 \int_\Rd \ytrigo g(\xi)\,d\xi = 0, \eqno{(\chapno3.5)}
$$
for every zero-summing sequence $(y_j)_\Zd$ and for any choice of points $(x_j)_\Zd$.
Then $g(\xi) = 0$ for every $\xi \in \Rdm0$.

\pf For any integer $k \in \{1, \ldots, d\}$ and for any positive real number
$\lambda$, let $h$ be the symmetric function
$$ h(\xi) = g(\xi) \sin^2 \lambda \xi_k, \qquad \xi \in \Rdm0. $$
The relation
$$ h(\xi) = g(\xi)\, \Bigl|\half \exp(i\lambda \xi_k) - \half \exp(-i\lambda\xi_k)\Bigl|^2
$$
and condition (\chapno3.4) imply that $h$ is absolutely integrable.

Let $(a_j)_\Zd$ be any real finitely supported sequence and let $(b_j)_\Zd$
be any sequence of points in $\Rd$. We define a real sequence $(y_j)_\Zd$
and points $(x_j)_\Zd$ in $\Rd$ by the equation
$$ \sum_{j\in\Zd} y_j \exp(ix_j\xi) = 
     \sin \lambda\xi_k \,\sum_{j\in\Zd} a_j \exp(ib_j\xi). $$
Thus $(y_j)_\Zd$ is a sequence of finite support. Further, setting $\xi = 0$,
we deduce that $\sum_{j\in\Zd} y_j = 0$, so $(y_j)_\Zd$ is a zero-summing
sequence. By condition (\chapno3.5), we have 
$$
 0  = \int_\Rd \Bigl|\sum_{j\in\Zd} y_j \exp(ix_j\xi)\Bigl|^2 g(\xi)\,d\xi 
    = \int_\Rd \Bigl|\sum_{j\in\Zd} a_j \exp(ib_j\xi)\Bigl|^2 h(\xi)\,d\xi .
$$
Therefore we can apply Lemma \chapno3.8 to $h$, finding that it vanishes
almost everywhere. Hence the continuity of $g$ for nonzero argument implies that
$g(\xi)\sin^2 \lambda\xi_k = 0$ for $\xi \ne 0$. But for every nonzero $\xi$ there exist
$k \in \{1, \ldots, d\}$ and $\lambda > 0$ such that $\sin \lambda\xi_k \ne 0$. 
Consequently $g$ vanishes on $\Rdm0$.~\qed

\noindent We now complete the proof of Theorem \chapno3.6.

\vskip  1mm  
\noindent
{\it Proof of Theorem \chapno3.6.\ } Let $(y_j)_\Zd$ be a zero-summing sequence and let
$(x_j)_\Zd$ be any set of points in $\Rd$. 
Then Theorem \chapno3.5 provides the expression
$$ \sum_{j,k\in\Zd} y_j y_k \phi(\|x_j - x_k\|)
  = - \int_0^\infty \Bigl(\sum_{j,k\in\Zd} y_j y_k \exp(-t\|x_j-x_k\|^2) \Bigr)
t^{-1}\,d\alpha(t), $$
this integral being well-defined because of the condition $\sum_{j\in\Zd} y_j =
0$. Therefore, using Proposition \chapno2.3 with $f(\cdot) = \exp(-t\|\cdot\|^2)$ 
in order to restate the Gaussian quadratic form in the
integrand, we find the equation

$$ \sum_{j,k\in\Zd} y_j y_k \phi(\|x_j - x_k\|) \hbox{\hskip 7cm} $$
$$\eqalign{  
  &=     -\!\int_0^\infty \Bigl\lbrack (2\pi)^{-d} \int_\Rd \ytrigo
(\pi/t)^{d/2}\exp(-\|\xi\|^2/4t) \,d\xi \Bigr\rbrack t^{-1}\,d\alpha(t)     \cr   
  &=     (2\pi)^{-d} \int_\Rd \ytrigo \psi(\xi)\,d\xi,
}$$   
where we have used Fubini's theorem to exchange the order of integration and
where $\psi$ is the function defined in (\chapno3.3).
By comparing this equation with the assertion of Proposition \chapno3.3,
we see that the difference $g(\xi) = \phihat(\|\xi\|)-\psi(\xi)$
satisfies the conditions of Corollary \chapno3.9. Hence $\phihat(\|\xi\|) =
\psi(\xi)$ for all $\xi \in \Rdm0$. The proof is complete.~\qed

\noindent{\sc Remark. }An immediate consequence of this theorem is that the
generalized Fourier transform of an admissible CND1 function cannot change
sign. 
\vskip  3mm  

The appearance of the Gaussian quadratic form in the proof of Theorem \chapno3.6
enables us to use the bounds of Lemma \chapno2.8, which gives the following result. 

\proclaim {Theorem \chapno3.10}. Let $\phi\colon[0,\infty) \to \R1$ be an
admissible CND1 function and let $(y_j)_\Zd$ be a zero-summing sequence.
Then we have the inequality
$$ \Bigl| \sum_{j,k\in\Zd} y_j y_k \phi(\|j-k\|) \Bigl| \ge
   |\sigma(\pi e)| \sum_{j\in\Zd} y_j^2 , $$
where $e = [1, \ldots, 1]^T$.

\pf Applying (\chapno3.1) and dissecting $\Rd$ into integer translates
of $\Td$,
we obtain the equations
$$\eqalign{ 
   \Bigl| \sum_{j,k\in\Zd} y_j y_k \phi(\|j-k\|) \Bigl|
 &=     (2\pi)^{-d} \int_\Rd \ytrig2 |\phihat(\|\xi\|)|\,d\xi     \cr   
 &=     \tintd \ytrig2 |\sigma(\xi)|\,d\xi,
} \eqno{(\chapno3.6)}$$  
where the interchange of summation and integration is justified by Fubini's theorem,
and where we have used the fact that $\phihat$ does not change sign. Here the symbol function
has the usual form (\chapno2.5).
Further, using (\chapno3.2), we 
again apply Fubini's theorem to deduce the formula
$$\eqalign{   
|\sigma(\xi)|
  &=      \sum_{k\in\Zd} |\phihat(\|\xi + 2\pi k\|)|     \cr   
  &=      \int_0^\infty  \Bigl(\sum_{k\in\Zd} \exp(-\|\xi+2\pi k\|^2/4t)\Bigr)\ 
(\pi/t)^{d/2} t^{-1} d\alpha(t) .
}$$   
It follows from Lemma \chapno2.8 that we have the bound
$$\eqalign{   
|\sigma(\xi)|
 &\ge \int_0^\infty  
         \Bigl(\sum_{k\in\Zd} \exp(-\|\pi e +2\pi k\|^2/4t)\Bigr)\  (\pi/t)^{d/2} t^{-1}
d\alpha(t)     \cr   
 &=      |\sigma(\pi e)|.
}  \eqno{(\chapno3.7)}$$   
The required inequality is now a consequence of (\chapno3.6) and the Parseval relation
$$ (2\pi)^{-d} \int_\Td \ytrig2 \,d\xi = \sum_{j\in\Zd} y_j^2.$$
\qed

When the symbol function is continuous on $\Rd \setminus 2\pi\Zd$, we can
show that the previous inequality is optimal using a modification of the proof
of Proposition \chapno2.5. Specifically, we construct a set $\{(\ynj)_\Zd:
n=1,2,\ldots\}$ of zero-summing sequences such that $\lim_{n\to\infty}
\sum_{j\in\Zd} (y_j^{(n)})^2 = 1$ and 
$$\lim_{n\to\infty}
  \Bigl|\sum_{j,k\in\Zd} y^{(n)}_j y^{(n)}_k \phi(\|j-k\|)\Bigl|
    =  |\sigma(\pi e)|,
$$   
which implies that we cannot replace $|\sigma(\pi e)|$ by any larger number
in Theorem \chapno3.10.

\proclaim {Corollary \chapno3.11}. Let $\phi\colon [0,\infty) \to \R1$ satisfy the
conditions of Theorem \chapno3.10 and let the symbol function be continuous
in the set $\Rd \setminus 2\pi\Zd$. Then the
bound of Theorem \chapno3.10 is optimal.

\pf  Let $m$ be an integer such that $4m \ge d+1$ and let $S_m$ be the 
trigonometric polynomial
$$ S_m(\xi) = [d^{-1}\sum_{j=1}^d \sin^2(\xi_j/2)]^{2m} , \qquad \xi \in \Rd .$$
Recalling from Lemma \chapno2.2 that the multivariate Fej\'er kernel is the square of 
the modulus of a trigonometric polynomial with real coefficients, 
we choose a finitely supported sequence $(\ynj)_\Zd$ satisfying the equations
$$ \Bigl|\sum_{j\in\Zd} \ynj \exp(ij\xi)\Bigl|^2 = K_n(\xi-\pi e) S_m(\xi), 
     \quad \xi \in \Rd. \eqno{(\chapno3.8)}
$$
Further, setting $\xi = 0$ we see that $(\ynj)_\Zd$ is a zero-summing
sequence. Applying (\chapno3.6), we find the relation
$$
  \Bigl|\sum_{j,k\in\Zd} y^{(n)}_j y^{(n)}_k \phi(\|j-k\|)\Bigl|
    = (2\pi)^{-d} \int_\Td K_n(\xi-\pi e) S_m(\xi)\,|\sigma(\xi)|\,d\xi. \eqno{(\chapno3.9)}
$$
Moreover, because the second condition of Definition \chapno3.2 implies that $S_m |\sigma|$ is a
continuous function, Lemma \chapno2.2 provides the equations
$$
 \lim_{n\to\infty}
  (2\pi)^{-d} \int_\Td K_n(\xi-\pi e) S_m(\xi)\,|\sigma(\xi)|\,d\xi
    = S_m(\pi e)\,|\sigma(\pi e)| = |\sigma(\pi e)|.
$$
It follows from (\chapno3.9) that we have the limit
$$
 \lim_{n\to\infty}
  \Bigl|\sum_{j,k\in\Zd} y^{(n)}_j y^{(n)}_k \phi(\|j-k\|)\Bigl|
   = |\sigma(\pi e)|.
$$

Finally, since $S_m$ is a continuous function, another application of Lemma \chapno2.2
yields the equation
$$
 \lim_{n\to\infty}  (2\pi)^{-d} \int_\Td K_n(\xi-\pi e)S_m(\xi)\,d\xi
   = S_m(\pi e) = 1. $$
By substituting expression (\chapno3.8) into the left hand side and employing the
Parseval relation
$$ (2\pi)^{-d} \int_\Td \Bigl|\sum_{j\in\Zd} \ynj \exp(ij\xi)\Bigl|^2 \,d\xi 
    =  \sum_{j\in\Zd} (\ynj)^2
$$
we find the relation $\lim_{n\to\infty} \sum_{j\in\Zd} (y_j^{(n)})^2 = 1$.~\qed

\sect{\chapno4. Applications}
This section relates the optimal inequality given in Theorem \chapno3.10 to the spectrum
of the distance matrix, using an approach due to Ball (1989). We apply
the following theorem.

\proclaim {Theorem \chapno4.1}. Let $A \in \R1^{n\times n}$ be a symmetric matrix with
eigenvalues $\lambda_1 \ge \cdots \ge \lambda_n$. Let $E$ be any subspace of
$\R1^n$ of dimension $m$. Then we have the inequality
$$ \max \{ x^T A x : x^T x = 1, x \perp E \} \ge \lambda_{m+1} .$$

\pf This is the Courant-Fischer minimax theorem. See Wilkinson (1965), pages
99ff.~\qed

For any finite subset $N$ of $\Zd$, let $A_N$ be the distance matrix
$(\phi(\|j-k\|))_{j,k\in N}$. Further, 
let the eigenvalues of $A_N$ be $\lambda_1 \ge \cdots \ge \lambda_{|N|}$, 
where $|N|$ is the cardinality of $N$, and let $\lambda^N_{\min}$
be the smallest eigenvalue in modulus.

\proclaim {Proposition \chapno4.2}. Let $\phi\colon [0,\infty) \to \R1$ be a
CND1 function that is not identically zero. Let $\phi(0) \ge 0$ and
let $\mu$ be a positive constant such that
$$ \sum_{j,k\in\Zd} y_j y_k \phi(\|j-k\|) \le -\mu \sum_{j\in\Zd} y_j^2,
\eqno{(\chapno4.1)} $$
for every zero-summing sequence $(y_j)_\Zd$. Then for every finite subset
$N$ of $\Zd$ we have the bound $$|\lambda^N_{\min}| \ge \mu . $$

\pf Equation (\chapno4.1) implies that
$$ y^T\!A_N y \le - \mu\,y^T y, $$
for every vector $(y_j)_{j\in N}$ such that $\sum_{j\in N} y_j = 0$.
Thus Theorem \chapno4.1 implies that the eigenvalues of $A_N$ satisfy $-\mu \ge \lambda_2 \ge \cdots
\ge \lambda_{|N|}$, where the subspace $E$ of that theorem is simply
the span of the vector $[1, 1, \ldots, 1]^T \in \R1^N$. 
In particular, $0 > \lambda_2 \ge \cdots \ge \lambda_{|N|}$. 
This observation and the condition $\phi(0) \ge 0$ 
provide the expressions
$$ 0 \le \hbox{trace}\,A_N = \lambda_1 + \sum_{j=2}^{|N|} \lambda_j
       = \lambda_1 - \sum_{j=2}^{|N|} |\lambda_j|. $$
Hence we have the relations $\lambda^N_\min = \lambda_2 \le -\mu$. The 
proof is complete.~\qed

We now turn to the case of the multiquadric $\phi_c(r)=(r^2+c^2)^{1/2}$, in
order to furnish a practical example of the above theory. This is a
non-negative CND1 function (see Micchelli (1986)) and its generalized Fourier
transform is the expression
$$ \phichat(\|\xi\|) = -\pi^{-1} (2\pi c/\|\xi\|)^{(d+1)/2} K_{(d+1)/2}(c\|\xi\|), $$
for nonzero $\xi$, which may be found in Jones (1982).
Here $\{ K_\nu(r): r > 0 \}$ is a modified
Bessel function which is positive and smooth in $\R1^+$, has a pole at the
origin, and decays exponentially (Abramowitz and Stegun (1970)).
Consequently, $\phi_c$ is a non-negative admissible CND1 function. 
Further, the exponential decay of $\phihat_c$ ensures that the symbol function 
$$ \sigma_c(\xi) = \sum_{k \in \Zd} \phihat_c(\|\xi+2\pi k\|) \eqno{(\chapno4.2)}$$
is continuous for $\xi \in \Rd \setminus 2\pi\Zd$. Therefore, given any finite
subset $N$ of $\Zd$, Theorem~\chapno3.10 and Proposition~4.2 imply
that the distance matrix $A_N$ has every
eigenvalue bounded away from zero by at least 
$$ \mu_c = \sum_{k\in\Zd} |\phihat_c(\|\pi e + 2\pi k\|)|, \eqno{(\chapno4.3)} $$
where $e = [1,1,\ldots,1]^T \in \Rd$. Moreover, Corollary \chapno3.11 shows that this bound is
optimal.

It follows from (\chapno4.3) that $\mu_c \to 0$ as $c \to
\infty$, because of the exponential decay of the modified Bessel
functions for large argument. For example, in the univariate case we have
the formula
$$
 \mu_c = (4c/\pi)\Bigl[K_1(c\pi) + K_1(3c\pi)/3 + K_1(5c\pi)/5 + \cdots \ \Bigr] ,
$$
and Table 4.1 displays some values of $\mu_c$. Of course, a practical
implication of this result is that we cannot expect accurate direct solution of
the interpolation equations for even quite modest values of $c$, at least
without using some special technique.
$$\vbox{
\settabs\+\indent&c\qquad\hskip 24pt&Optimal bound\quad&\cr
\+&{\it c}&{\rm Optimal bound}\cr
\+& $1.0$& $4.319455 \times 10^{-2}$\cr
\+& $2.0$& $2.513366 \times 10^{-3}$\cr
\+& $3.0$& $1.306969 \times 10^{-4}$\cr
\+& $4.0$& $6.462443 \times 10^{-6}$\cr
\+& $5.0$& $3.104941 \times 10^{-7}$\cr
\+& $10.0$& $6.542373 \times 10^{-14}$\cr
\+& $15.0$& $2.089078 \times 10^{-20}$\cr
}$$
\centerline{{\bf Table 4.1:}\quad The optimal bound on the smallest 
eigenvalue as $c \to \infty$}
\smallskip
The optimal bound is achieved only when the numbers of centres is
infinite. Therefore it is interesting
to investigate how rapidly $|\lambda^N_\min|$ converges to the
optimal lower bound as $|N|$ increases. Table 4.2
displays $|\lambda^N_\min| = \mu_c(n)$, say, for the distance
matrix $(\phi_c(\|j-k\|))_{j,k=0}^{n-1}$ for several values of $n$ when
$c = 1$. The third
column lists close estimates of $\mu_c(n)$ obtained using a theorem of Szeg\H{o}
(see Section~5.2 of Grenander
and Szeg\H{o} (1984)). Specifically, Szeg\H{o}'s theorem
provides the approximation
$$ \mu_c(n) \approx \sigma_c(\pi+\pi/n), $$
where $\sigma_c$ is the function defined in (\chapno4.2). This theorem
of Szeg\H{o} requires the fact that the minimum value of the symbol function
is attained at $\pi$, which is inequality (\chapno3.7).
Further, it provides the estimates 
$$ \lambda_{k+1} \approx \sigma_c(\pi+k\pi/n), \qquad k = 1, \ldots, n-1,$$
for {\it all} the negative eigenvalues of the distance matrix.
Figure \chapno1 displays the numbers $\{ -1/\lambda_k : k = 2, \ldots, n \}$
and their estimates $\{-1/\sigma(\pi+k\pi/n): k = 1, \ldots, n-1 \}$ in the case
when $n = 100$. We see that the agreement is excellent.
Furthermore, this modification of the classical theory of Toeplitz forms
also provides an interesting and useful perspective on the construction of
efficient preconditioners for the conjugate gradient solution of the
interpolation equations. We include no further information on these topics,
this last paragraph being presented as an {\it ap\'eritif} to the paper of
Baxter (1992c).

$$\vbox{
\settabs\+\indent&n\hskip 24pt&$4.324685 \times 10^{-2}$\hskip 24pt&$4.324685 \times 10^{-2}$
\qquad&\cr
\+&{\it n}& \hfill$\mu_1(n)$\hfill & \hfill$\sigma_1(\pi+\pi/n)$\hfill &\cr
\+& $100$ & \hfill $4.324685 \times 10^{-2}$ \hfill & \hfill $4.324653 \times 10^{-2}$ \hfill& \cr
\+& $150$ & \hfill $4.321774 \times 10^{-2}$ \hfill & \hfill $4.321765 \times 10^{-2}$ \hfill& \cr
\+& $200$ & \hfill $4.320758 \times 10^{-2}$ \hfill & \hfill $4.320754 \times 10^{-2}$ \hfill& \cr
\+& $250$ & \hfill $4.320288 \times 10^{-2}$ \hfill & \hfill $4.320286 \times 10^{-2}$ \hfill& \cr
\+& $300$ & \hfill $4.320033 \times 10^{-2}$ \hfill & \hfill $4.320032 \times 10^{-2}$ \hfill& \cr
\+& $350$ & \hfill $4.319880 \times 10^{-2}$ \hfill & \hfill $4.319879 \times 10^{-2}$ \hfill& \cr
}$$
\centerline{{\bf Table 4.2:}\quad Some calculated and estimated values
of $\lambda^N_\min$ when $c=1$}
\smallskip

\vfill
\eject
\centerline{\psfig{file=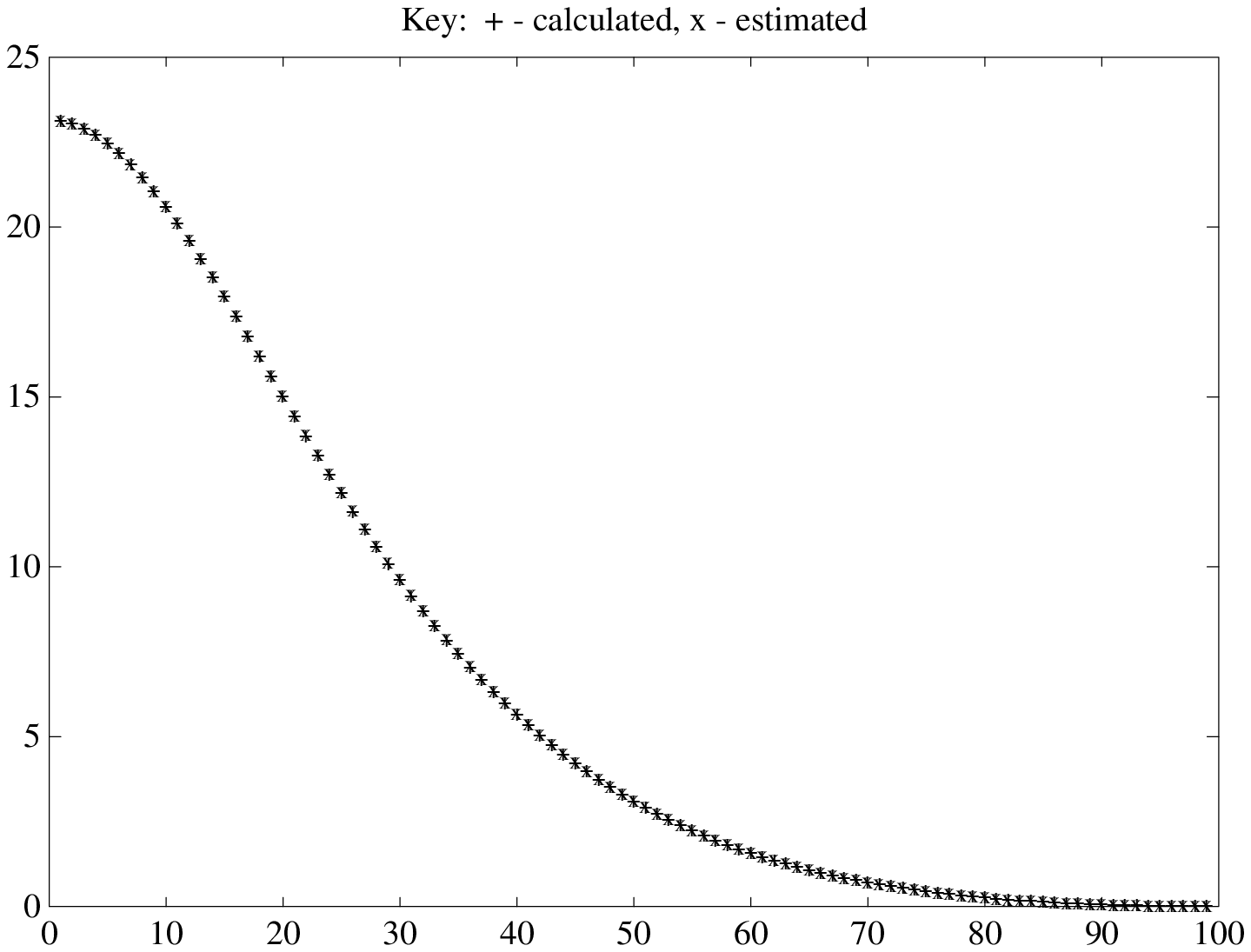}}
\centerline{{\bf Figure \chapno1.}\  Spectral estimates for a 
distance matrix of order 100}

\def\trig{{|\sum_{j\in\Zd} y_j \exp(ij\xi)|^2}}
\def\Td{{[0,2\pi]^d}}
\def\fhat{{\widehat f}}

\sect{\chapno5. A stability estimate}
The purpose of this last note is to derive an optimal inequality of the form
$$ \int_\Rd \left|\sum_{j\in\Zd} y_j \phi(\|x-j\|)\right|^2 \,dx 
       \ge C_\phi \, \sum_{j\in\Zd} y_j^2, $$
where $(y_j)_{j\in\Zd}$ is a real sequence of finite support such that
$\sum_{j\in\Zd} y_j = 0$, and $\phi\colon[0,\infty) \to \R1$ belongs to a certain class of functions
including the multiquadric. Specifically, this is the class of {\it
admissible CND1} functions. 
These functions have 
generalized Fourier transforms given by
$$ \phihat(\|\xi\|) = - \int_0^\infty \exp(-\|\xi\|^2/t) \,d\mu(t), \eqno{(\chapno5.1)}, $$
where $d\mu$ is a positive (but not finite) Borel measure on $[0,\infty)$. A derivation of this
expression may be found in Theorem \chapno2.6 above. 

\proclaim {Lemma \chapno5.1}. Let  $(y_j)_{j\in\Zd}$ be a zero-summing
sequence and let $\phi\colon[0,\infty) \to \R1$ be an admissible CND1
function. Then we have the equation
$$ \int_\Rd \left|\sum_{j\in\Zd} y_j \phi(\|x-j\|)\right|^2 \,dx 
       = (2\pi)^{-d} \int_\Td \trig \sigma(\xi) \,d\xi, $$
where $\sigma(\xi) = \sum_{k\in\Zd} |\phihat(\|\xi+2\pi k\|)|^2$.

\pf Applying the Parseval theorem and dissecting $\Rd$ into copies of the cube $\Td$, we 
obtain the equations
$$
 \int_\Rd |\sum_{j\in\Zd} y_j \phi(\|x-j\|)|^2 \,dx \hbox{\hskip 2in}
$$
$$\eqalign{
   &= (2\pi)^{-d} \int_\Rd \trig |\phihat(\|\xi\|)|^2 \,d\xi \cr
   &= \sum_{k\in\Zd} (2\pi)^{-d} \int_\Td \trig |\phihat(\|\xi+2\pi k\|)|^2 \,d\xi \cr
   &= (2\pi)^{-d} \int_\Td \trig \sigma(\xi) \,d\xi, 
}$$
where the interchange of summation and integration is justified by Fubini's theorem.\qed

If $\sigma(\xi) \ge m$ for almost every point $\xi$ in $\Td$, then the import of
Lemma \chapno5.1 is the bound
$$ \int_\Rd \left|\sum_{j\in\Zd} y_j \phi(\|x-j\|)\right|^2 \,dx
        \ge m \sum_{j\in\Zd} y_j^2 .$$
We shall prove that we can take $m = \sigma(\pi e)$, where $e = [1,1,\ldots,1]^T \in \Rd$.
Further, we shall show that the inequality is optimal if the function $\sigma$ is continuous
at the point $\pi e$.

Equation (\chapno5.1) is the key to this analysis, just as before. We see that
$$ |\phihat(\|\xi\|)|^2 =
   \int_0^\infty \int_0^\infty \exp( -\|\xi\|^2 (t_1^{-1}+t_2^{-1}) ) \,d\mu(t_1)\,d\mu(t_2),$$
whence, 
$$ \sigma(\xi) = 
 \int_0^\infty \int_0^\infty \sum_{k\in\Zd}
             \exp( -\|\xi+2\pi k\|^2 (t_1^{-1}+t_2^{-1}) )
                                          \,d\mu(t_1)\,d\mu(t_2), \eqno{(\chapno5.2)},
$$
where the interchange of summation and integration is justified by Fubini's theorem.

Now it is proved in Lemma \chapno1.8 that
$$ \sum_{k\in\Zd} \exp(-\lambda\|\xi+2\pi k\|^2)
      \ge \sum_{k\in\Zd} \exp(-\lambda\|\pi e+2\pi k\|^2), $$
for any positive constant $\lambda$. Therefore equation (\chapno5.2) provides the inequality
$$\eqalign{
 \sigma(\xi) &\ge 
         \int_0^\infty \int_0^\infty \sum_{k\in\Zd}
             \exp( -\|\pi e+2\pi k\|^2 (t_1^{-1}+t_2^{-1}) )
                 \,d\mu(t_1)\,d\mu(t_2), \cr
             &= \sigma(\pi e),
}$$
which is the promised value of the lower bound $m$ on $\sigma$ mentioned above.
Thus we have proved the following theorem.

\proclaim {Theorem \chapno5.2}. Let $(y_j)_{j\in\Zd}$, $\phi$ and $\sigma$ be as defined in Lemma 1.
Then we have the inequality
$$ \int_\Rd \left|\sum_{j\in\Zd} y_j \phi(\|x-j\|)\right|^2 \,dx 
       \ge \sigma(\pi e)\, \sum_{j\in\Zd} y_j^2 .$$

The proof that this bound is optimal uses the technique of Theorem \chapno2.11.

\def\ynj{{y_j^{(n)}}}
\proclaim {Theorem \chapno5.3}. The inequality of Theorem 2 is optimal if $\sigma$ is continuous
at $\pi e$.

\pf The condition that $\phi$ be admissible requires the existence of
the limit $\lim_{\|\xi\|\to 0} \|\xi\|^{d+1} \phihat(\|\xi\|) $.Let $m$ be a positive
integer such that $2m \ge d+1$ and let us define a sequence $\{ (\ynj)_{j\in\Zd} :
n = 1, 2, \ldots \}$ by
$$ \left| \sum_{j\in\Zd} \ynj \exp(ij\xi) \right|^2 = 
       \left( d^{-1} \sum_{j=1}^d \sin^2(\xi_j/2) \right)^{2m} \, K_n(\xi - \pi e), $$
where $K_n$ denotes the multivariate Fej\'er kernel. The standard properties of the
Fej\'er kernel needed for this proof are described in Lemma \chapno1.2. They allow
us to deduce that $(\ynj)_{j\in\Zd}$ is a zero-summing for every $n$. Further, we see that
$$\eqalign{
      \sum_{j\in\Zd} |\ynj|^2 &= (2\pi)^{-d} \int_\Td K_n(\xi-\pi e)
             (d^{-1} \sum_{j=1}^d \sin^2(\xi_j/2) )^{2m} \,d\xi \cr
         &= 1, \qquad\hbox{ for } n \ge 4m .
}$$
Finally, $m$ has been chosen so that the function
$$ \{ \left(d^{-1} \sum_{j=1}^d \sin^2(\xi_j/2)\right)^{2m} \, \sigma(\xi) :
          \xi \in \Td \} $$
is continuous. Therefore, 
we have
$$
    \lim_{n \to \infty}  \int_\Rd \left|\sum_{j\in\Zd} y_j
\phi(\|x-j\|)\right|^2 \,dx \hbox{\hskip 2in}
$$
$$\eqalign{
       &= \lim_{n\to\infty} (2\pi)^{-d} \int_\Td K_n(\xi-\pi e)
              (d^{-1} \sum_{j=1}^d \sin^2(\xi_j/2) )^{2m} \sigma(\xi) \,d\xi \cr
       &= \sigma(\pi e),
}$$
using the fact that $\sigma$ is continuous at $\pi e$ and standard properties of the
Fej\'er kernel.\qed

\sect {\chapno6. Scaling the infinite grid}
Here we consider the behaviour of the norm estimate given above when
we scale the infinite regular grid.

\proclaim {Proposition \chapno6.1}. Let $r$ be a positive number and let $(a_j)_{j\in\Zd}$ be a real
sequence of finite support. Then
$$
 \sum_{j,k\in\Zd} a_j a_k \exp(-\|r j - r k \|^2)
  =  (2\pi)^{-d} \int_\Td \Bigl| \sum_{j\in\Zd} a_j \exp ij\xi \Bigr|^2
        E_r^{(d)}(\xi)\,d\xi, 
\eqno{(\chapno6.1)}$$
where
$$ E_r^{(d)}(\xi) = \sum_{k\in\Zd} e^{-\|rk\|^2} \exp ik\xi, 
\qquad \xi \in \Rd. \eqno{(\chapno6.2)}$$

\pf Section 4.2 provides the equation
$$
 \sum_{j,k\in\Zd} a_j a_k \exp(-\|r j - r k \|^2) \hbox{\hskip 2 truein}$$
$$ 
  =  (2\pi)^{-d} \int_\Td \Bigl| \sum_{j\in\Zd} a_j \exp ij\xi \Bigr|^2
        (\pi/r^2)^{d/2} \sum_{k\in\Zd} \exp(-\|\xi+2\pi k\|^2/4r^2)\,d\xi. 
\eqno{(\chapno6.3)}$$
Further, the Poisson summation formula gives the relation
$$
 (2\pi)^d \sum_{k\in\Zd} \exp(-\|\xi+2\pi k\|^2/4r^2)
   = (4\pi r^2)^{d/2} \sum_{k\in\Zd} e^{-\|rk\|^2} \exp ik\xi.
\eqno{(\chapno6.4)}
$$ 
Substituting (\chapno6.4) into (\chapno6.3) yields equations
(\chapno6.1) and (\chapno6.2).~\qed

The functions $E_r^{(1)}$ and $E_r^{(d)}$ are related in a simple way.

\proclaim {Lemma \chapno6.2}. We have the expression 
$$ E_r^{(d)}(\xi) = \prod_{k=1}^d E_r^{(1)}(\xi_k). \eqno{(\chapno6.5)}$$

\pf This is a straightforward consequence of (\chapno6.2).~\qed

Applying the theta function formulae of Section 4.2 yields the
following result.

\proclaim {Lemma \chapno6.3}. 
$$ E_r^{(1)}(\xi) = 
  \prod_{k=1}^\infty (1- e^{-2kr^2})(1 + 2 e^{-(2k-1)r^2}\cos \xi
+ e^{-(4k-2)r^2}). \eqno{(\chapno6.6)}
$$

\pf The Theta function $\theta_3$ of Jacobi type is given by
$$\eqalign{
  \theta_3(z,q) 
  &= 1 + 2 \sum_{k=1}^\infty q^{k^2} \cos 2kz, \qquad q,z \in \CC, |q|
< 1,\cr
  &= \prod_{k=1}^\infty (1-q^{2k})(1 + 2q^{2k-1}\cos 2z + q^{4k-2}),
}\eqno{(\chapno6.7)}$$
which equations are discussed in greater detail in Section 4.2.
Setting $q = e^{-r^2}$ we have the expressions
$$E_r^{(1)}(\xi) = \theta_3(\xi/2,q)
   = \prod_{k=1}^\infty (1- e^{-2kr^2})(1 + 2 e^{-(2k-1)r^2}\cos \xi
+ e^{-(4k-2)r^2}).
\eqno{(\chapno6.8)}$$
The proof is complete.~\qed

Now Section 4.3 provides the inequality
$$ \sum_{j,k\in\Zd} a_j a_k \exp(-\|rj-rk\|^2)
 \ge E_r^{(d)}(\pi e) \sum_{j\in\Zd} a_j^2, \eqno{(\chapno6.9)}$$
where $e = [1, \ldots, 1]^T \in \Rd$. Using equation (\chapno6.6), we see that
$$\eqalign{
 E_r^{(1)}(\pi) 
  &= \prod_{k=1}^\infty (1 - e^{-2kr^2})(1 - 2 e^{-(2k-1)r^2} + e^{-(4k-2)r^2}) \cr
  &= \prod_{k=1}^\infty (1 - e^{-2kr^2})(1 - e^{-(2k-1)r^2})^2,
}\eqno{(\chapno6.10)}$$
which implies that $\{E_r^{(\pi)} : r > 0 \}$ is an increasing function.
Further, it is a consequence of (\chapno6.5) that $\{ E_r^{(d)}(\pi e): r > 0 \}$ is also an increasing
function. We state these results formally.

\proclaim {Theorem \chapno6.4}. Let $r > s > 0$. Then we have the inequality
$$
\inf \sum_{j,k\in\Zd} a_j a_k \exp(-\|rj-rk\|^2)
  \ge \inf \sum_{j,k\in\Zd} a_j a_k \exp(-\|sj-sk\|^2),
\eqno{(\chapno6.11)}$$
where the infima are taken over the set of real sequences of finite support.

In fact we extend the given analysis to a class of functions including
the multiquadric. The appropriate definitions and theorems form
Section 4.3, but the key result is Theorem 4.3.6: Under
suitable conditions, the function $\phi \colon [0,\infty) \to \R1$ possesses
the generalized Fourier transform
$$\phi(\|\xi\|) = 
 - \int_0^\infty \exp(-\|\xi\|^2/4t) t^{-1}\,d\mu(t),
\eqno{(\chapno6.11)}$$
where
$$ \phi(r) = \phi(0) 
  + \int_0^\infty (1-e^{-r^2t}) t^{-1}\,d\mu(t), \eqno{(\chapno6.12)}$$
and $\mu$ is a positive Borel measure such that $\int_0^1 d\mu(t) <
\infty$ and $\int_1^\infty t^{-1}\,d\mu(t) < \infty$. Now the function
$\phi_r \colon x \mapsto \phi(\|rx\|)$ has the Fourier transform
$$ \hat\phi_r(\|\xi\|) = \phihat(\|\xi\|/r) r^{-d}. \eqno{(\chapno6.13)}$$
Further, the associated symbol function is defined by the equation
$$ \sigma_r(\xi) 
  = \sum_{k\in\Zd} |\hat\phi_r(\|\xi+2\pi k\|)|, \eqno{(\chapno6.14)}$$
and so (\chapno6.13) implies the expression
$$ \sigma_r(\xi)
 = \int_0^\infty r^{-d} \sum_{k\in\Zd} 
     \exp(-\|\xi+2\pi k\|^2/4tr^2) (\pi/t)^{d/2} t^{-1}\,d\mu(t).
\eqno{(\chapno6.15)}$$
Using the Poisson summation formula, we have
$$ (2\pi)^d \sum_{k\in\Zd}  \exp(-\|\xi+2\pi k\|^2/4tr^2)
  = (4tr^2 \pi)^{d/2} \sum_{k\in\Zd} e^{-\|k\|^2tr^2} e^{ik\xi}.
\eqno{(\chapno6.16)}$$
Consequently we have 
$$ r^{-d}(\pi/t)^{d/2} \sum_{k\in\Zd} \exp(-\|\xi+2\pi k\|^2/4tr^2)
 = \sum_{k\in\Zd} e^{-\|k\|^2r^2t} e^{ik\xi} =
E_{rt^{1/2}}^{(d)}(\xi),
\eqno{(\chapno6.18)}$$
providing the equation
$$ \sigma_r(\pi e) = \int_0^\infty E_{r t^{1/2}}^{(d)}(\pi e) t^{-1}
\,d\mu(t),
\eqno{(\chapno6.18)}$$
and so $\{\sigma_r(\pi e) : r > 0 \}$ is an increasing function.

\sect {Appendix}
I do not like stating integral representations such as Theorem
\chapno3.5 without including some explicit examples. Therefore this
appendix calculates $d\alpha$ for $\phi(r) = r$ and $\phi(r)=(r^2 +
c^2)^{1/2}$, where $c$ is positive and we are using the notation of
4.3.5.

For $\phi(r) = r$ the key integral is
$$
 \Gamma(-\half) = \int_0^\infty \left( e^{-u} - 1 \right)
u^{-3/2}\,du, \eqno{(A1)}$$
which is derived in Whittaker and Watson (1927), Section 12.21. Making
the substitution $u = r^2 t$ in (A1) and using the equations
$\pi^{1/2} = \Gamma(1/2) = -\Gamma(-1/2)/2$ we have
$$
- (4\pi)^{1/2} = r^{-1} \int_0^\infty \left(e^{-r^2 t} - 1 \right)
t^{-3/2}\,dt, $$
that is 
$$
 r = \int_0^\infty \left( 1 - e^{-r^2 t} \right) t^{-1}\ (4\pi
t)^{-1/2}\,dt. \eqno{(A2)}$$
Thus the Borel measure is $d\alpha_1(t) = (4\pi t)^{-1/2}\,dt$ and
$\int_0^1 d\alpha_1(t) = \int_1^\infty t^{-1}\,d\alpha_1(t) =
\pi^{1/2}$.

The representation for the multiquadric is an easy consequence of
(A2). Substituting $(r^2 + c^2)^{1/2}$ and $c$ for $r$ in (A2) we
obtain
$$
 (r^2 + c^2)^{1/2} = \int_0^\infty \left(1 - e^{-(r^2+c^2)t} \right)
t^{-1}\ (4\pi t)^{-1/2}\,dt \eqno{(A3)}$$
and 
$$
 c = \int_0^\infty \left(1 - e^{-c^2 t} \right)
t^{-1}\ (4\pi t)^{-1/2}\,dt ,\eqno{(A4)}$$
respectively. Subtracting (A4) from (A3) provides the formula
$$ (r^2 + c^2)^{1/2} = c + \int_0^\infty \left(1 - e^{-r^2 t} \right)
t^{-1}\ e^{-c^2 t} (4\pi t)^{-1/2}\,dt. \eqno{(A5)}$$
Hence the measure is $d\alpha_2(t) = e^{-c^2 t} (4\pi t)^{-1/2}\,dt$.

\vfill \eject
%
\def\chapno{5.}

\headline{\ifnum\pageno=70\hfil\else{\headfont\hfil\chptitle\hfil}\fi}
\def\chptitle{Norm estimates for Toeplitz distance matrices II}
\bigskip
\centerline{\bf 5 : Norm estimates for Toeplitz distance matrices II}



\def\pts{{(x_j)_{j=1}^n}}

\def\sjk{{\sum_{j,k=1}^n y_j y_k\,}}
\def\djk{{\| x_j - x_k \|}}
\def\xjk{{x_j-x_k}}
\def\yjk{{|\sum_{j=1}^n y_j e^{ix_j t}|^2}}
\def\trig{{|\sum_{j=1}^n y_j e^{ijt}|^2}}

\def\skzd{{\sum_{k\in\Zd}}}

\def\rint{{\int_{-\infty}^{\infty}}}
\def\hint{{\int_0^\infty}}

\def\Sdm1{{S^{d-1}}}

\def\ysqr{{\Vert y \Vert^2}}

\def\ainv{{A_n^{-1}}}
\def\aninv{{\Vert A_n^{-1} \Vert_2}}

\def\zsum{{\sum_{k=-\infty}^\infty}}

\def\hi{{(2\pi)^{-1}}}
\def\tint{{\hi\int_0^{2\pi}}}
\def\tintd{{(2\pi)^{-d}\int_{\Td}}}
\def\fhat{{\hat f}}
\def\ajxi{{\Bigl|\sum_{j\in\Zd} a_j \exp(ij\xi)\Bigl|^2}}
\def\Rdm0{{\Rd\setminus\{0\}}}
\def\ghat{{\hat g}}
\def\ytrigo{{\Bigl|\sum_{j\in\Zd} y_j \exp(ix_j \xi)\Bigl|^2}}
\def\ytrig2{{\Bigl|\sum_{j\in\Zd} y_j \exp(ij \xi)\Bigl|^2}}
\def\ynj{{y_j^{(n)}}}

\def\phichat{{\phihat_c}}
\def\chichat{{{\hat \chi}_c}}

\def\min{{\rm min}}

\font\sc=cmcsc10

\font\headfont=cmsl9

\def\pf{{\noindent{\it Proof.\ }}}
\def\qed{{ \vrule height7pt width7pt depth0pt}\par\bigskip}

\def\Rd{{\R1^d}}
\def\Rn{{\R1^n}}
\def\Z0{{\sum_{j=1}^n y_j = 0}}
\def\Zd{{\ZZ^d}}
\def\T1{{[0,2\pi]}}
\def\Td{{[0,2\pi]^d}}
\def\CC{{\cal C}}

\def\half{{1\over2}}
\def\OO{{\cal O}}
\def\phi{{\varphi}}
\def\phihat{{\hat\varphi}}

\outer\def\beginsection#1\par{\vskip0pt plus.3\vsize\penalty-150
 \vskip0pt plus-.3\vsize\bigskip\vskip\parskip
 \message{#1}\centerline{\bf#1}\nobreak\smallskip\noindent}

\def\proclaim #1. #2\par{\medbreak
\noindent{\bf#1.\enspace}{\it#2}\par\medbreak}

\def\proclaimdefn #1. #2\par{\medbreak
\noindent{\bf#1.\enspace}#2\par\medbreak}



\def\sect#1{\goodbreak\bigskip\smallskip\centerline{\bf #1}
\medskip\noindent\ignorespaces}

\def\nullx{\hfill}



\def\authors{B. J. C. Baxter and C. A. Micchelli}

\def\fhat{{\hat f}}
\def\phic{{\phi_c}}
\def\phichat{{\phihat_c}}
\def\chic{{\chi_c}}
\def\chichat{{{\hat \chi}_c}}
\def\vol{\hbox{vol}}
\def\Cd{\CC^d}

\sect {\chapno1. Introduction}
Let $\phi \colon \Rd \to \R1$ be an even continuous function of at
most polynomial growth. Associated with this function is a symmetric
bi-infinite multivariate Toeplitz matrix
$$
 \Phi = \left( \phi(j-k) \right)_{j,k\in\Zd}.
\eqno{(\chapno1.1)} $$
Every finite subset $I = (i_j)_{j=1}^n$ of $\Zd$ determines a finite submatrix of
$\Phi$ given by
$$
 \Phi_I := \left( \phi(i_j-_kk) \right)_{j,k=1}^n.
\eqno{(\chapno1.2)}$$
We are interested in upper bounds on the $\ell^2$-norm of the inverse
matrix $\Phi^{-1}$, that is the quantity
$$
 \|\Phi_I^{-1}\| = 1 \Bigl/ \min\{ \|x\|_2 : 
   \|\Phi_I x\|_2 = 1, \quad x \in \R1^I \}, 
\eqno{(\chapno1.3)}$$
where $\|x\|_2^2 = \sum_{j\in I} x_j^2$ for $x = (x_j)_{j\in I}$. The type of bound we seek
follows the pattern of results in the previous chapter. Specifically, we let
$\phihat$ be the distributional Fourier transform of $\phi$ in the
sense of Schwartz (1966), which we assume to be a measurable function
on $\Rd$. We let $e := (1, \ldots, 1)^T \in \Rd$ and set 
$$
 \tau_\phihat := \sum_{j\in\Zd} |\phihat(\pi e + 2\pi j)|
\eqno{(\chapno1.4)}$$
whenever the right hand side of this equation is meaningful. Then, for
a certain class of {\it radially symmetric} functions, we proved in
Chapter 4 that
$$
 \|\Phi_I^{-1} \| \le 1 / \tau_\phihat
\eqno{(\chapno1.5)}$$
for every finite subset $I$ of $\Zd$. Here we extend this
bound to a wider class of functions which need not be radially
symmetric. For instance, we show that (\chapno1.5) holds for the class of
functions
$$
 \phi(x) = ( \|x\|_1 + c)^\gamma, \qquad x\in\Rd,
$$
where $\|x\|_1 = \sum_{j=1}^d |x_j|$ is the $\ell_1$-norm of $x$, $c$
is non-negative, and $0 < \gamma < 1$.

Our analysis develops the methods of Chapter 4. However, here
we emphasize the importance of certain properties of P\'olya frequency
functions and P\'olya frequency sequences (due to I.~J. Schoenberg) in
order to obtain estimates like (\chapno1.5).

In Section 2 we consider Fourier
transform techniques which we need to prove our bound. Further, the
results of this section improve on the treatment of the last chapter, in
that the condition of admissibility (see Definition \chapno3.2) is shown to be unnecessary. Section 3 contains a
discussion of the class of functions $\phi$ for which we will prove
the bound (\chapno1.4). The final section contains the proof of our main
result.

\def\Fhat{\hat F}
\def\intrd{{(2\pi)^{-d} \int_\Rd}}
\def\inttd{{(2\pi)^{-d} \int_\Td}}
\def\trig{{\Bigl| \sum_{j\in\Zd} y_j e^{ij\xi} \Bigr|^2}}

\sect {\chapno2. Preliminary facts}
We begin with a rather general framework. Let $\phi \colon \Rd \to
\R1$ be a continuous function of polynomial growth. Thus $\phi$
possesses a distributional Fourier transform in the sense of Schwartz
(1966). We shall assume $\phihat$ is almost everywhere equal to a Lebesgue measurable
function on $\Rd$, that is we assume $\phihat$ to be the sum of a
measurable function and a tempered distribution whose support is a
set of Lebesgue measure zero. Given a nonzero real sequence $(y_j)_{j\in\Zd}$ of finite
support and points $(x^j)_{j\in\Zd}$ in $\Rd$, we introduce the
function $F \colon \Rd \to \R1$ given by
$$
 F(x) = \sum_{j,k\in\Zd} y_j y_k \phi(x + x^j-x^k), \qquad x \in \Rd.
\eqno{(\chapno2.1)}$$
Thus 
$$
 F(0) = \sum_{j,k\in\Zd} y_j y_k \phi(x^j - x^k),
\eqno{(\chapno2.2)}$$
which is the quadratic form whose study is the object of much of this dissertation.
We observe that the Fourier transform of $F$ is the
tempered distribution
$$ \Fhat(\xi) = \Bigl| \sum_{j\in\Zd} y_j e^{ix^j \xi} \Bigr|^2
\phihat(\xi), \qquad \xi \in \Rd.
\eqno{(\chapno2.3)}$$
Further, if $\Fhat$ is an absolutely integrable function, then we have
the equation
$$
 F(0) = \intrd \Fhat(\xi)\,d\xi,
\eqno{(\chapno2.4)}$$
since $F$ is the inverse distributional Fourier transform of $\Fhat$
and this coincides with the classical inverse transform when $\Fhat
\in L^1(\Rd)$. In  other words, we have the equation
$$
 \sum_{j,k\in\Zd} y_j y_k \phi(x^j - x^k) = \intrd
   \Bigl| \sum_{j\in\Zd} y_j e^{ix^j \xi} \Bigr|^2 \phihat(\xi)\,d\xi
\eqno{(\chapno2.5)}$$
when $\Fhat$ is absolutely integrable. If we make the further
assumption that $\phihat$ is {\it one-signed} almost everywhere on $\Rd$, and the points
$(x^j)_{j\in\Zd}$ 
form a subset of the integers $\Zd$, then it is possible to
improve (\chapno2.5). First observe that dissecting $\Rd$ into $2\pi$-integer
translates of the cube $\Td$ provides the relations
$$\eqalign{
 \sum_{j,k\in\Zd} y_j y_k \phi(j - k)
   &= \intrd \trig \phihat(\xi)\,d\xi \cr
   &= \sum_{k\in\Zd} \inttd \trig \phihat(\xi+2\pi k)\,d\xi \cr
   &= \inttd \trig \sigma(\xi)\,d\xi,
}\eqno{(\chapno2.6)}$$
where 
$$
 \sigma(\xi) = \sum_{k\in\Zd} \phihat(\xi + 2\pi k), \qquad \xi \in
\Rd ,
\eqno{(\chapno2.7)}$$
and the monotone convergence theorem justifies the exchange of
summation and integration. Further, we see that another consequence of
the condition that $\phihat$ be one-signed is the bound
$$ \trig |\phihat(\xi)| < \infty $$
for almost every point $\xi \in \Td$, because the left hand side of
(\chapno2.6) is {\it a fortiori} finite. This implies that $\sigma$ is almost
everywhere finite, since the set of all zeros of a nonzero trigonometric polynomial 
has measure zero. This last result is well-known, but we include its
short proof for completeness. Following Rudin (1973), we shall say that a
continuous function $f \colon \Cd \to \CC$ is an {\it entire
function of $d$ complex variables} if, for every point $(w_1, \ldots,
w_d) \in \Cd$ and for every $j \in \{ 1, \ldots, d \}$, the mapping
$$ \CC \ni z \mapsto f(w_1, \ldots, w_{j-1}, z, w_{j+1}, \ldots, w_d)$$
is an entire function of one complex variable.

\proclaim {Lemma \chapno2.1}. Given complex numbers $(a_j)_{j=1}^n$ and a set of
distinct points $(x^j)_{1}^n$ in $\Rd$,
we let $p \colon \Rd \to \CC$ be the
function
$$ p(\xi) = \sum_{j=1}^n a_j e^{ix^j \xi}, \qquad \xi \in \Rd.$$
Then $p$ enjoys the following properties:
\item{(i)} $p$ is
identically zero if and only if $a_j = 0$, $1 \le j
\le n$.
\item{(ii)} $p$ is nonzero almost everywhere unless $a_j=0$, $1 \le j
\le n$.

\pf

\def\vol{\hbox{vol}}
\noindent(i) Suppose $p$ is identically zero.
Choose any $j \in \{ 1, \ldots, n\}$ and let $f_j \colon \Rd \to
\R1$ be a continuous function of compact support such that $f_j(x^k) =
\delta_{jk}$ for $1 \le k \le n$. Then
$$ a_j = \sum_{k=1}^n a_k f_j(x^k) = \intrd \sum_{k=1}^n a_k
e^{ix^k\xi} {\hat f_j}(\xi)\,d\xi = 0.$$
The converse is obvious.

\noindent(ii) Let $f \colon \Cd \to \CC$ be an entire
function and let
$$ Z = \{ x\in \Rd : f(x) = 0 \}. $$
If $\vol_d Z$ is a set of positive Lebesgue measure in $\Rd$, then we
shall prove that $f$ is
identically zero, which implies the required result.

We proceed by induction on the dimension $d$. If $d=1$ and $\vol_1 Z >
0$, then $f$ is an entire function of one complex variable with
uncountably many zeros. Such a function must vanish everywhere,
because every uncountable subset of $\CC$ possesses a limit point. Now
suppose that the result is true for $d-1$ for some $d \ge 2$. Fubini's
theorem provides the relation
$$ 0 < \vol_d Z = \int_{\R1^{d-1}} \vol_1 Z(x_2, \ldots, x_d) \,dx_2
\ldots dx_d,$$
where 
$$ Z(x_2, \ldots, x_d) = \{ x_1 \in \R1 : (x_1, \ldots, x_d) \in Z
\}.$$
Thus there is a set, $X$ say, in $\R1^{d-1}$ of positive
$(d-1)$-dimensional Lebesgue measure such that $\vol_1 Z(x_2, \ldots,
x_d)$ is positive for every $(x_2, \ldots, x_d) \in X$, and therefore
the entire function $\CC \ni z \mapsto f(z,x_2, \ldots, x_d)$ vanishes
for all $z \in \CC$, because $Z(x_2, \ldots, x_d)$ is an uncountable
set. Thus, choosing any $z_1 \in \CC$, we see that the entire function
of $d-1$ complex variables defined by 
$$ (z_2, \ldots, z_d) \mapsto f(z_1, z_2, \ldots, z_d), \qquad (z_2,
\ldots, z_d) \in \CC^{d-1}, $$
vanishes for all $(z_2, \ldots, z_d)$ in $X$, which is a set of
positive $(d-1)$- dimensional Lebesgue measure. By induction
hypothesis, we deduce that 
$$ f(z_1, z_2, \ldots, z_d) = 0 \hbox{ for all } z_2, \ldots, z_d \in
\CC,$$
and since $z_1$ can be any complex number, we conclude that $f$ is
identically zero. Therefore the lemma is true.
~\qed

\def\esssup{\hbox{ess sup }}
\def\essinf{\hbox{ess inf }}

We can now derive our first bounds on the quadratic form (\chapno2.2). For
any measurable function $g \colon \Td \to \R1$, we recall the
definitions of the essential supremum
$$
 \esssup g = \inf \{ c \in \R1 : g(x) \le c \hbox{ for almost every }
x \in \Td \}
\eqno{(\chapno2.8)}$$
and the essential infimum
$$
 \essinf g = \sup \{ c \in \R1 : g(x) \ge c \hbox{ for almost every }
x \in \Td \}.
\eqno{(\chapno2.9)}$$
Thus (\chapno2.6) and the Parseval relation provide the inequalities
$$
 \essinf \sigma \sum_{j\in\Zd} y_j^2
   \le \sum_{j,k\in\Zd} y_j y_k \phi(x^j-x^k) 
   \le \esssup \sigma \sum_{j\in\Zd} y_j^2.
\eqno{(\chapno2.10)}$$
Let $V$ be the vector space of real sequences $(y_j)_{j\in\Zd}$ of
finite support for which the function $\Fhat$ of (\chapno2.3) is
absolutely integrable.  We have seen that (\chapno2.10) is valid for every
element $(y_j)_{j\in\Zd}$ of $V$. Of course, at this stage there is
no guarantee that $V \ne \{0\}$ or that the bounds are
finite. Nevertheless, we identify below a case when the bounds (\chapno2.10)
cannot be improved. This will be of relevance later.

\def\II{{\cal I}}

\proclaim {Proposition \chapno2.2}. Let $P$ be a nonzero trigonometric
polynomial such that the principal ideal $\II$ generated by $P$, that is
the set
$$
 \II = \{ P T : T \hbox{ a real trigonometric polynomial} \},
\eqno{(\chapno2.11)}$$
consists of trigonometric polynomials whose Fourier coefficient
sequences are elements of $V$. Further, suppose that there is a
point $\eta$ at which $\sigma$ is continuous and $P(\eta) \ne 0$. Then
we can find a sequence $\{ (y_j^{(n)})_{j\in\Zd} : n = 1, 2, \ldots
\}$ in $V$ such that
$$
 \lim_{n\to\infty} \sum_{j,k\in\Zd} y_j^{(n)} y_k^{(n)} \phi(j-k)
 \Bigl/ \sum_{j\in\Zd} [y_j^{(n)}]^2 = \sigma(\eta).
\eqno{(\chapno2.12)}$$

\pf We recal Section 4 and recall that the $n$\/th degree tensor
product Fej\'er kernel is defined by
$$
 K_n(\xi) := \prod_{j=1}^d {\sin^2 n\xi_j/2 \over n \sin^2 \xi_j/2}
           = \Bigl| n^{-d/2} \sum_{k\in\Zd \atop 0 \le k < e n} e^{ik\xi}
\Bigr|^2 =: |L_n(\xi)|^2, \qquad \xi \in \Rd,
\eqno{(\chapno2.13)}$$
where $e = (1, \ldots, 1)^T \in \Rd$ and $L_n(\xi) = n^{-d/2} \sum_{0
\le k < en} e^{ik\xi}$. Then the function $P(\cdot)
L_n(\cdot-\eta)$ is a member of $\II$ and we choose
$(y_j^{(n)})_{j\in\Zd}$ to be its Fourier coefficient sequence. The
Parseval relation provides the equation
$$
  \sum_{j\in\Zd} [y_j^{(n)}]^2 = \inttd P^2(\xi)K_n(\xi-\eta)\,d\xi
\eqno{(\chapno2.14)}$$
and the approximate identity property of the Fej\'er kernel (Zygmund
(1988), p.86) implies that
$$\eqalign{
 P^2(\eta) 
  &= \lim_{n\to\infty} \inttd P^2(\xi) K_n(\xi-\eta)\,d\xi \cr
  &= \lim_{n\to\infty} \sum_{j\in\Zd} [y_j^{(n)}]^2.
}\eqno{(\chapno2.15)}$$
Further, because $\sigma$ is continuous at $\eta$, we also have
the relations
$$\eqalign{
 P^2(\eta)\sigma(\eta)
  &= \lim_{n\to\infty} \inttd
      P^2(\xi)K_n(\xi-\eta)\sigma(\xi)\,d\xi \cr
  &= \lim_{n\to\infty} \sum_{j,k\in\Zd} y_j^{(n)} y_k^{(n)} \phi(j-k),
}\eqno{(\chapno2.16)}$$
the last line being a consequence of (\chapno2.6).
Hence (\chapno2.15) and (\chapno2.16) provide equation (\chapno2.12).~\qed

\proclaim {Corollary \chapno2.3}. If $\sigma$ attains its essential infimum
(resp. supremum) at a point of continuity, and if we can find a
trigonometric polynomial $P$ satisfying the conditions of Proposition
\chapno2.2, then the lower (resp. upper) bound of
(\chapno2.10) cannot be improved.

\pf This is an obvious consequence of Proposition \chapno2.2.~\qed
 
\def\Ghat{\hat G}
\def\AA{{\cal A}}
\def\GG{{\cal G}}

We now specialize this general setting to the following case.

\proclaimdefn {Definition \chapno2.4}. Let $G \colon \Rd \to \R1$ be a
continuous absolutely integrable function such that $G(0) = 1$ for
which the Fourier
transform is non-negative and absolutely integrable. Further, we
require that there exist
non-negative constants $C$ and $\kappa$ for which
$$
 | 1 - G(x) | \le C \|x\|^\kappa, \qquad x \in \Rd.
\eqno{(\chapno2.17)}$$
We let $\GG$ denote the class of all such functions $G$.

Clearly the Gaussian $G(x) = \exp(-\|x\|^2)$ provides an
example of such a function. The next lemma mentions some salient
properties of $\GG$ which do not, however, require (\chapno2.17).

\proclaim {Lemma \chapno2.5}. Let $G \in \GG$.

{\it
\noindent\item{(i)} $G$ is a symmetric function, that is
$$ G(x) = G(-x), \qquad x \in \Rd. \eqno{(\chapno2.18)}$$

\noindent\item{(ii)} $$|G(x)| \le 1, \qquad x \in \Rd. \eqno{(\chapno2.19)}$$

\noindent\item{(iii)} $G$ is a positive definite function in the sense
of Bochner. In other words, for any real sequence $(y_j)_{j\in\Zd}$ of
finite support, and for any points $(x^j)_{j\in\Zd}$ in $\Rd$, we have
$$ \sum_{j,k\in\Zd} y_j y_k G(x^j - x^k) \ge 0. \eqno{(\chapno2.20)}$$
}

\pf

\noindent\item{(i)} Since $\Ghat$ is real-valued we have
$$ 2i \int_\Rd G(x) \sin x\xi\,dx = \Ghat(\xi) - \Ghat(-\xi) \in \R1,
\qquad \xi \in \Rd, $$
which is a contradiction unless both sides vanish. Thus $\Ghat$ is a
symmetric function. However, $G$ must inherit this symmetry, by the
Fourier inversion theorem.

\noindent\item{(ii)} The non-negativity of $\Ghat$ provides the
relations
$$ |G(x)| 
   = \Bigl| \intrd \Ghat(\xi) e^{-ix\xi}\,d\xi\Bigr|
   \le \intrd \Ghat(\xi)\,d\xi = G(0) = 1.$$

\noindent\item{(iii)} The condition $\Ghat \in L^1(\Rd)$ implies the
validity of (\chapno2.5) for $\phi$ replaced by $G$, whence
$$
 \sum_{j,k\in\Zd} y_j y_k G(x^j-x^k) = 
   \intrd \Bigl|\sum_{j\in\Zd} y_j e^{ix^j\xi}\Bigr|^2 \Ghat(\xi)\,d\xi
\ge 0, $$
as required.~\qed

We remark that the first two parts of Lemma \chapno2.5 are usually deduced from the
requirement that $G$ be a positive definite function in the Bochner
sense
(see Katznelson (1976), p.137). We have presented our material in this order
because it is the non-negativity condition on $\Ghat$ which forms our
starting point.

Given any $G \in \GG$, we define the set $\AA(G)$ of functions of the form
$$
  \phi(x) = c 
   + \int_0^\infty [ 1 - G(t^{1/2}x)] t^{-1}\,d\alpha(t), \qquad x \in
\Rd,
\eqno{(\chapno2.21)}$$
where $c$ is a constant and $\alpha \colon [0,\infty) \to \R1$ is a non-decreasing function
such that 
$$\int_1^\infty t^{-1}\,d\alpha(t) < \infty \hbox{ and } \int_0^1
t^{\kappa/2-1}\,d\alpha(t) < \infty. \eqno{(\chapno2.22)}$$

Let us show that (\chapno2.21) is well-defined. Inequality
(\chapno2.19) implies the bound
$$
 \int_1^\infty \Bigl| 1 - G(t^{1/2}x) \Bigr| t^{-1}\,d\alpha(t) 
 \le 2 \int_1^\infty t^{-1}\,d\alpha(t) < \infty. 
$$
Moreover, applying condition (\chapno2.17) we obtain
$$
 \int_0^1 \Bigl| 1 - G(t^{1/2}x)\Bigr| t^{-1}\,d\alpha(t) 
   \le C \|x\|^\kappa \int_0^1 t^{\kappa/2-1}\,d\alpha(t) < \infty. 
\eqno{(\chapno2.23)}$$
Therefore the integral of (\chapno2.21)
is finite and $\phi$ is a function of polynomial growth. A
simple application of the dominated convergence theorem
reveals that $\phi$ is also continuous, so that we may view it as a
tempered distribution.

The following definition is convenient.

\proclaimdefn {Definition \chapno2.6}. We shall say that a real sequence
$(y_j)_{j\in\Zd}$ of finite support is {\it zero-summing} if
$\sum_{j\in\Zd} y_j = 0$.

An important property of $\AA(G)$ is that it consists of {\it
conditionally negative definite functions of order 1 on $\Rd$}, that
is whenever $\phi \in \AA(G)$
$$ \sum_{j,k\in\Zd} y_j y_k \phi(x^j - x^k) \le 0 \eqno{(\chapno2.24)}$$
for every zero-summing sequence $(y_j)_{j\in\Zd}$ and for any points
$(x_j)_{j\in\Zd}$ in $\Rd$. Indeed, (\chapno2.21) provides the equation
$$ \sum_{j,k\in\Zd} y_j y_k \phi(x^j-x^k)
  = -\int_0^\infty \sum_{j,k\in\Zd} y_j y_k G(t^{1/2}(x^j-x^k)) \,
        t^{-1}\,d\alpha(t) ,
\eqno{(\chapno2.25)}$$
and the right hand side is non-positive because $G$ is positive
definite in the Bochner sense (Lemma \chapno2.5 (iii)).

We now fix attention on a particular element $G \in \GG$ and a function
$\phi \in \AA(G)$.

\proclaim {Theorem \chapno2.7}. Let $(y_j)_{j\in\Zd}$ be a zero-summing
sequence that is not identically zero. Then, for any points $(x^j)_{j\in\Zd}$ in $\Rd$, we have the
equation
$$
 \sum_{j,k\in\Zd} y_j y_k \phi(x^j-x^k) = -\intrd \Bigl| \sum_{j\in\Zd}
y_j e^{ix^j \xi} \Bigr|^2 H(\xi)\,d\xi,
\eqno{(\chapno2.26)}$$
where
$$
 H(\xi) = \int_0^\infty \Ghat(\xi/t^{1/2}) t^{-d/2-1}\,d\alpha(t),
\qquad \xi \in \Rd. 
\eqno{(\chapno2.27)}$$
Furthermore, this latter integral is finite for almost every $\xi \in \Rd$.

\pf Applying the Fourier inversion theorem to $G$ in (\chapno2.25), we obtain
$$ \sum_{j,k\in\Zd} y_j y_k \phi(x^j-x^k) \hbox{\hskip 7cm} $$
$$\eqalign{
  &= - (2\pi)^{-d} \int_0^\infty \int_\Rd \Bigl|\sum_{j\in\Zd} y_j 
           \exp(it^{1/2}\eta x^j)\Bigr|^2
           \Ghat(\eta)t^{-1}\,d\eta\,d\alpha(t) \cr
  &= - (2\pi)^{-d} \int_0^\infty \int_\Rd \Bigl|\sum_{j\in\Zd} y_j e^{ix^j\xi}\Bigr|^2
         \Ghat(\xi/t^{1/2}) t^{-d/2-1}\, d\xi\, d\alpha(t),
}\eqno{(\chapno2.28)}$$
where we have used the substitution $\xi = t^{1/2} \eta$. Because the 
integrand in the last line is
non-negative, we can exchange the order
of integration to obtain (\chapno2.26). Of course the left hand side
of (\chapno2.26) is finite, which implies that  the integrand of 
(\chapno2.26) is an
absolutely integrable function, and hence finite almost everywhere.
But, by Lemma \chapno2.1, $|\sum_j y_j e^{ix^j \xi} |^2 \ne 0$ for almost
every $\xi \in \Rd$ if the sequence $(y_j)_{j\in\Zd}$ is non-zero.
Therefore $H$ is finite almost everywhere.~\qed

\proclaim {Corollary \chapno2.8}. The hypotheses of Theorem \chapno2.7
imply the equation
$$ F(x) = -\intrd \Bigl|\sum_{j\in\Zd} y_j 
              e^{i x^j\xi} \Bigr|^2 H(\xi) e^{ix\xi}\,d\xi,
\qquad x \in \Rd, \eqno{(\chapno2.29)}$$
where $F$ is given by (\chapno2.1). Consequently, $\phihat(\xi) = - H(\xi)$
for almost every $\xi \in \Rd$, that is
$$ \phihat(\xi) =  - \int_0^\infty \Ghat(\xi/t^{1/2})
t^{-d/2-1}\,d\alpha(t).
\eqno{(\chapno2.30)}$$

\pf It is straightforward to deduce the relation
$$ F(x) = - (2\pi)^{-d} \int_0^\infty \int_\Rd 
               \Bigl|\sum_{j\in\Zd} y_j e^{ix^j\xi}\Bigr|^2 e^{ix\xi}
                 \Ghat(\xi/t^{1/2}) t^{-d/2-1}\, d\xi\, d\alpha(t),
$$            
which is analogous to (\chapno2.28). Now the absolute value of this integrand
is precisely the integrand in the second line of (\chapno2.28). Thus we may
apply Fubini's theorem to exchange the order of integration, obtaining
(\chapno2.29).

\def\psihat{\hat\psi}

Next, we prove that 
$\xi \mapsto - |\sum_j y_j e^{ix^j\xi}|^2 H(\xi)$ is the 
Fourier transform of $F$. Indeed, let $\psi\colon\Rd\to\R1$ be any
smooth function whose partial derivatives enjoy supra-algebraic decay.
It is sufficient (see Rudin (1973)) to show that
$$ \int_\Rd \psihat(x) F(x)\,dx = -\int_\Rd \psi(\xi) 
     \Bigl|\sum_{j\in\Zd} y_j e^{ix^j\xi}\Bigr|^2 H(\xi)\,d\xi.
$$
Applying (\chapno2.29) and Fubini's theorem, we get
$$\eqalign{
  \int_\Rd \psihat(x)F(x)\,dx 
 &= -(2\pi)^{-d} \int_\Rd \int_\Rd \psihat(x)  
        \Bigl|\sum_{j\in\Zd} y_j e^{ix^j\xi}\Bigr|^2 e^{ix\xi}
            H(\xi)\,d\xi\,dx \cr
 &= -\int_\Rd   \Bigl|\sum_{j\in\Zd} y_j e^{ix^j\xi}\Bigr|^2 
            H(\xi) 
      \left((2\pi)^{-d}\int_\Rd \psihat(x)e^{ix\xi}\,dx\right) \,d\xi
\cr
 &= -\int_\Rd  \Bigl|\sum_{j\in\Zd} y_j e^{ix^j\xi}\Bigr|^2 
            H(\xi)\psi(\xi)\,d\xi,
}$$
which establishes (\chapno2.30).
However, we already know that the Fourier transform 
$\Fhat(\xi)$ is almost everywhere
equal to $|\sum_j y_j e^{ix^j\xi}|^2 \phihat(\xi)$. By Lemma \chapno2.1, we
know that $\sum_j y_j e^{ix^j\xi} \ne 0$ for almost all $\xi \in \Rd$,
which implies that $\phihat = -H$ almost everywhere.~\qed

\def\fhat{{\hat f}}
\def\phichat{{\phihat_c}}
\def\chichat{{{\hat \chi}_c}}
\def\Td{{[0,2\pi]^d}}
\def\Lambdahat{{\hat \Lambda}}
\def\CCC{{\cal C}}
\def\EE{{\cal E}}
\def\PP{{\cal P}}

\sect {\chapno3. P\'olya frequency functions}
For every real sequence $(a_j)_{j=1}^\infty$ and any non-negative constant
$\gamma$ such that $0 < \gamma + \sum_{j=1}^\infty a_j^2 < \infty$, we set
$$
  E(z) = e^{-\gamma z^2} \prod_{j=1}^\infty (1 - a_j^2 z^2), \qquad z
\in \CC.
\eqno{(\chapno3.1)}$$
This is an entire function which is nonzero in the vertical strip
$$
 | \Re z | < \rho := 1 / \sup \{ |a_j| : j = 1, 2, \ldots \}.
$$
It can be shown (Karlin (1968), Chapter 5) that there exists a
continuous function $\Lambda \colon \R1 \to \R1$ such that
$$
 \int_\R1 \Lambda(t) e^{-zt}\,dt = {1 \over E(z)}, \qquad |\Re z| < \rho.
\eqno{(\chapno3.2)}$$
This function $\Lambda$ is what Schoenberg (1951) calls a P\'olya
frequency function. We have restricted ourselves to functions
$\Lambda$ which are even, that is
$$
   \Lambda(t) = \Lambda(-t), \qquad t \in \R1.
\eqno{(\chapno3.3)}$$
Also, $E(0)=1$ implies that
$$
 \int_\R1 \Lambda(t)\,dt = 1. 
\eqno{(\chapno3.4)}$$
According to (\chapno3.1) the Fourier transform of $\Lambda$ is given by
$$
 \Lambdahat(\xi) = {1 \over E(i\xi)} 
    = {e^{-\gamma\xi^2} \over \prod_{j=1}^\infty (1 + a_j^2 \xi^2)},
\quad \xi \in \R1 .
\eqno{(\chapno3.5)}$$
We see that $\Lambda(\cdot)/\Lambda(0)$  is a member of the set $\GG$
described in Definition \chapno2.4 for $d=1$, and therefore Lemma \chapno2.5 is applicable. In
particular,
$$ |\Lambda(t)| \le \Lambda(0), \qquad t \in \R1.
\eqno{(\chapno3.6)}$$
However, much more than (\chapno3.6) is true. Schoenberg (1951) proved that
$$
 \det (\Lambda(x_j-y_k))_{j,k=1}^n \ge 0
\eqno{(\chapno3.7)}$$
whenever $x_1 < \cdots < x_n$ and $y_1 < \cdots < y_n$. This fact will
be used in an essential way in Section \chapno4. For the moment we use it to
improve (\chapno3.6) to 
$$
   \Lambda(t) \in [0, \Lambda(0)], \qquad t \in \R1.
\eqno{(\chapno3.8})$$
Let $\PP$ denote the class of functions
$\Lambda \colon \R1 \to \R1$ that satisfy (\chapno3.2)
for some $\gamma \ge 0$ and sequence $(a_j)_{j=1}^\infty$ satisfying
$0 < \gamma + \sum_{j=1}^\infty a_j^2 < \infty$. 
For any positive $a$ the function
$$
 S_a(t) = {1 \over 2|a|} e^{-|t/a|}, \qquad t \in \R1,
\eqno{(\chapno3.9)}$$
is in $\PP$ since 
$$
  \int_\R1 S_a(t) e^{-zt}\,dt = {1 \over 1 - a^2 z^2}, \qquad |\Re z|
< 1/a.
\eqno{(\chapno3.10)}$$
Let $\EE = \{ S_a : a > 0 \}$.
These are the only elements of $\PP$ that are {\it not} in $C^2(\R1)$,
because all other members of $\PP$ have the property that
$\Lambdahat(t) = \OO(t^{-4})$ as $|t| \to \infty$. Hence there exists
a constant $\kappa$ such that
$$
 |\Lambda(0) - \Lambda(t)| \le \kappa t^2, \quad \hbox{ for } t \in
\R1 \hbox{ and } \Lambda \in \PP\setminus\EE,
\eqno{(\chapno3.11)}$$
or
$$
 |\Lambda(0) - \Lambda(t)| \le \kappa |t|, \quad t \in \R1, \quad \Lambda
\in \EE.
\eqno{(\chapno3.12)}$$
We note also that every element of $\PP$ decays exponentially for
large argument (see Karlin (1968), p. 332).

We are now ready to define the multivariate class of functions which
interest us. Choose any $\Lambda_1, \ldots, \Lambda_d \in \PP$ and
define
$$
 G(x) = \prod_{j=1}^d {\Lambda_j(x_j)\over\Lambda_j(0)}, \quad x = (x_1,
\ldots, x_d) \in \Rd.
\eqno{(\chapno3.13)}$$
According to (\chapno3.11) and (\chapno3.12), there is a constant $C \ge 0$ such
that
$$
  1 - G(x)  \le C \|x\|_2^2, \eqno{(\chapno3.14)}
$$
when $\Lambda_j \notin \EE$ for every factor $\Lambda_j$ in
(\chapno3.11). However, if $\Lambda_j \in \EE$ for every $j$, then we
only have
$$
  1 - G(x) \le C \|x\|_2 , \eqno{(\chapno3.15)}
$$
for some constant $C$. We are unable to study the general behaviour at
this time.
Remarking that the Fourier transform of $G$ is given by 
$$
 \Ghat(\xi) = \prod_{j=1}^d 
       {\Lambdahat_j(\xi_j) \over \Lambda_j(0)}, \qquad \xi = (\xi_1,
\ldots, \xi_d) \in \Rd, \eqno{(\chapno3.16)}$$
we conclude that $G$ is a member of the class $\GG$ of Definition \chapno2.4.
Moreover, we can now construct the set $\AA(G)$. To this end,
let $\alpha \colon [0,\infty) \to \R1$ be a non-decreasing function such
that
$$ \int_1^\infty t^{-1}\,d\alpha(t) < \infty, \eqno{(\chapno3.17)}$$
and for any constant $c\in\R1$ define $\phi \colon [0,\infty) \to \R1$ by
(\chapno2.21).
Thus we see that as long as we require the measure $d\alpha$ to
satisfy the extra condition
$$
 \int_0^1 t^{-1/2}\,d\alpha(t) < \infty
\eqno{(\chapno3.18)}$$
whenever one of the factors in (\chapno3.11) is an element of $\EE$, then $\phi$ is
a continuous function of polynomial growth and the results of Section
2 apply. We let $\CCC$ denote the class of all such functions, for all
$G \in \GG$.

Let us note that $\CCC$
contains the following important subclass of functions. In
1938, I.~J. Schoenberg proved that a continuous radially symmetric
function $\phi \colon \Rd \to \R1$ is conditionally negative definite of
order 1 on every $\Rd$ if and only if it has the form
$$ \phi(x) = \phi(0) + 
\int_0^\infty \Bigl(1 - \exp(-t\|x\|^2)\Bigr) t^{-1}\,d\alpha(t), 
\qquad x \in \Rd,$$
where $\alpha \colon [0,\infty) \to \R1$ is a non-decreasing function
satisfying (\chapno3.17). In this case, the Gaussian is clearly of the form (\chapno3.13),
implying that we do indeed
have a subclass of $\CCC$. Thus we have established Theorem \chapno2.7 and Corollary
\chapno2.8 under weaker conditions than those assumed in Chapter 4.

Our class $\CCC$ also contains functions of the form
$$ \phi(x) = c + \int_0^\infty 
      \Bigl( 1 - \exp(-t^{1/2}\|x\|_1) \Bigr) t^{-1}\,d\alpha(t), 
\qquad x \in \Rd,$$
where $\alpha \colon [0,\infty) \to \R1$ is a non-decreasing function
satisfying (\chapno3.17) and (\chapno3.18), and $\|x\|_1 = \sum_{j=1}^d |x_j|$ for
$x = (x_1, \ldots, x_d) \in \Rd$. For instance, using the easily
verified formula
$$ {\gamma \over \Gamma(1+2\gamma)}
     \int_0^\infty \Bigl(1 - e^{-t^{1/2}\sigma}\Bigr) t^{\gamma-1}
e^{-\delta t^{1/2}}\,dt = \delta^{-2\gamma} - (\delta + \sigma)^{-2\gamma}, $$
which is valid for $\delta \ge 0$ and $-1/2 < \gamma < 0$, we see that
$\phi(x) = (\delta + \|x\|_1)^\tau$, for $\delta \ge
0$ and $0 < \tau < 1$, is in our class $\CCC$. 


Although it is not central to our interests in this section, we will
discuss some additional properties of the Fourier transform of a
function $\phi \in \CCC$. First, observe that (\chapno3.5) implies that
$\Lambdahat$ is a decreasing function on $[0,\infty)$ for every $\Lambda$ in $\PP$.
Consequently every $G \in \GG$ satisfies the inequality $\Ghat(\xi)
\le \Ghat(\eta)$ for $\xi \ge \eta \ge 0$. This property is
inherited by the function $H$ of (\chapno2.27), that is
$$ H(\xi) \le H(\eta) \hbox{ whenever } \xi \ge \eta \ge 0,
\eqno{(\chapno3.19)}$$
which allows us to strengthen Theorem \chapno2.7.

\def\Rm0d{(\R1\setminus\{0\})^d}

\proclaim {Proposition \chapno3.1}. $H$ is continuous on
$\Rm0d$.

\pf We first show that $H$ is finite on $\Rm0d$. 
We already know that $\phihat = -H$ almost everywhere, which
implies that every set of positive measure contains a point at which
$H$ is finite. In particular, let
$\delta$ be a positive number and set 
$U_\delta = \{\xi \in \Rd : 0 < \xi_j < \delta,
\quad j = 1, \ldots, d \}$. Thus there is a point $\eta \in U_\delta$
such that $H(\eta) < \infty$. Applying (\chapno3.19) and recalling that $H$
is a symmetric function, we deduce the inequality
$$ H(\xi ) \le H(\eta) < \infty, \qquad \xi \in F_\delta, \eqno{(\chapno3.20)}$$
where $F_\delta := \{\xi \in \Rd: \left|\xi_j\right| \ge \delta, \quad
j = 1, \ldots, d\,\}$. Since $\delta > 0$ is arbitrary, we see that
$H$ is finite in $\Rm0d$.

To prove that $H$ is continuous in $F_\delta$, let
$(\xi_n)_{n=1}^\infty$ be a convergent sequence in $F_\delta$ with
limit $\xi_\infty$. By (\chapno3.20), the functions
$$ \{ t \mapsto \Ghat(\xi_n t^{-1/2}) t^{-d/2-1} : n = 1, 2, \ldots
\}$$
are absolutely integrable on $[0,\infty)$ with respect to the measure
$d\alpha$. Moreover, they are dominated by the $d\alpha$-integrable
function
$t \mapsto \Ghat(\eta t^{-1/2})t^{-d/2 - 1}$. Finally, the continuity
of $\Ghat$ provides the equation
$$ \lim_{n \to \infty} \Ghat(\xi_n t^{-1/2}) t^{-d/2-1} = \Ghat(\xi_\infty
t^{-1/2}) t^{-d/2-1}, \qquad t \in [0,\infty), $$
and thus
$\lim_{n\to\infty} H(\xi_n) = H(\xi_\infty)$ by the dominated
convergence theorem. Since $\delta$ was an arbitrary positive number, 
we conclude that $H$ is continuous on $\Rm0d$.~\qed

The remainder of this section requires a distinction of cases. The
first case (Case I) is the nicest. This occurs when
every factor $\Lambda_j$ in (\chapno3.13) has a positive
exponent $\gamma_j$ in the Fourier transform formula (\chapno3.5). We let
Case II denote the contrary case. Our investigation of Case II is not
yet complete, so we shall concentrate on Case I for the remainder of
this section.

For Case I we have the bound
$$ \Ghat(\xi) \le e^{-\left(\gamma_1 \xi_1^2 + \cdots + \gamma_d
\xi_d^2 \right)}, \qquad \xi \in \Rd, $$
which implies the limit 
$$\lim_{t\to 0} \Ghat(\xi t^{-1/2}) t^{-d/2-1} = 0, \qquad \xi \ne
0.$$
Thus the function $t \mapsto \Ghat(\xi t^{-1/2}) t^{-d/2 - 1}$ is
continuous for $t \in [0,\infty)$ when $\xi$ is nonzero, which implies that
$$\int_0^1 \Ghat(\xi t^{-1/2}) t^{-d/2 - 1}\,d\alpha(t) < \infty ,
\qquad \xi \ne 0.$$
Moreover, since
$$ \int_1^\infty \Ghat(\xi t^{-1/2}) t^{-d/2 - 1}\,d\alpha(t) \le
\int_1^\infty t^{-1}\,d\alpha(t) < \infty, $$
we have $H(\xi) < \infty$ for every $\xi \in \Rdm0$. 
Finally, a simple extension
of the proof of Proposition \chapno3.1 shows that $H$ is 
continuous on $\Rdm0$. 

In fact, we can prove that for $H \in
C^\infty(\Rdm0)$ in Case I. We observe that it is sufficient to show that
every derivative of $\Ghat(\xi t^{-1/2}) t^{-d/2 - 1}$ with respect to
$\xi$ is an absolutely integrable function with respect to the measure
$d\alpha$ on $[0,\infty)$, because then we are justified in
differentiating under the integral sign. Next, the form of $\Ghat$
implies that we only need to show that every derivative of
$\Lambdahat$, where $\Lambdahat$ is given by (\chapno3.13) and $\gamma > 0$,
enjoys faster than algebraic decay for large argument. To this end we
claim that for every $C < \rho := 1/ \sup \{ |a_j| : j = 1, 2, \ldots
\}$ there is a constant $D$ such that
$$ \left| \Lambdahat(\xi + i\eta) \right| \le D e^{-\gamma \xi^2}, \qquad \xi
\in \R1, \quad |\eta| \le C. \eqno{(\chapno3.21)}$$
To verify the claim, observe that when $|\eta| \le C \le |\xi|$ we have the inequalities
$$ \left|e^{-\gamma(\xi+i\eta)^2}\right| \le e^{C^2\gamma}
e^{-\gamma\xi^2} \hbox{ and }
|1 + a_j^2(\xi+i\eta)^2| \ge 1 +
a_j^2(\xi^2-\eta^2) \ge 1. $$ 
Thus, setting $M = \max \{|\Lambdahat(\xi+i\eta)|e^{\gamma\xi^2} :
|\xi|\le C, |\eta| \le C \}$, we conclude that $D := \max\{M,
e^{C^2\gamma}\}$ is suitable in (\chapno3.21).
Finally, we apply the Cauchy
integral formula to estimate the $k$th derivative. We have
$$ \Lambdahat^{(k)}(\xi) = {1 \over 2\pi i} \int_\Gamma
{\Lambdahat(\zeta) \over (\zeta-\xi)^{k+1}}\,d\zeta, $$
where $\Gamma : [0,2\pi] \to \CC$ is given by $\Gamma(t) = r e^{it}$
and $r < C$ is a constant. Consequently we have the bound
$$ \left| \Lambdahat^{(k)}(\xi) \right| \le (D/\alpha^k) e^{-\gamma\,
\min\{(\xi-r)^2, (\xi+r)^2\}}, \qquad \xi \in \R1, $$
and the desired supra-algebraic decay is established. We now state
this formally.

\proclaim {Proposition \chapno3.2}. In Case I, the function $H$ of (\chapno2.27) is
smooth for nonzero argument.

Next, to identify $-H$ with $\phihat$ on $\Rdm0$ in Case I, we let
$\psi \colon \Rd \to \R1$ be a smooth function whose support is a
compact subset of $\Rdm0$. By definition we have
$$ \langle \phihat, \psi \rangle
   = \int_\Rd \psihat(x) \phi(x)\,dx, \eqno{(\chapno3.22)}$$
where $\langle \cdot, \cdot \rangle$ denotes the action of a tempered
distribution on a test function (see Schwartz (1966)). Substituting
the expression for $\phi$ given by (\chapno2.21) into the right hand side of
(\chapno3.22) and using the fact that 
$$0 = \psi(0) = (2\pi)^{-d} \int_\Rd
\psihat(\xi)\,d\xi \eqno{(\chapno3.23)}$$ 
gives
$$ \langle \phihat, \psi \rangle 
= - \int_\Rd \left(\int_0^\infty \psihat(x) (1 - G(t^{1/2}x))
t^{-1}\,d\alpha(t) \right)\,dx. $$
We want to swap the order of integration here. This will be justified
by Fubini's theorem if we can show that 
$$ \int_\Rd \left(\int_0^\infty |\psihat(x)| (1 - G(t^{1/2}x))
t^{-1}\,d\alpha(t) \right)\,dx < \infty. \eqno{(\chapno3.24)}$$
We defer the proof of (\chapno3.24) to Lemma \chapno3.3 below and press on. Swapping
the order of integration and recalling (\chapno3.23) yields
$$\eqalign{
 \langle \phihat, \psi \rangle 
 &= - \int_0^\infty \left( \int_\Rd \psihat(x) G(t^{1/2}x)\,dx \right)
t^{-1}\,d\alpha(t) \cr
 &= - \int_0^\infty \left( \int_\Rd \psi(\xi) \Ghat(\xi
t^{-1/2})\,d\xi \right) t^{-d/2-1}\,d\alpha(t) \cr
}$$
using Parseval's relation in the last line. Once again, we want to
swap the order of integration and, as before, this is justified by
Fubini's theorem if a certain integral is finite, specifically
$$ \int_0^\infty \left( \int_\Rd \left|\psi(\xi)\right| \Ghat(\xi
t^{-1/2})\,d\xi \right) t^{-d/2-1}\,d\alpha(t) < \infty. \eqno{(\chapno3.25)}$$
The proof of (\chapno3.25) will also be found in Lemma \chapno3.3 below. After
swapping the order of integration we have
$$ \langle \phihat, \psi \rangle = - \int_\Rd \psi(\xi)
H(\xi)\,d\xi, \eqno{(\chapno3.26)}$$
which implies that $\phihat = - H$ in $\Rdm0$.

\proclaim {Lemma \chapno3.3}. Inequalities (\chapno3.24) and (\chapno3.25) are valid in
Case I.

\pf For (\chapno3.24), we have
$$
 \int_\Rd \left( \int_0^\infty |\psihat(x)| (1 - G(t^{1/2}x))
t^{-1}\,d\alpha(t) \right)\,dx \hbox{\hskip 2.5truein}$$
$$\eqalign{
 {}
 &\le \int_\Rd \left( \kappa \int_0^1 |\psihat(x)| \|x\|^2\,d\alpha(t)
\right)\,dx + \int_\Rd \left( \int_1^\infty |\psihat(x)|
t^{-1}\,d\alpha(t) \right)\,dx \cr
 &= \kappa (\alpha(1) - \alpha(0)) \int_\Rd |\psihat(x)| \|x\|^2\,dx +
\left(\int_1^\infty t^{-1}],d\alpha(t) \right) \left( \int_\Rd
|\psihat(x)|\,dx \right) \cr
 &< \infty, 
}$$
recalling that $\psihat$ must enjoy faster than algebraic decay
because $\psi$ is a smooth function.

For (\chapno3.25), the substitution $\eta = \xi t^{-1/2}$ provides the
integral
$$ I := \int_0^\infty \left( \int_\Rd |\psi(\eta t^{1/2})|
\Ghat(\eta)\,d\eta \right) t^{-1}\,d\alpha(t). $$
Now there is a constant $D$ such that $|\psi(y)| \le D \|y\|^2$ for
every $y \in \Rd$, because the support of $\psi$ is a closed subset of
$\Rdm0$. Hence 
$$\eqalign{
  I &\le \int_0^1 D \left(\int_\Rd \Ghat(\eta) \|\eta\|^2\,d\eta
\right)\,d\alpha(t) + (2\pi)^d G(0) \|\psi\|_\infty \int_1^\infty
t^{-1}\,d\alpha(t) \cr
    &< \infty.
}$$
The proof is complete.~\qed

\def\trig{\Bigl| \sum_{j\in\Zd} y_j e^{ij\xi} \Bigr|^2}

\sect {\chapno4. Lower bounds on eigenvalues}
Let $\phi\colon\Rd\to\R1$ be a member of $\CCC$ and let
$(y_j)_{j\in\Zd}$ be a zero-summing sequence. An immediate consequence
of (\chapno2.26) is the equation
$$
 \sum_{j,k\in\Zd} y_j y_k \phi(j-k)
   =  \intrd \trig \phihat(\xi)\,d\xi,
\eqno{(\chapno4.1)}$$
where $\phihat(\xi) = -H(\xi)$ for almost all $\xi \in \Rd$ and $H$ is given
by (\chapno2.27). Moreover, (\chapno2.6) is valid, that is
$$ \sum_{j,k\in\Zd} y_j y_k \phi(j-k) = 
     \inttd \trig \sigma(\xi)\,d\xi, \eqno{(\chapno4.2)}
$$
where $\sigma$ is given by (\chapno2.7). Applying (\chapno2.30), we have
$$\eqalign{
   |\sigma(\xi)| 
  &= \sum_{k\in\Zd} \Bigl|\phihat(\xi+2\pi k)\Bigr| \cr
  &= \int_0^\infty 
        \sum_{k\in\Zd} \Ghat(t^{-1/2}(\xi+2\pi k)) 
            \ t^{-d/2-1}\,d\alpha(t). 
}\eqno{(\chapno4.3)}$$
As in Section 2, we consider essential upper and lower bounds on
$\sigma$. Let us begin this study by fixing $t > 0$ and analysing 
the function
$$ \tau (\xi) = \sum_{k\in\Zd} \Ghat(t^{-1/2}(\xi+2\pi k)), \qquad \xi
\in \Rd. \eqno{(\chapno4.4)}$$
By (\chapno3.14), we have
$$ \tau(\xi) = \prod_{j=1}^d {E_j(\xi_j) \over \Lambda_j(0)}, \qquad \xi \in \Rd,
\eqno{(\chapno4.5)}$$
where
$$ E_j(x) = \sum_{k\in\ZZ} \Lambdahat_j((x+2\pi k)t^{-1/2}), \qquad x \in
\R1, \quad j = 1, \ldots, d. \eqno{(\chapno4.6)}$$
We now employ the following key lemma.

\proclaim {Lemma \chapno4.1}. Let $\Lambda \in \PP$ and let 
$$E(x) = \sum_{k\in\ZZ} \Lambdahat((x+2\pi k)t^{-1/2}), \qquad x \in
\R1.$$ 
Then $E$ is an even function and $E(0) \ge E(x) \ge E(y) \ge E(\pi)$ 
for every $x$ and $y$ in $\R1$ with $0 \le x \le y \le \pi$.

\pf The exponential decay of $\Lambda$ and the absolute integrability
of $\Lambdahat$ imply that the Poisson summation formula is valid,
which gives the relation
$$ E(x) = t^{1/2} \sum_{k\in\ZZ} \Lambda(kt^{1/2}) e^{ikx}, \qquad x
\in \R1. \eqno{(\chapno4.7)}$$
Now the sequence $a_k := \Lambda(kt^{1/2})$, $k \in \ZZ$, is an even,
exponentially decaying P\'olya frequency sequence, that is every minor
of the Toeplitz matrix $(a_{j-k})_{j,k\in\ZZ}$ is non-negative
definite (and we see that this is a consequence of (\chapno3.7)). By a 
result of Edrei (1953), $\sum_{k\in\ZZ} a_k z^k$ is a meromorphic
function on an annulus $\{z\in\CC : 1/R \le |z| \le R \}$, for some $R
> 1$, and enjoys an infinite product expansion of
the form
$$ \sum_{k\in\ZZ} a_k z^k 
  = C e^{\lambda(z + z^{-1})} 
   \prod_{j=1}^\infty { (1+\alpha_j z)(1+\alpha_j z^{-1}) \over
(1-\beta_j z)(1-\beta_j z^{-1})}, \qquad z \ne 0, \eqno{(\chapno4.8)}$$ 
where $C \ge 0$, $\lambda \ge 0$, $0 < \alpha_j, \beta_j < 1$ and
$\sum_{j=1}^\infty \alpha_j + \beta_j < \infty$. Hence
$$ E(x) = C t^{1/2} e^{2\lambda\cos x} \prod_{j=1}^\infty {1 +
2\alpha_j\cos x + \alpha_j^2 \over 1 - 2\beta_j \cos x + \beta_j^2},
\qquad x\in\R1. \eqno{(\chapno4.9)}$$
Observe that each term in the product is an even function which is
decreasing on $[0,2\pi]$, which provides the required inequality.~\qed

In particular, $E_j(x) \ge E_j(\pi)$ for $j = 1, \ldots, d$, where
$E_j$ is given by (\chapno4.6). Hence 
$$ \tau(\xi) \ge \tau(\pi e), \qquad \xi \in \Rd,
\eqno{(\chapno4.10)}$$
and applying (\chapno4.3) we get
$$ |\sigma(\xi)| \ge |\sigma(\pi e)|, \qquad \xi \in \Rd.
\eqno{(\chapno4.11)}$$
We now come to our principal result.

\proclaim {Theorem \chapno4.2}. Let $(y_j)_{j\in\Zd}$ be a zero-summing
sequence and let $\phi \in \CCC$. Then we have the inequality
$$
  \Bigl| \sum_{j,k\in\Zd} y_ky_k \phi(j-k) \Bigr|
     \ge |\sigma(\pi e)| \sum_{j\in\Zd} y_j^2.
\eqno{(\chapno4.12)}$$

\pf Equation (\chapno4.2) and the Parseval relation provide the inequality
$$
   \Bigl| \sum_{j,k\in\Zd} y_ky_k \phi(j-k) \Bigr|
     \ge
   |\sigma(\pi e)| \inttd \trig\,d\xi = |\sigma(\pi e)| \sum_{j\in\Zd} y_j^2,
$$
as in inequality (\chapno2.10).~\qed

Of course, we are interested in showing that (\chapno4.12) cannot be
improved, that is $|\sigma(\pi e)|$ cannot be replaced by a larger
number independent of $(y_j)_{j\in\Zd}$. Recalling Proposition \chapno2.2,
this is true if $\sigma$ is continuous at $\pi e$. In fact, we can use
Lemma \chapno4.1 to prove that $\sigma$ is continuous everywhere in the set
$(0,2\pi)^d$. We first collect some necessary preliminary results.

\proclaim {Lemma \chapno4.3}. The function $\tau$ given by (\chapno4.4) is continous
for every $t > 0$ and satisfies the inequality
$$ \tau(\xi) \le \tau(\eta) 
\hbox{ for } 0 \le \eta \le \xi \le \pi e.
\eqno{(\chapno4.13)}$$
Furthermore,
$$
 \tau(\pi e + \xi) = \tau(\pi e - \xi) \hbox{ for all } \xi \in (-\pi,\pi)^d.
\eqno{(\chapno4.14)}$$   

\pf The definition of $G$, (\chapno4.5) and (\chapno4.7) provide the Fourier series
$$ \tau(\xi) = t^{d/2} \sum_{k\in\Zd} G(kt^{1/2}) e^{ik\xi}, \qquad
\xi \in \Rd,
\eqno{(\chapno4.15)}$$
and the exponential decay of $G$ implies the uniform convergence of
this series. Hence $\tau$ is continuous, being the uniform limit of the
finite sections of (\chapno4.15).

Applying the product formula (\chapno4.5) and Lemma \chapno4.1, we obtain
(\chapno4.13) and (\chapno4.14).~\qed

\proclaim {Proposition \chapno4.4}.  $\sigma$ is continuous on $(0,2\pi)^d$.

\pf Equation (\chapno4.2) implies that $\trig |\sigma(\xi)| < \infty$ for
almost every $\xi \in \Td$. Consequently, $\sigma$ is finite almost
everywhere, by Lemma \chapno2.1. Thus every non-empty open subset of $\Td$
contains a point at which $\sigma$ is finite. Specifically, let
$\delta \in (0,\pi)$ and define the closed set $K_\delta :=
[\delta,2\pi-\delta]^d$. Thus the open set $\Td \setminus K_\delta$
contains a point, $\eta$ say, for which
$$
 \infty > |\sigma(\eta)| 
   = \int_0^\infty \sum_{k\in\Zd} \Ghat((\eta+2\pi k)t^{-1/2})\,
t^{-d/2 - 1} \, d\alpha(t).
\eqno{(\chapno4.16)}$$
Let us show that $\sigma$ is continuous in $K_\delta$. To this end,
choose any convergent sequence $(\xi_n)_{n=1}^\infty$ in $K_\delta$ and let
$\xi_\infty$ denote its limit. We must prove that 
$$\lim_{n\to\infty}
\sigma(\xi_n) = \sigma(\xi_\infty).$$ 
Now Lemma \chapno4.3 and (\chapno4.16) supply
the bound
$$ |\sigma(\xi_n)| \le |\sigma(\eta)| < \infty, \qquad n = 1, 2, \ldots, $$
that is the functions 
$$\{ t \mapsto \sum_{k\in\Zd} \Ghat((\xi_n+2\pi k)t^{-1/2})
t^{-d/2-1}\,d\alpha(t) : n = 1, 2, \ldots \,\}$$
are absolutely integrable on $[0,\infty)$ with respect to
the measure $d\alpha$. Moreover, they are dominated by the absolutely
integrable function $t \mapsto \sum_{k\in\Zd} \Ghat((\eta+2\pi k)t^{-1/2})
t^{-d/2-1}$. However, the continuity of $\tau$ proved in Lemma \chapno4.3
allows to deduce that
$$ \lim_{n\to\infty} \sum_{k\in\Zd} \Ghat((\xi_n+2\pi k)t^{-1/2}) t^{-d/2-1}
     = \sum_{k\in\Zd} \Ghat((\xi_\infty + 2\pi k)t^{-1/2}) t^{-d/2-1},
$$
for all positive $t$.
Thus the dominated convergence theorem implies that
$\sigma(\xi_n) \to \sigma(\xi_\infty)$ as $n$ tends to infinity. Since $\delta \in
(0,\pi)$ was arbitrary, we conclude that $\sigma$ is continuous in all
of $(0,2\pi)^d$.~\qed

\proclaim {Corollary \chapno4.5}. Inequality (\chapno4.12) cannot be improved for
$\phi \in \CCC$ if we can find a trigonometric polynomial $P$
satisfying the conditions of Proposition 5.2.2 at the point $\pi e$.

\pf We simply apply Proposition 5.5.2.~\qed

\sect {\chapno5. Total positivity and the Gaussian cardinal function}
This material is not directly related to the earlier sections of this
chapter, but it does use a total positivity property to deduce an
interesting fact concerning infinity norms of Gaussian distance
matrices generated by infinite regular grids.

Let $\lambda$ be a positive constant and let $\phi \colon \R1
\to \R1$ be the Gaussian 
$$
 \phi(x) = \exp(-\lambda x^2), \qquad x \in \R1. \eqno{(\chapno5.1)}
$$
It is known (see Buhmann (1990)) that there exists a real sequence $(c_k)_{k\in\ZZ}$ such
that $\sum_{k\in\ZZ} c_k^2 < \infty$ and the function $\chi \colon \R1
\to \R1$ given by
$$ 
 \chi(x) = \sum_{k\in\ZZ} c_k \phi(x-k), \qquad x \in \R1, 
\eqno{(\chapno5.2)}$$
satisfies the equation
$$
 \chi(j) = \delta_{0j}, \qquad j \in \ZZ .
$$
Thus $\chi$ is the cardinal function of interpolation for the Gaussian
radial basis function.

\proclaim {Proposition \chapno5.1}. The coefficients $(c_k)_{k\in\ZZ}$ of the
cardinal function $\chi$ alternate in sign, that is $(-1)^k c_k \ge 0$
for every integer $k$.

\pf For each non-negative integer $n$, we let
$$ 
 A_n = (\phi(j-k))_{j,k=-n}^n . \eqno{(\chapno5.3)}
$$
Now $A_n$ is an invertible totally positive matrix, which implies that
$A_n^{-1}$ enjoys the ``chequerboard'' property, that is the elements
of the inverse matrix satisfy $(-1)^{j+k} (A_n^{-1})_{jk} \ge 0$, for
$j,k = -n, \ldots, n$. In particular, if we let
$$
  c_k^{(n)} = (A_n^{-1})_{0k}, \qquad k = -n, \ldots, n, 
\eqno{(\chapno5.4)}$$
then $(-1)^k c_k^{(n)} \ge 0$ and the definition of $A_n^{-1}$
provides the equations
$$
 \sum_{k=-n}^n c_k^{(n)} \phi(j-k) = \delta_{0j}, \qquad j = -n,
\ldots, n. 
\eqno{(\chapno5.5)}$$
In other words, the function $\chi_n \colon \R1 \to \R1$ defined by
$$
 \chi_n(x) = \sum_{k=-n}^n c_k^{(n)} \phi(x-k), \qquad x \in \R1,
\eqno{(\chapno5.6)}$$
provides the cardinal function of interpolation for the finite set
$\{-n, \ldots, n\}$.

Now Theorem 9 of Buhmann and Micchelli (1991) provides the following
useful fact relating the coefficients of $\chi_n$ and $\chi$:
$$
 \lim_{n\to\infty} c_k^{(n)} = c_k, \qquad k \in \ZZ.
$$
Thus the property $(-1)^k c_k^{(n)} \ge 0$ implies the required
condition
$(-1)^k c_k \ge 0$.~\qed

We now consider the bi-infinite symmetric Toeplitz matrix $A =
(\phi(j-k))_{j,k\in\ZZ}$ as a bounded linear operator $A \colon
\ell^p(\ZZ) \to \ell^p(\ZZ)$ when $p \ge 1$. Thus 
$A^{-1} = (c_{j-k})_{j,k\in\ZZ}$, where the $(c_j)_{j\in\ZZ}$ are
given by (\chapno5.2), and a theorem of Buhmann (1990) provides the equation
$$
 c_k = (2\pi)^{-1} \int_0^{2\pi} {1 \over \sigma(\xi)}\,
e^{-ik\xi}\,d\xi, \qquad k \in \ZZ, 
\eqno{(\chapno5.7)}$$
where 
$$
 \sigma(\xi) = \sum_{k\in\ZZ} \phihat(\xi+2\pi k), \qquad \xi \in \R1.
\eqno{(\chapno5.8)}$$
Therefore, using standard results of Toeplitz operator theory
(Grenander and Szeg\H{o} (1984)), 
we obtain the expression
$$ \| A^{-1} \|_2 
   = \max\{ {1\over \sigma(\xi)} : \xi \in [0,2\pi]\}. 
$$
Applying Lemma 4.2.7, we get
$$
 \| A^{-1} \|_2 = {1 \over \sigma(\pi)} 
    = \sum_{k\in\ZZ} (-1)^k c_k. \eqno{(\chapno5.9)}
$$
But Proposition \chapno5.1 and the symmetry of $A$ provide the relations
$$ \|A^{-1}\|_1 = \|A^{-1}\|_\infty = \sum_{k\in\ZZ} |c_k| =
\sum_{k\in\ZZ} (-1)^k c_k, \eqno{(\chapno5.10)},$$
so that $A^{-1}$ provides a nontrivial linear operator on
$\ell^p(\ZZ)$, for $p = 1, 2$, and $\infty$, whose norms agree on each
of these sequence spaces. Further, we recall that $\log \|A^{-1}\|_p$
is a convex function of $1/p$ for $p \ge 1$, which is a consequence of
the Riesz-Thorin theorem (Hardy {\it et al} (1952), pp. 214, 219). Hence we have proved the
interesting fact that $\|A^{-1}\|_p = \|A^{-1}\|_1$ for all $p \ge 1$.

In the multivariate case, the cardinal function is given by
expressions analogous to (\chapno5.2) and (\chapno5.7). Specifically, we let $\phi(x) =
\exp(-\lambda\|x\|^2)$, $x \in \Rd$, and then $\chi \colon \Rd \to
\Rd$ is defined by
$$
  \chi(x) = \sum_{k\in\Zd} c_k^{(d)} \phi(x-k), \qquad x \in \Rd,
\eqno{(\chapno5.11)}$$
where
$$ c_k^{(d)} = (2\pi)^{-d} \int_\Td {1 \over \sigma^{(d)}(\xi)}\,
e^{-ik\xi}\,d\xi, \qquad k = (k_1, \ldots, k_d) \in \Zd,
\eqno{(\chapno5.12)}$$
and 
$$
 \sigma^{(d)}(\xi) = \sum_{k\in\Zd} \phihat(\xi+2\pi k).
\eqno{(\chapno5.13)}$$
The key point is that $\phi$ is a tensor product of univariate functions, which
implies the relation
$$
 \sigma^{(d)}(\xi) = \prod_{j=1}^d \sigma(\xi_j), \qquad \xi = (\xi_1,
\ldots, \xi_d) \in \Rd,
\eqno{(\chapno5.14)}$$
where $\sigma$ is given by (\chapno5.8). Consequently the coefficients of the
multivariate cardinal function are related to those of the univariate
cardinal function by the formula
$$
 c_k^{(d)} = \prod_{j=1}^d c_{k_j}, \qquad k = (k_1, \ldots, k_d) \in
\Zd.
\eqno{(\chapno5.15)}$$
In particular, the following corollary is an immediate consequence of
Proposition \chapno5.1.

\proclaim {Corollary \chapno5.2}. $(-1)^{k_1+\cdots+k_d} c_k^{(d)} \ge 0$ for
every integer $k = (k_1, \ldots, k_d) \in \Zd$.

\vfill \eject
\def\chapno{6.}
\headline{\ifnum\pageno=95\hfil\else{\headfont\hfil\chptitle\hfil}\fi}
\def\chptitle{Preconditioned conjugate gradients}
\bigskip
\centerline{\bf 6 : Norm Estimates and Preconditioned Conjugate Gradients}


\def\pts{{(x_j)_{j=1}^n}}

\def\sjk{{\sum_{j,k=1}^n y_j y_k\,}}
\def\djk{{\| x_j - x_k \|}}
\def\xjk{{x_j-x_k}}
\def\yjk{{|\sum_{j=1}^n y_j e^{ix_j t}|^2}}
\def\trig{{|\sum_{j=1}^n y_j e^{ijt}|^2}}

\def\skzd{{\sum_{k\in\Zd}}}

\def\rint{{\int_{-\infty}^{\infty}}}
\def\hint{{\int_0^\infty}}

\def\Sdm1{{S^{d-1}}}

\def\ysqr{{\Vert y \Vert^2}}

\def\ainv{{A_n^{-1}}}
\def\aninv{{\Vert A_n^{-1} \Vert_2}}

\def\zsum{{\sum_{k=-\infty}^\infty}}

\def\hi{{(2\pi)^{-1}}}
\def\tint{{\hi\int_0^{2\pi}}}
\def\tintd{{(2\pi)^{-d}\int_{\Td}}}
\def\fhat{{\hat f}}
\def\ajxi{{\Bigl|\sum_{j\in\Zd} a_j \exp(ij\xi)\Bigl|^2}}
\def\Rdm0{{\Rd\setminus\{0\}}}
\def\ghat{{\hat g}}
\def\ytrigo{{\Bigl|\sum_{j\in\Zd} y_j \exp(ix_j \xi)\Bigl|^2}}
\def\ytrig2{{\Bigl|\sum_{j\in\Zd} y_j \exp(ij \xi)\Bigl|^2}}
\def\ynj{{y_j^{(n)}}}

\def\phichat{{\phihat_c}}
\def\chichat{{{\hat \chi}_c}}

\def\min{{\rm min}}

\font\sc=cmcsc10

\def\pf{{\noindent{\it Proof.\ }}}
\def\qed{{\hfill$\square$\vskip 10pt}}

\def\R1{{\cal R}}
\def\Rd{{\R1^d}}
\def\Rn{{\R1^n}}
\def\ZZ{{{\cal Z} }}
\def\Z0{{\sum_{j=1}^n y_j = 0}}
\def\Zd{{\ZZ^d}}
\def\T1{{[0,2\pi]}}
\def\Td{{[0,2\pi]^d}}
\def\CC{{\cal C}}

\def\half{{1\over2}}
\def\OO{{\cal O}}
\def\phi{{\varphi}}
\def\phihat{{\hat\varphi}}

\outer\def\beginsection#1\par{\vskip0pt plus.3\vsize\penalty-150
 \vskip0pt plus-.3\vsize\bigskip\vskip\parskip
 \message{#1}\centerline{\sc#1}\nobreak\smallskip\noindent}






\def\nullx{\hfill}

\def\rightheadline{\ifnum\pageno=1 \nullx%
  \else\sc\hfil\chptitle\hfil\fi}
  \def\leftheadline{\ifnum\pageno=1 \nullx%
  \else\hfil\sc\authors\hfil\fi}  

\def\theta{{\vartheta}}
\def\psihat{\hat\psi}
\def\Rnn{\R1^{n \times n}}
\def\Image{\hbox{Im\ }}
\def\Sp{\hbox{Sp\ }}
\def\vol{\hbox{vol}}
\def\Cd{\CC^d}


\sect {\chapno1. Introduction}
Let $n$ be a positive integer and let $A_n$ be the symmetric Toeplitz
matrix given by
$$ A_n = \left(\phi(j-k)\right)_{j,k=-n}^n, \eqno{(\chapno1.1)}$$
where $\phi \colon \R1 \to \R1$ is either a Gaussian ($\phi(x) =
\exp(-\lambda x^2)$ for some positive constant $\lambda$) or a
multiquadric ($\phi(x) = (x^2 + c^2)^{1/2}$ for some real constant
$c$). In this section we construct efficient preconditioners for the
conjugate gradient solution of the linear system
$$A_n x = f, \qquad f \in \R1^{2n+1}, \eqno{(\chapno1.2)}$$
when $\phi$ is a Gaussian, or the augmented linear system
$$\eqalign{ A_n x + e y &= f, \cr
            e^T x &= 0,\cr} \eqno{(\chapno1.3)}$$
when $\phi$ is a multiquadric. Here $e = [1, 1, \ldots, 1]^T \in \R1^{2n+1}$ and $y \in \R1$.
Section \chapno2 describes the construction for the Gaussian and Section \chapno3
deals with the multiquadric. Of course, we exploit the Toeplitz
structure of $A_n$ to perform a matrix-vector multiplication in $\OO(n
\log n)$ operations whilst storing $\OO(n)$ real numbers. Further, we
shall see numerically that the number of iterations required to
achieve a solution of (\chapno1.2) or (\chapno1.3) to within a given tolerance is
independent of $n$.

Our method applies to many other radial basis functions, such as the
inverse multiquadric ($\phi(x) = (x^2 + c^2)^{-1/2}$) and the thin
plate spline ($\phi(x) = x^2 \log |x|$). However, we concentrate on
the Gaussian and the multiquadric because they exhibit most of the
important features of our approach in a concrete setting. Similarly we
only touch briefly on the $d$-dimensional analogue of (\chapno1.1), that is 
$$ A_n^{(d)} = \left( \phi(j-k) \right)_{j,k\in [-n,n]^d}. \eqno{(\chapno1.4)}$$
We shall still call $A_n^{(d)}$ a Toeplitz matrix. Moreover
the matrix-vector multiplication
$$ A_n^{(d)} x = \left( \sum_{k \in [-n,n]^d} \phi(\|j-k\|) x_k
\right)_{j\in [-n,n]^d}, \eqno{(\chapno1.5)}$$
where $\|\cdot\|$ is the Euclidean norm and $x = (x_j)_{j\in [-n,n]^d}$, can
still be calculated in $\OO(N \log N)$ operations, where $N = (2n+1)^d$, whilst
requiring $\OO(N)$ real numbers to be stored. This trick is a simple extension
of the Toeplitz matrix-vector multiplication method when $d = 1$, but seems
to be  less familiar for $d$ greater than one. This will be dealt with
in detail in Baxter (1992c).

\sect {\chapno2. The Gaussian}
Our treatment of the preconditioned conjugate gradient (PCG) method follows
Section 10.3 of Golub and Van Loan (1989), and we begin with a general
description. We let $n$ be a positive integer and $A \in \Rnn$ be a symmetric
positive definite matrix. For any nonsingular symmetric matrix $P \in \Rnn$
and $b \in \R1^n$ we can use the following iteration to solve the linear system
$PAPx =Pb$.

\proclaim {Algorithm \chapno2.1}. Choose any $x_0$ in $\Rn$. Set $r_0 = Pb - PAP x_0$
and $d_0 = r_0$.

{\obeylines 
For $k=0, 1, 2, \ldots$ do begin
\quad $a_k = r_k^Tr_k / d_k^T PAP d_k$
\qquad $x_{k+1} = x_k + a_k d_k$
\qquad $r_{k+1} = r_k - a_k PAP d_k$
\quad $b_k = r_{k+1}^T r_{k+1} / r_k^T r_k$
\qquad $d_{k+1} = r_{k+1} + b_k d_k$
\quad Stop if $\|r_{k+1}\|$ or $\|d_{k+1}\|$ is sufficiently small.
end.\par}
\medskip

\noindent In order to simplify Algorithm \chapno2.1 define
$$C = P^2, \qquad
\xi_k = P x_k, \qquad
r_k = P \rho_k \qquad
\hbox{ and } \qquad \delta_k = P d_k.
\eqno{(\chapno2.1)}$$
Substituting in Algorithm \chapno2.1 we obtain the following method.
\vfill\eject

\proclaim {Algorithm \chapno2.2}. Choose any $\xi_0$ in $\Rn$. Set
$\rho_0 = b - A \xi_0$, $\delta_0 = C \rho_0$.

{\obeylines 
For $k=0, 1, 2, \ldots$ do begin
\quad $a_k = \rho_k^T C \rho_k / \delta_k^T A \delta_k$
\qquad $\xi_{k+1} = \xi_k + a_k \delta_k$
\qquad $\rho_{k+1} = \rho_k - a_k A \delta_k$
\quad $b_k = \rho_{k+1}^T C \rho_{k+1} / \rho_k^T C \rho_k$
\qquad $\delta_{k+1} = C \rho_{k+1} + b_k \delta_k$
\quad Stop if $\|\rho_{k+1}\|$ or $\|\delta_{k+1}\|$ is sufficiently small.
end.\par}
\medskip
It is Algorithm \chapno2.2 that we shall consider as our PCG method in this
section, and we shall call $C$ the preconditioner. We see that the
only restriction on $C$ is that it must be a symmetric positive
definite matrix, but we observe that the spectrum of $CA$ should
consist of a small number of clusters, preferably one cluster
concentrated at one. At this point, we also mention that the condition
number of $CA$ is not a reliable guide to the efficacy of our
preconditioner. For example, consider the two cases when (i) $CA$ has
only two different eigenvalues, say $1$ and $100,000$, and (ii) when
$CA$ has eigenvalues uniformly distributed in the interval $[1, 100]$.
The former has the larger condition number but, in exact
arithmetic, the answer will be achieved in two steps, whereas the
number of steps can be as high as $n$ in the latter case. Thus the
term ``preconditioner'' is sometimes inappropriate, although its usage
has become standard.

We can shed no light on the problem of constructing
preconditioners for the general case.Accordingly, we let $A$ be the matrix $A_n$ of (\chapno1.1) and let $\phi(x) = \exp(-
x^2)$. Thus $A_n$ is positive definite and can be embedded in the bi-infinite symmetric
Toeplitz matrix
$$ A_\infty = \left(\phi(j-k) \right)_{j,k\in\ZZ}. \eqno{(\chapno2.2)}$$
The classical theory of Toeplitz operators (see, for instance,
Grenander and Szeg\H{o} (1984)) and the work of Section 4 provide the
relations   
$$ \Sp A_n \subset \Sp A_\infty = [\sigma(\pi), \sigma(0)] \subset (0,\infty),
\eqno{(\chapno2.3)}$$
where $\sigma$ is the symbol function
$$ \sigma(\xi) = \sum_{k\in\ZZ} \phihat(\xi+2\pi k), \qquad \xi \in
\R1. \eqno{(\chapno2.4)}$$
Further, Theorem 9 of Buhmann and Micchelli (1991) allows us to
conclude that, for any fixed integers $j$ and $k$, we have
$$ \lim_{n \to \infty} (A_n^{-1})_{j,k} = (A_\infty^{-1})_{j,k}.
\eqno{(\chapno2.5)}$$
It was equations (\chapno2.3) and (\chapno2.5) which led us to investigate the possibility of
using some of the elements of $A_n^{-1}$ for a relatively small value of $n$
to construct preconditioners for $A_N$, where $N$ may be much larger
than $n$. Specifically, let us choose integers $0 < m \le n$ and define the
sequence $$c_j = (A_n^{-1})_{j0}, \qquad j = -m, \ldots, m. \eqno{(\chapno2.6)}$$
We now let $C_N$ be the $(2N+1) \times (2N+1)$ banded symmetric Toeplitz
matrix
$$ C_N = \pmatrix{ c_0 & \ldots & c_m & & & \cr
                   \vdots & \ddots & & \ddots & & \cr
                   c_m & & & & & \cr
                   & \ddots & & & & c_m \cr
                   & & & & & \vdots \cr
                   & & & c_m & \ldots & c_0 \cr}. \eqno{(\chapno2.7)}
$$
We claim that, for sufficiently large $m$ and $n$, $C_N$ provides an
excellent preconditioner when $A = A_N$ in Algorithm \chapno2.2. Before
discussing any theoretical motivation for this choice of
preconditioner, we present an example. We let $n=64$, $m=9$ and $N =
32,768$. Constructing $A_n$ and calculating the elements $\{
(A_n^{-1})_{j0}: j=0, 1, \ldots, m \}$ we find that
$$ \pmatrix{ c_0 \cr c_1 \cr \vdots \cr c_9 } 
  = \pmatrix{    
  \ \,\,1.4301 \times 10^0 \cr
  -5.9563 \times 10^{-1}\cr
  \ \,\, 2.2265 \times 10^{-1}\cr
  -8.2083 \times 10^{-2}\cr
  \ \,\, 3.0205 \times 10^{-2}\cr
  -1.1112 \times 10^{-2}\cr
  \ \,\, 4.0880 \times 10^{-3}\cr
  -1.5039 \times 10^{-3}\cr
  \ \,\, 5.5325 \times 10^{-4}\cr
  -2.0353 \times 10^{-4}\cr}.
\eqno{(\chapno2.8)}$$

\vfill
\eject
\centerline{\psfig{file=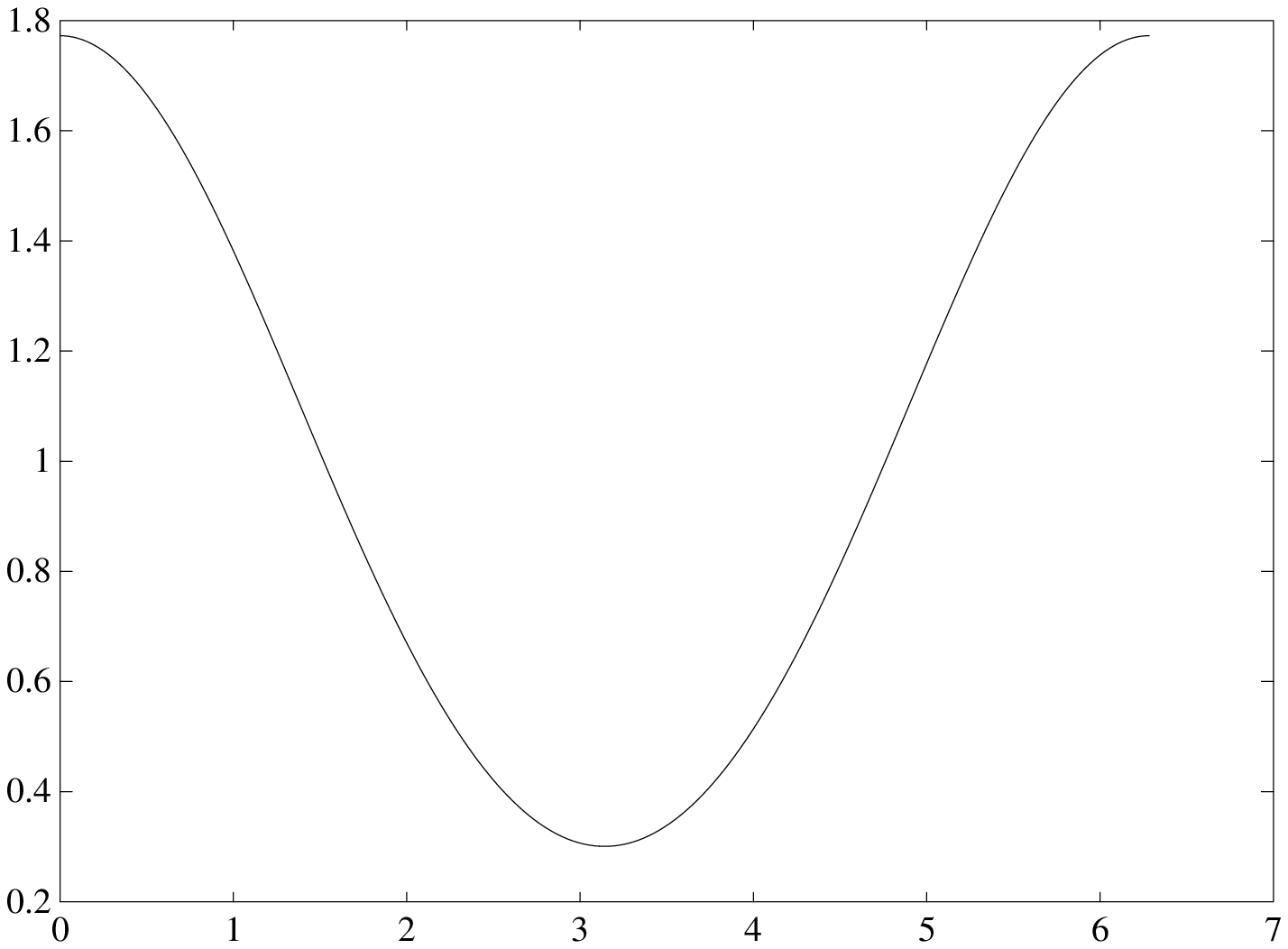}}
\centerline{{\bf FIGURE \chapno1:} The symbol function for $C_\infty$.}

Now $C_N$ can be embedded in the bi-infinite Toeplitz matrix
$C_\infty$ defined by
$$ (C_\infty)_{jk} = \cases{ c_{j-k}, \qquad |j-k| \le m, \cr 0,
\qquad\qquad |j-k| > m,} \eqno{(\chapno2.9)}$$
and the symbol for this operator is the trigonometric polynomial
$$ \sigma_{C_\infty}(\xi) = \sum_{j=-m}^m c_j e^{ij\xi}, \qquad \xi
\in \R1. \eqno{(\chapno2.10)}$$

In Figure 6.1 we display a graph of $\sigma_{C_\infty}$ for $0 \le \xi
\le 2\pi$, and it is clearly a positive function. Thus the relations
$$ \Sp C_N \subset \Sp C_\infty = \{\sigma_{C_\infty}(\xi) : \xi \in
[0,2\pi] \} \subset (0, \infty) \eqno{(\chapno2.11)}$$
imply that $C_N$ is positive definite. Hence it is suitable to use
$C_N$ as the preconditioner in Algorithm \chapno2.2. Our aim in this example
is to compare this choice of preconditioner with the use of the
identity matrix as the preconditioner. To this end, we let the
elements of the righthandside vector $b$ of Algorithm \chapno2.2 be random
real numbers uniformly distributed in the interval $[-1,1]$. Applying
Algorithm \chapno2.2 using the identity matrix as the preconditioner provides
the results of Table \chapno1. Table \chapno2 contains the analogous results
using (\chapno2.7) and (\chapno2.8). In both cases the iterations were stopped when
the residual vector satisfied the bound $\|r_{k+1}\| /\|b\| < 10^{-13}$.
The behaviour shown in the tables is typical; we find that the number of steps
required is independent of $N$ and $b$.

$$\vbox{
\settabs\+\indent&c\qquad\hskip 24pt&Optimal bound\quad&\cr
\+&{\rm Iteration}&{\rm Error}\cr
\+& $1$& $2.797904 \times 10^{1}$\cr
\+& $10$& $1.214777 \times 10^{-2}$\cr
\+& $20$& $1.886333 \times 10^{-6}$\cr
\+& $30$& $2.945903 \times 10^{-10}$\cr
\+& $33$& $2.144110 \times 10^{-11}$\cr
\+& $34$& $8.935534 \times 10^{-12}$\cr

}$$
\centerline{{\bf Table \chapno1:} No preconditioning}
\smallskip

$$\vbox{
\settabs\+\indent&c\qquad\hskip 24pt&Optimal bound\quad&\cr
\+&{\rm Iteration}&{\rm Error}\cr
\+& $1$& $2.315776 \times 10^{-1}$\cr
\+& $2$& $1.915017 \times 10^{-3}$\cr
\+& $3$& $1.514617 \times 10^{-7}$\cr
\+& $4$& $1.365228 \times 10^{-11}$\cr
\+& $5$& $1.716123 \times 10^{-15}$\cr

}$$
\centerline{{\bf Table \chapno2:} Using (\chapno2.7) and
(\chapno2.8) as the preconditioner}
\smallskip

Why should (\chapno2.7) and (\chapno2.8) provide a good preconditioner? Let us
consider the bi-infinite Toeplitz matrix $C_\infty A_\infty$. The
spectrum of this operator is given by
$$ \Sp C_\infty A_\infty = \{ \sigma_{C_\infty}(\xi) \sigma(\xi) : \xi
\in [0,2\pi] \}, \eqno{(\chapno2.12)}$$
where $\sigma$ is given by (\chapno2.4) and $\sigma_{C_\infty}$ by (\chapno2.10).
Therefore in order to concentrate $\Sp C_\infty A_\infty$ at unity we
must have
$$ \sigma_{C_\infty}(\xi) \sigma(\xi) \approx 1, \qquad \xi \in
[0,2\pi]. \eqno{(\chapno2.13)}$$
In other words, we want $\sigma_{C_\infty}$ to be a trigonometric
polynomial approximating the continuous function $1/\sigma$. Now if
the Fourier series of $1/\sigma$ is given by
$$ \sigma^{-1}(\xi) = \sum_{j\in\ZZ} \gamma_j e^{ij\xi}, \qquad \xi
\in \R1, \eqno{(\chapno2.14)}$$
then its Fourier coefficients $(\gamma_j)_{j\in\ZZ}$ are the
coefficients of the cardinal function $\chi$ for the integer grid,
that is
$$ \chi(x) = \sum_{j\in\ZZ} \gamma_j \phi(x-j), \qquad x \in \R1,
\eqno{(\chapno2.15)}$$
and 
$$   \chi(k) = \delta_{0k}, \qquad k \in \ZZ. \eqno{(\chapno2.16)}$$
(See, for instance, Buhmann (1990).) Recalling (\chapno2.5), we deduce that
one way to calculate approximate values of the coefficients
$(\gamma_j)_{j\in\ZZ}$ is to solve the linear system
$$ A_n c^{(n)} = e^0, \eqno{(\chapno2.17)}$$
where $e^0 = (\delta_{j0})_{j=-n}^n \in \R1^{2n+1}$. This observation
is not new; indeed Buhmann and Powell (1990) used precisely this idea
to calculate approximate values of the cardinal function $\chi$. We
now set
$$ c_j = c_j^{(n)}, \qquad 0 \le j \le m, \eqno{(\chapno2.18)}$$
and we observe that the symbol function $\sigma$ for the Gaussian is a
theta function (see Section 4.2). Thus $\sigma$ is a positive
continuous function whose Fourier series is absolutely convergent.
Hence $1/\sigma$ is a positive continuous function and Wiener's lemma
(Rudin (1973)) implies the absolute convergence, and therefore the uniform
convergence, of its Fourier series. We deduce that the symbol function
$\sigma_{C_\infty}$ can be chosen to approximate $1/\sigma$ to within
any required accuracy. More formally we have the 

\proclaim {Lemma \chapno2.3}. Given any $\epsilon > 0$, there are positive
integers $m$ and $n_0$ such that
$$ \Bigl| \sigma(\xi)\sum_{j=-m}^m c_j^{(n)} e^{ij\xi} - 1 \Bigr| \le
\epsilon, \qquad \xi \in [0,2\pi],$$
for every $n \ge n_0$, where $c^{(n)} = (c_j^{(n)})_{j=-n}^n$ is given
by (\chapno2.17).

\pf The uniform convergence of the Fourier series for $\sigma$ implies
that we can choose $m$ such that
$$ \Bigl|\sigma(\xi)\sum_{j=-m}^m  \gamma_j e^{ij\xi} -1 \Bigr| \le
\epsilon, \qquad \xi \in [0,2\pi].$$
By (\chapno2.5), we can also choose $n_0$ such that $|\gamma_j - c_j^{(n)}|
\le \epsilon$ for $j = -m, \ldots, m$ and $n \ge n_0$. Then we have
$$\eqalign{
 \Bigl| \sigma(\xi) \sum_{j=-m}^m c_j^{(n)} e^{ij\xi} - 1 \Bigr| 
   &\le \Bigl| \sigma(\xi)\sum_{j=-m}^m \gamma_j e^{ij\xi} -1 \Bigr| +
\sigma(\xi) \Bigl| \sum_{j=-m}^m (\gamma_j - c_j^{(n)}) e^{ij\xi}| \cr
   &\le \epsilon [ 1 + (2m+1)\sigma(0)], \cr}$$
remembering from Chapter 4 that $0 < \sigma(\pi) \le \sigma(\xi) \le
\sigma(0)$. Since $\epsilon$ is arbitrary the lemma is true.~\qed

\sect{\chapno3. The Multiquadric}
The multiquadric interpolation matrix
$$ A = \Bigl( \phi(\|x_j - x_k\|) \Bigr)_{j,k=1}^n, $$
where $\phi(r) = (r^2 + c^2)^{1/2}$ and $(x_j)_{j=1}^n$ are points in
$\Rd$, is not positive definite. We recall from Chapter 2 that it is
{\it almost negative definite}, that is for any real numbers
$(y_j)_{j=1}^n$ satisfying $\sum y_j = 0$ we have
$$ \sum_{j,k=1}^n y_j y_k \phi(\|x_j - x_k\|) \le 0. \eqno{(\chapno3.1)}$$
Furthermore, inequality (\chapno3.1) is strict whenever $n \ge 2$ and the
points $(x_j)_{j=1}^n$ are all different, and we shall assume this for the rest
of the section. In other words, $A$ is negative definite on the subspace
$<e>^\perp$, where $e = [1, 1, \ldots, 1]^T \in \R1^n$. 

Of course we cannot apply Algorithms \chapno2.1 and \chapno2.2 in this case. However
we can use the almost negative definiteness of $A$ to solve a closely
related linearly constrained quadratic programming problem:
$$\eqalign{ &\ \hbox{ minimize } \qquad \half \xi^T A \xi - b^T \xi \cr
            &\ \hbox{ subject to } \qquad e^T \xi = 0, }\eqno{(\chapno3.2)}$$
where $b$ can be any element of $\Rn$. It can be shown that 
the standard theory of Lagrange
multipliers guarantees the existence of a unique pair of vectors $\xi^* \in \R1^n$
and $\eta^* \in \R1^m$ satisfying the equations
$$\eqalign{ A \xi^* + e \eta^* &= b \cr
\hbox{ and } e^T \xi^* &= 0,}
\eqno{(\chapno3.3)}$$
where $\eta^*$ is the Lagrange multiplier vector for the constrained
optimization problem (\chapno3.2). We do not go into further detail on this point
because the nonsingularity of the matrix
$$ \pmatrix{ A & e \cr
e^T & 0} \eqno{(\chapno3.4)}$$
is well-known (see, for
instance, Powell (1990)). Instead we observe that one way to solve
(\chapno3.3) is to apply the following modification of Algorithm \chapno2.1 to (\chapno3.2).

\proclaim {Algorithm \chapno3.1}. Let $P$ be any symmetric $n \times n$
matrix such that $\ker P = \langle e \rangle$.

{\obeylines 
Set $x_0 = 0$, $r_0 = Pb - PAP x_0$, $d_0 = r_0$.
For $k=0, 1, 2, \ldots$ do begin
\quad $a_k = r_k^Tr_k / d_k^T PAP d_k$
\qquad $x_{k+1} = x_k + a_k d_k$
\qquad $r_{k+1} = r_k - a_k PAP d_k$
\quad $b_k = r_{k+1}^T r_{k+1} / r_k^T r_k$
\qquad $d_{k+1} = r_{k+1} + b_k d_k$
\quad Stop if $\|r_{k+1}\|$ or $\|d_{k+1}\|$ is sufficiently small.
\quad end.\par}
\medskip

We observe that Algorithm \chapno3.1 solves the linearly constrained
optimization problem
$$\eqalign{ &\ \hbox{ minimize } \quad \half x^T PAP x - b^T Px \cr
            &\ \hbox{ subject to } \quad e^T x = 0. }\eqno{(\chapno3.5)}$$
Moreover, the following elementary lemma implies that the solutions $\xi^*$of
(\chapno3.3) and $x^*$ of (\chapno3.5) are related by the equations $\xi^* = P x^*$.

\proclaim {Lemma \chapno3.2}. Let $S$ be any symmetric $n \times n$ matrix
and let $K = \ker S$. The $S : K^\perp \to K^\perp$ is a bijection.
In other words, given any $b \in K^\perp$ there is precisely one $a \in
K^\perp$ such that
$$ S a = b. \eqno{(\chapno3.6)}$$

\pf For any $n \times n$ matrix $M$ we have the equation
$$ \R1^n = \ker M \oplus \Image M^T.$$
Consequently the symmetric matrix $S$ satisfies
$$ \R1^n = \ker S \oplus \Image S, $$
whence $\Image S = K^\perp$. Hence for every $b \in K^\perp$ there exists
$\alpha \in \Rn$ such that $S \alpha = b$. Now we can write $\alpha =
a + \beta$, where $a \in K^\perp$ and $\beta \in K$ are uniquely
determined by $\alpha$. Thus $S a = S \alpha = b$, and (\chapno3.6) has a
solution. If $a^\prime \in K^\perp$ also satifies (\chapno3.6), then their
difference $a - a^\prime$ lies in the intersection $K \cap K^\perp =
\{0\}$, which settles the uniqeuness of $a$.~\qed

Setting $P = S$ and $K = \langle e \rangle$ in Lemma \chapno3.2 we deduce that there is
exactly one $x^* \in \langle e \rangle^\perp$ such that
$$ PAP x^* = Pb, $$
and $PAP$ is negative definite when restricted to the subspace
$\langle e \rangle^\perp$.
Follwing the development of Section \chapno2, we define
$$C = P^2, \qquad
\xi_k = Px_k, \qquad
\hbox{ and } \qquad \delta_k = P d_k,\eqno{(\chapno3.7)}$$
as in equation (\chapno2.1).
However we cannot define $\rho_k$ by (\chapno2.1) because $P$ is singular.
One solution, advocated by Dyn, Levin and Rippa (1986), is to use the
recurrence for $(\rho_k)$ embodied in Algorithm \chapno2.1 without further
ado.

\vfill\eject
\proclaim {Algorithm \chapno3.3a}. Choose any $\xi_0$ in $\langle e \rangle^\perp$. Set
$\rho_0 = b - A \xi_0$ and $\delta_0 = C \rho_0$.

{\obeylines 
For $k=0, 1, 2, \ldots$ do begin
\quad $a_k = \rho_k^T C \rho_k / \delta_k^T A \delta_k$
\qquad $\xi_{k+1} = \xi_k + a_k \delta_k$
\qquad $\rho_{k+1} = \rho_k - a_k A \delta_k$
\quad $b_k = \rho_{k+1}^T C \rho_{k+1} / \rho_k^T C \rho_k$
\qquad $\delta_{k+1} = C \rho_{k+1} + b_k \delta_k$
\quad Stop if $\|\rho_{k+1}\|$ or $\|\delta_{k+1}\|$ is sufficiently small.
end.\par}
\medskip

However this algorithm is unstable in finite precision arithmetic, as
we shall see in our main example below. One modification that
sucessfully avoids instability is to force the condition 
$$\rho_k \in \langle e \rangle^\perp, \eqno{(\chapno3.8)}$$
to hold for all $k$.
Now Lemma \chapno3.2 implies the existence of exactly one vector
$\rho_k \in \langle e \rangle^\perp$ for which $P \rho_k = r_k$. Therefore,
defining $Q$ to be the orthogonal projection onto $\langle e \rangle^\perp$, that
is $Q : x \mapsto x - e (e^T x)/(e^T e)$, we obtain

\proclaim {Algorithm \chapno3.3b}. Choose any $\xi_0$ in $\langle e \rangle^\perp$. Set
$\rho_0 = Q(b - A \xi_0)$, $\delta_0 = C \rho_0$.

{\obeylines 
For $k=0, 1, 2, \ldots$ do begin
\quad $a_k = \rho_k^T C \rho_k / \delta_k^T A \delta_k$
\qquad $\xi_{k+1} = \xi_k + a_k \delta_k$
\qquad $\rho_{k+1} = Q(\rho_k - a_k A \delta_k)$
\quad $b_k = \rho_{k+1}^T C \rho_{k+1} / \rho_k^T C \rho_k$
\qquad $\delta_{k+1} = C \rho_{k+1} + b_k \delta_k$
\quad Stop if $\|\rho_{k+1}\|$ or $\|\delta_{k+1}\|$ is sufficiently small.
end.\par}
\medskip

We see that the only restriction on $C$ is that it must be a
non-negative definite symmetric matrix such that $\ker C = \langle e \rangle$. It is
easy to construct such a matrix given a positive definite symmetric
matrix D by adding a rank one matrix:
$$ C = D - {(De)(De)^T \over e^T De}. \eqno{(\chapno3.9)}$$
The Cauchy-Schwarz inequality implies that $x^T C x \ge 0$ with
equality if and only if $x \in \langle e \rangle$. Of course we do not need to form
$C$ explicitly, since $C : x \mapsto Dx - (e^T Dx / e^T De) De$.
Before constructing $D$ we consider the spectral properties of
$A_\infty = (\phi(j-k))_{j,k\in\ZZ}$ in more detail.

A minor
modification to Proposition 5.2.2 yields the following useful result.
We recall the definition of a zero-summing sequence from Definition
4.3.1 and that of the symbol function from (5.2.7).

\proclaim {Proposition \chapno3.4}. For every $\eta \in (0,2\pi)$ we
can find a set $\{ (y_j^{(n)})_{j\in\ZZ} : n = 1, 2, \ldots \}$ of
zero-summing sequences such that
$$ \lim_{n\to\infty} \sum_{j,k\in\ZZ} y_j^{(n)} y_k^{(n)} \phi(j-k)
\Bigl/ \sum_{j\in\ZZ} [y_j^{(n)}]^2 = \sigma(\eta). \eqno{(\chapno3.10)}$$

\pf We adopt the proof technique of Proposition 5.2.2. For each
positive integer $n$ we define the trigonometric polynomial
$$ L_n(\xi) = n^{-1/2} \sum_{k=0}^{n-1} e^{ik\xi}, \qquad \xi \in
\R1,$$
and we recall from Section 4.2 that 
$$
 K_n(\xi) = {\sin^2 n\xi/2 \over n \sin^2 \xi/2} = \left| L_n(\xi)
\right|^2, \eqno{(\chapno3.11)}$$
where $K_n$ is the $n$th degree Fej\'er kernel. We now choose
$(y^{(n)}_j)_{j\in\ZZ}$ to be the Fourier coefficients of the
trigonometric polynomial $\xi \mapsto L_n(\xi - \eta) \sin \xi/2$,
which implies the relation
$$ \Bigl| \sum_{j\in\ZZ} y^{(n)}_j e^{ij\xi} \Bigr|^2 = \sin^2 \xi/2\ 
K_n(\xi-\eta), $$
and we see that $(y_j^{(n)})_{j\in\ZZ}$ is a zero-summing sequence. By
the Parseval relation we have
$$ \sum_{j\in\ZZ} [y_j^{(n)}]^2 = (2\pi)^{-1} \int_0^{2\pi} \sin^2
\xi/2 \ K_n(\xi - \eta)\,d\xi \eqno{(\chapno3.12)}$$
and the approximate identity property of the Fej\'er kernel (Zygmund
(1988), p. 86) implies that
$$\eqalign{
 \sin^2 \eta/2 
  &= \lim_{n\to\infty} (2\pi)^{-1} \int_0^{2\pi} \sin^2
\xi/2 \ K_n(\xi - \eta)\,d\xi \cr
  &= \lim_{n\to\infty} \sum_{j\in\ZZ} [y_j^{(n)}]^2.
}$$
Further, because $\sigma$ is continuous on $(0,2\pi)$ (see Section
4.4), we have 
$$\eqalign{
 \sin^2 \eta/2 \ \sigma(\eta)
  &= \lim_{n\to\infty} (2\pi)^{-1} \int_0^{2\pi} \sin^2
\xi/2 \ K_n(\xi - \eta) \sigma(\xi)\,d\xi \cr
  &= \lim_{n\to\infty} \sum_{j,k\in\ZZ} y^{(n)}_j y^{(n)}_k \phi(j-k), 
}$$
the last line being a consequence of (4.3.6).~\qed

Thus we have shown that, just as in the classical theory of Toeplitz
operators (Grenander and Szeg\H{o} (1984)), everything depends on the
range of values of the symbol function $\sigma$. Because $\sigma$
inherits the double pole that $\phihat$ enjoys at zero, we have
$\sigma \colon (0,2\pi) \mapsto (\sigma(\pi), \infty)$. In Figure
\chapno2 we display the function $[0,2\pi] \ni \xi \mapsto 1/\sigma(\xi)$.

Now let $m$ be a positive integer and let $(d_j)_{j=-m}^m$ be an even
sequence of real numbers. We define a bi-infinite banded symmetric
Toeplitz matrix $D_\infty$ by the equations
$$ (D_\infty)_{jk} = \cases{ d_{j-k},  \qquad |j-k| \le m, \cr 0,
\qquad \hbox{ otherwise }.} \eqno{(\chapno3.13)}$$
Thus $(D_\infty A_\infty)_{jk} = \psi(j-k)$ where $\psi(x) =
\sum_{l=-m}^m d_l \phi(x-l)$. Further
$$ \sum_{j,k\in\ZZ} y_j y_k \psi(j-k) 
  = (2\pi)^{-1} \int_0^{2\pi} \Bigl| \sum_{j\in\ZZ} y_j e^{ij\xi}
\Bigr|^2 \sigma_{D_\infty}(\xi) \sigma(\xi)\,d\xi.
\eqno{(\chapno3.14)}$$
Now the function $\xi \mapsto \sigma_{D_\infty}(\xi) \sigma(\xi)$ is
continuous for $\xi \in (0,2\pi)$, so the argument of Proposition
\chapno3.4 also shows that, for every $\eta \in (0,2\pi)$, we can find
a set $\{ (y^{(n)}_j)_{j\in\ZZ} : n = 1, 2, \ldots\ \}$ of
zero-summing sequences such that
$$ \lim_{n\to\infty} \sum_{j,k\in\ZZ} y^{(n)}_j y^{(n)}_k \psi(j-k)
\Bigl/ \sum_{j\in\ZZ} [y^{(n)}_j]^2 = \sigma_{D_\infty}(\eta)
\sigma(\eta). \eqno{(\chapno3.15)}$$

\vfill
\eject
\centerline{\psfig{file=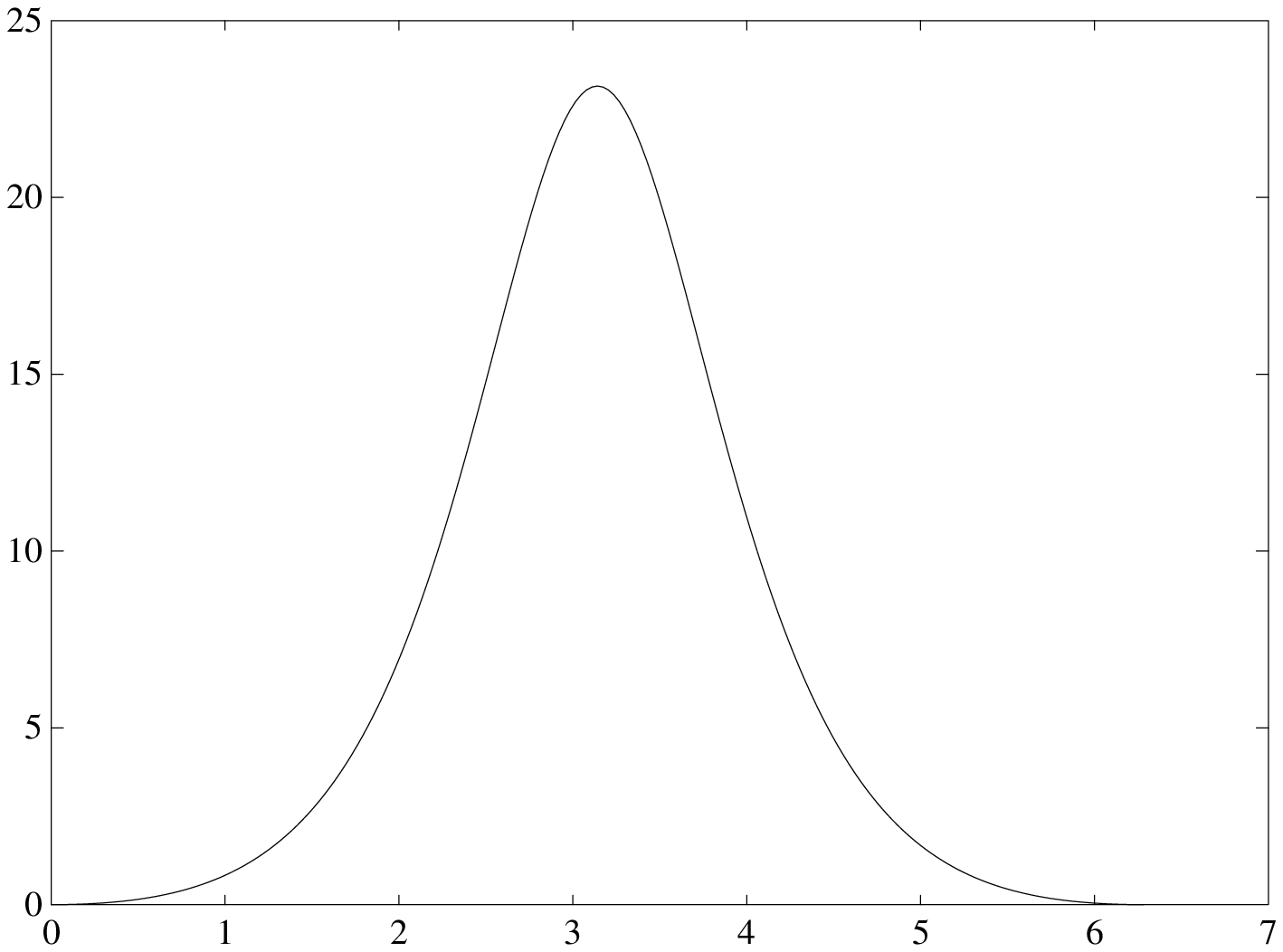}}
\centerline{{\bf Figure \chapno2:} The reciprocal symbol function $1/\sigma$ for
the multiquadric.}

A good preconditioner must ensure that $\{
\sigma_{D_\infty}(\xi) \sigma(\xi) : \xi \in (0,2\pi) \}$ is a bounded
set. Because of the form of $\sigma_{D_\infty}$ we have the equation
$$ \sum_{j=-m}^m d_j = 0. \eqno{(\chapno3.16)}$$
Moreover, as in Section 6.2, we want the approximation
$$ \sigma_{D_\infty}(\xi) \sigma(\xi) \approx 1, \qquad \xi \in
(0,2\pi), \eqno{(\chapno3.17)}$$
and we need $\sigma_{D_\infty}$ to be a non-negative trigonometric
polynomial which is positive almost everywhere, which ensures that
every one of its principal minors is positive definite.

Recalling Theorem 9 of Buhmann and Micchelli (1991), we let
$$ c^{(n)}_j = - \left( A^{-1}_n \right)_{j0}, \qquad j = -m, \ldots, m,
\eqno{(\chapno3.18)}$$
and to subtract a multiple of the vector $[1, \ldots, 1]^T \in
\R1^{2m+1}$ from $(c^{(n)}_j)_{j=-m}^m$ to form a new vector
$(d_j)_{j=-m}^m$ satisfying $\sum_{j=-m}^m d_j = 0$. Recalling that
$c^{(n)}_j \approx \gamma_j$ for suitable $m$ and $n$, where
$$ \sigma^{-1}(\xi) = \sum_{j\in\ZZ} \gamma_j e^{ij\xi}, \qquad \xi
\in \R1, \eqno{(\chapno3.19)}$$
and $\sum_{j \in \ZZ} \gamma_j = 0$ (since $\sigma$ inherits the double pole
of $\phihat$ at zero), we hope to achieve (\chapno3.17). Fortunately,
in several cases, 
we find that $\sigma_{D_\infty}$ is negative on $(0,2\pi)$, so that
$\sigma_{D_\infty}$ needs no further modifications. Unfortunately we
cannot explain this lucky fact at present, but perhaps one should not
always look a mathematical gift horse in the mouth. Therefore let
$n=64$ and $m=9$. Direct calculation yields
$$ \pmatrix{ c_0 \cr c_1 \cr \vdots \cr c_9 } 
  = - \pmatrix{    
  -6.8219 \times 10^0 \cr
  \ \,\, 4.9588 \times 10^0\cr
  -2.0852 \times 10^0\cr
  \ \,\, 7.2868 \times 10^{-1}\cr
  -2.5622 \times 10^{-1}\cr
  \ \,\, 8.8267 \times 10^{-1}\cr
  -3.1071 \times 10^{-2}\cr
  \ \,\, 1.0626 \times 10^{-2}\cr
  -3.7923 \times 10^{-3}\cr
  \ \,\, 1.2636\times 10^{-3}\cr},
\eqno{(\chapno3.20)}$$
and we then obtain
$$ - \pmatrix{ c_0 \cr c_1 \cr \vdots \cr c_9 } 
  = \pmatrix{    
  -6.8220 \times 10^0 \cr
  \ \,\, 4.9587 \times 10^0\cr
  -2.0852 \times 10^0\cr
  \ \,\, 7.2863 \times 10^{-1}\cr
  -2.5626 \times 10^{-1}\cr
  \ \,\, 8.8224 \times 10^{-1}\cr
  -3.1113 \times 10^{-2}\cr
  \ \,\, 1.0583 \times 10^{-2}\cr
  -3.8350 \times 10^{-3}\cr
  \ \,\, 1.2210\times 10^{-3}\cr}.
\eqno{(\chapno3.21)}$$
Figures \chapno3 and \chapno4 display the functions
$\sigma_{D_\infty}$ and $\xi \mapsto \sigma_{D_\infty}(\xi) / \sin^2
(\xi/2)$ on the domain $[0,2\pi]$ respectively. The latter is clearly
a positive function, which implies that the former is positive on the
open interval $(0,2\pi)$.

Thus, given 
$$ A_N = \Bigl( \phi(j-k) \Bigr)_{j,k=-N}^N$$
for any $N \ge n$, we let $D_N$ be any $(2N+1) \times (2N+1)$ principal
minor of $D_\infty$ and define the preconditioner $C_N$ by the equation
$$ C_N = D_N - {(D_N e) (D_N e)^T \over e^T D_N e},
\eqno{(\chapno3.22)}$$
where $e = [1, \ldots, 1]^T \in \R1^{2N+1}$. We reiterate that we
actually compute the matrix-vector product $C_N x$ by the operations
$x \mapsto D_N x - (e^T D_N x/e^T D_N e) e$ rather than by storing the
elements of $C_N$ in memory.

$C_N$ provides an excellent preconditioner. Tables 6.3 and 6.4
illustrate its use when Algorithm 6.3.3b is applied to the linear
system
$$\eqalign{ A_N x + e y &= b, \cr
            e^T x &= 0,\cr} \eqno{(\chapno3.23)}$$
when $N=2,048$ and $N = 32,768$ respectively.
Here $y \in \R1$, $e = [1, \ldots, 1]^T \in \R1^{2N+1}$ and $b \in
\R1^{2N+1}$ consists of pseudo-random real numbers uniformly
distributed in the interval $[-1,1]$. Again, this behaviour is typical
and all our numerical experiments indicate that the number of steps is
independent of $N$. We remind the reader that the error shown is
$\|\rho_{k+1}\|$, but that the iterations are
stopped when either $\|\rho_{k+1}\|$ or $\|\delta_{k+1}\|$ is less
than $10^{-13} \|b\|$, where we are using the notation of Algorithm 6.3.3b.

It is interesting to compare Table 6.3 with Table 6.5.
Here we have chosen $m=1$, and the preconditioner is essentially a multiple
of the second
divided difference preconditioner advocated by Dyn, Levin and Rippa
(1986). Indeed, we find that $d_0 = 7.8538$ and $d_1 = d_{-1} =
-3.9269$. 
We see that its behaviour is clearly inferior to the preconditioner
generated by choosing $m=9$. Furthermore, this is to be expected,
because we are choosing a smaller finite section to approximate the
reciprocal of the symbol function. However, because
$\sigma_{D_\infty}(\xi)$ is a multiple of $\sin^2 \xi/2$, this preconditioner
still possesses the property that $\{ \sigma_{D_\infty}(\xi)
\sigma(\xi) : \xi \in (0,2\pi) \}$ is a bounded set of real numbers.

$$\vbox{
\settabs\+\indent&c\qquad\hskip 24pt&Optimal bound\quad&\cr
\+&{\rm Iteration}&{\rm Error}\cr
\+& $1$& $3.975553 \times 10^{4}$\cr
\+& $2$& $8.703344 \times 10^{-1}$\cr
\+& $3$& $2.463390 \times 10^{-2}$\cr
\+& $4$& $8.741920 \times 10^{-3}$\cr
\+& $5$& $3.650521 \times 10^{-4}$\cr
\+& $6$& $5.029770 \times 10^{-6}$\cr
\+& $7$& $1.204610 \times 10^{-5}$\cr
\+& $8$& $1.141872 \times 10^{-7}$\cr
\+& $9$& $1.872273 \times 10^{-9}$\cr
\+& $10$& $1.197310 \times 10^{-9}$\cr
\+& $11$& $3.103685 \times 10^{-11}$\cr

}$$
\centerline{{\bf Table \chapno3:} Preconditioned CG -- $m=9$, $n=64$,
$N=2,048$}

$$\vbox{
\settabs\+\indent&c\qquad\hskip 24pt&Optimal bound\quad&\cr
\+&{\rm Iteration}&{\rm Error}\cr
\+& $1$& $2.103778 \times 10^{5}$\cr
\+& $2$& $4.287497 \times 10^{0}$\cr
\+& $3$& $5.163441 \times 10^{-1}$\cr
\+& $4$& $1.010665 \times 10^{-1}$\cr
\+& $5$& $1.845113 \times 10^{-3}$\cr
\+& $6$& $3.404016 \times 10^{-3}$\cr
\+& $7$& $3.341912 \times 10^{-5}$\cr
\+& $8$& $6.523212 \times 10^{-7}$\cr
\+& $9$& $1.677274 \times 10^{-5}$\cr
\+& $10$& $1.035225 \times 10^{-8}$\cr
\+& $11$& $1.900395 \times 10^{-10}$\cr

}$$
\centerline{{\bf Table \chapno4:} Preconditioned CG -- $m=9$, $n=64$, $N=32,768$}
\smallskip

It is also interesting to compare the spectra of $C_n A_n$ for $n=64$
and $m= 1$ and $m = 9$. Accordingly, Figures \chapno5 and \chapno6
display all but the largest nonzero eigenvalues of $C_n A_n$ for $m=1$ and
$m=6$ respectively. The largest eigenvalues are $502.6097.  $ and
$288.1872$, respectively, and these were omitted from the plots
in order to reveal detail at smaller scales. We see that the clustering of the spectrum
when $m=9$ is excellent.

$$\vbox{
\settabs\+\indent&c\qquad\hskip 24pt&Optimal bound\quad&\cr
\+&{\rm Iteration}&{\rm Error}\cr
\+& $1$& $2.645008 \times 10^{4}$\cr
\+& $10$& $8.632419 \times 10^{0}$\cr
\+& $20$& $9.210298 \times 10^{-1}$\cr
\+& $30$& $7.695337 \times 10^{-1}$\cr
\+& $40$& $3.187051 \times 10^{-5}$\cr
\+& $50$& $5.061053 \times 10^{-7}$\cr
\+& $60$& $7.596739 \times 10^{-9}$\cr
\+& $70$& $1.200700 \times 10^{-10}$\cr
\+& $73$& $3.539988 \times 10^{-11}$\cr
\+& $74$& $1.992376 \times 10^{-11}$\cr

}$$
\centerline{{\bf Table \chapno5:} Preconditioned CG -- $m=1$, $n=64$, $N=8,192$}
\smallskip

The final topic in this section demonstrates the instability of
Algorithm 6.3.3a when compared with Algorithm 6.3.3b. We refer the
reader to Table 6.6, where we have chosen $m=9$,
$n=N=64$, and setting $b = [1, 4, 9, \ldots, N^2]^T$. The iterations
for Algorithm 6.3.3b, displayed in Table 6.7, were stopped at
iteration $108$. For Algorithm 6.3.3a, iterations were stopped when
either
$\|\rho_{k+1}\|$ or $\|\delta_{k+1}\|$ became smaller than $10^{-13}
\|b\|$. It is useful to display the norm of $\|\delta_k\|$ rather than
$\|\rho_k\|$ in this case. We see that the two algorithms almost agree
on the early interations, but that Algorithm 6.3.3a soon begins
cycling, and no convergence seems to occur. Thus when $\rho_k$ can
leave the required subspace due to finite precision arithmetic, it is
possible to attain non-descent directions.

$$\vbox{
\settabs\+\indent&c\qquad\hskip 24pt&Optimal bound\quad&&Optimal bound\quad&\cr
\+&{\rm Iteration}&{$\|\delta_k\|$ }-- 6.3.3a&{$\|\delta_k\|$ }-- 6.3.3b \cr
\+& $1$& $4.436896 \times 10^{4}$& $4.436896 \times 10^{4}$\cr
\+& $2$& $2.083079 \times 10^{2}$& $2.083079 \times 10^{2}$\cr
\+& $3$& $2.339595 \times 10^{0}$& $2.339595 \times 10^{0}$\cr
\+& $4$& $1.206045 \times 10^{-1}$& $1.206041 \times 10^{-1}$\cr
\+& $5$& $1.698965 \times 10^{-3}$& $1.597317 \times 10^{-3}$\cr
\+& $6$& $6.537466 \times 10^{-2}$& $6.512586 \times 10^{-2}$\cr
\+& $7$& $1.879294 \times 10^{-4}$& $9.254943 \times 10^{-6}$\cr
\+& $8$& $2.767714 \times 10^{-2}$& $1.984033 \times 10^{-7}$\cr
\+& $9$& $3.453789 \times 10^{-4}$\cr
\+& $10$& $1.914126 \times 10^{-3}$\cr
\+& $20$& $4.628447 \times 10^{-1}$\cr
\+& $30$& $3.696474 \times 10^{-0}$\cr
\+& $40$& $8.061922 \times 10^{+3}$\cr
\+& $50$& $2.155310 \times 10^{0}$\cr
\+& $100$& $3.374467 \times 10^{-1}$\cr
\+& $101$& $1.121903 \times 10^{0}$\cr
\+& $102$& $1.920517 \times 10^{-1}$\cr
\+& $103$& $3.772007 \times 10^{-2}$\cr
\+& $104$& $3.170231 \times 10^{-2}$\cr
\+& $105$& $2.612073 \times 10^{-1}$\cr
\+& $106$& $2.236274 \times 10^{0}$\cr
\+& $107$& $8.875137 \times 10^{-1}$\cr
\+& $108$& $1.823607 \times 10^{-1}$\cr

}$$
\centerline{{\bf Table \chapno5:} Algorithms  6.3.3a \& b -- $m=1$, $n=64$,
$N=64$, $b = [1, 4, \ldots, N^2]^T$. }
\smallskip

\vfill
\eject
\centerline{\psfig{file=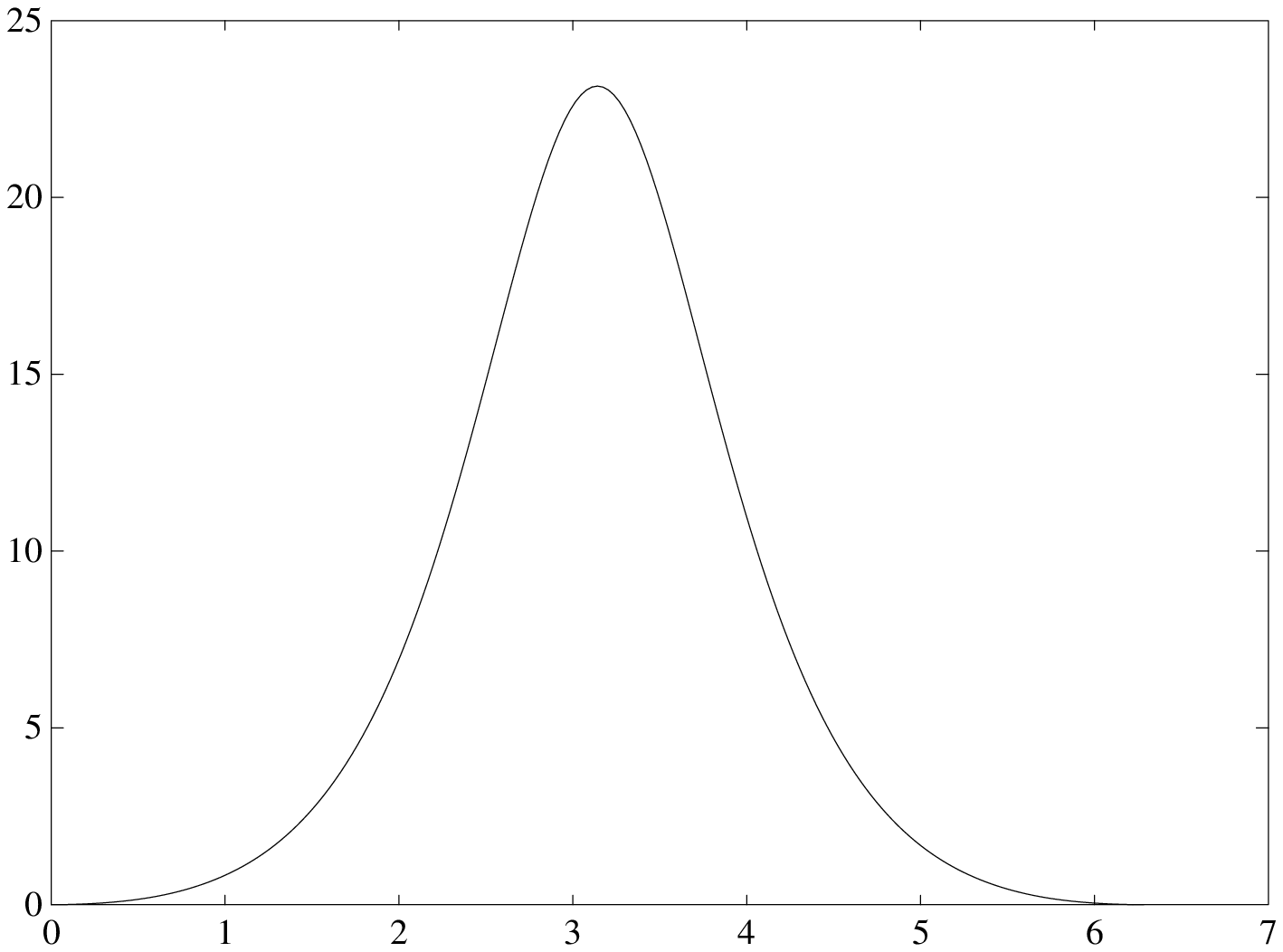}}
\centerline{{\bf Figure \chapno3:} The function $\sigma_{D_\infty}$. }

\vfill
\eject

\vfill
\eject
\centerline{\psfig{file=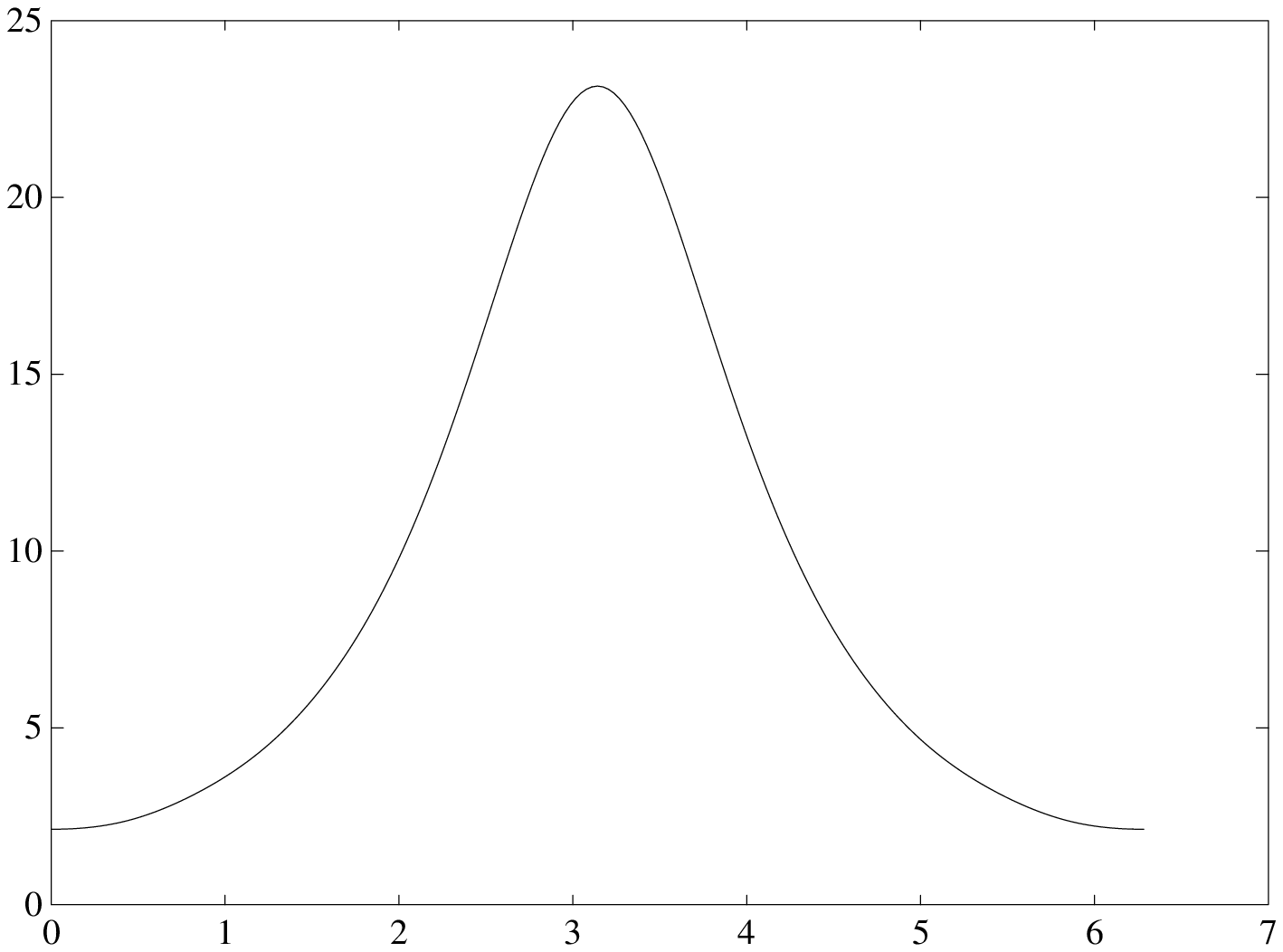}}
\centerline{{\bf Figure \chapno4:} The function $\xi \mapsto 
\sigma_{D_\infty}(\xi)/\sin^2(\xi/2)$.}

\vfill
\eject

\vfill
\eject
\centerline{\psfig{file=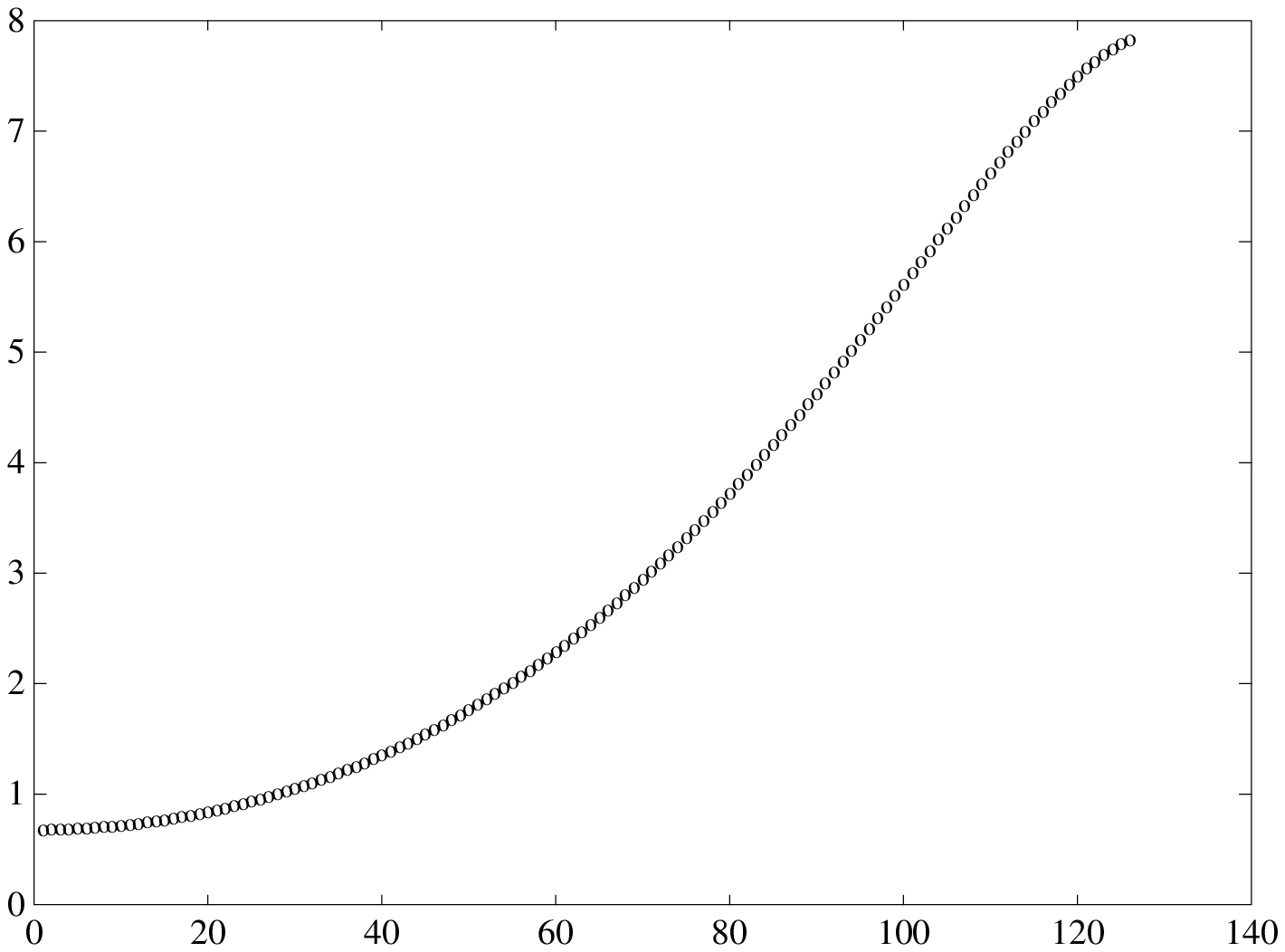}}
\centerline{{\bf Figure \chapno5:} The spectrum of $C_n A_n$ for $m=1$
and $n=64$.}

\vfill
\eject
\centerline{\psfig{file=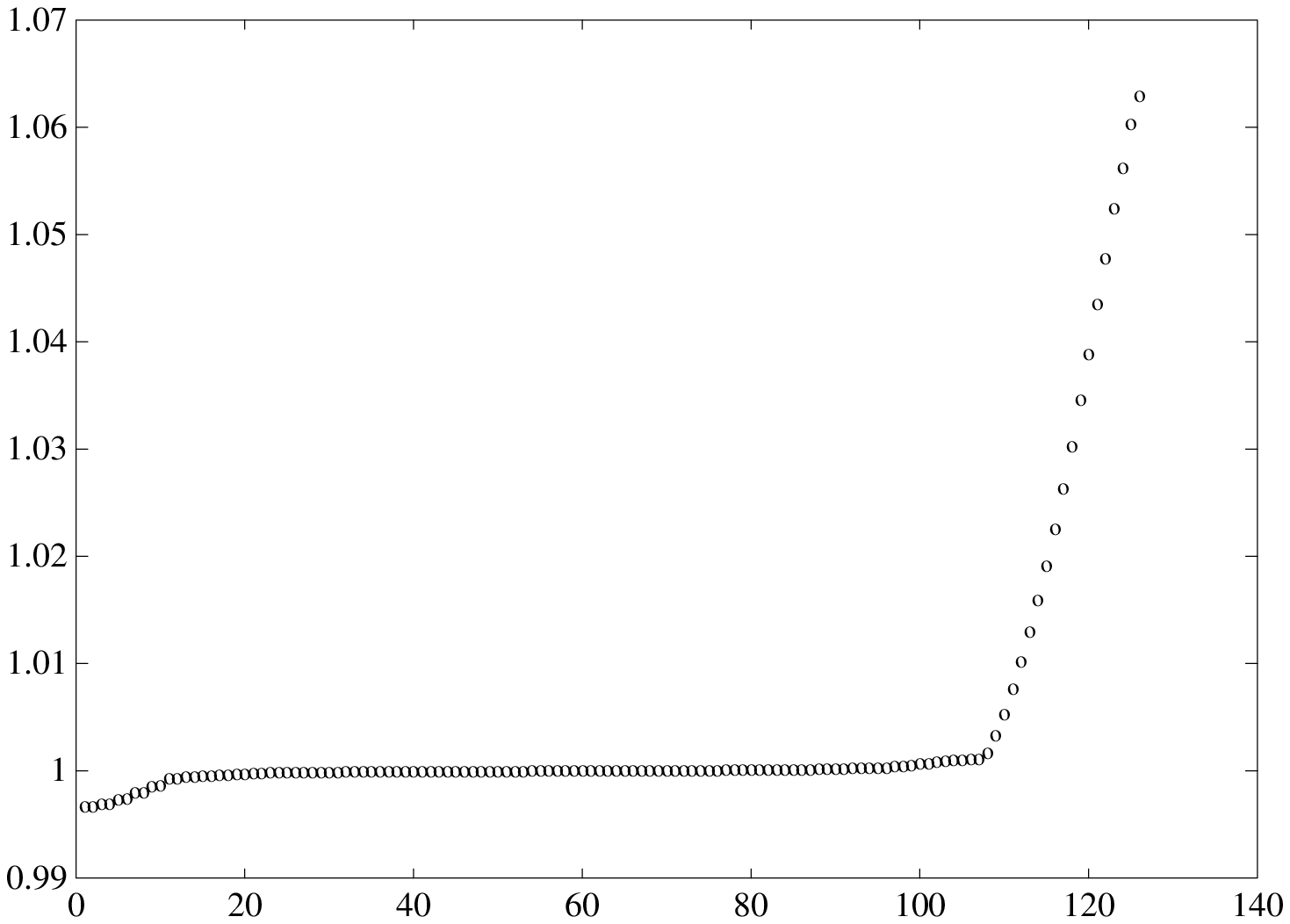}}
\centerline{{\bf Figure \chapno6:} The spectrum of $C_n A_n$ for
$m=9$ and $n=64$.}

\vfill\eject
\vfill\eject
\def\chapno{7.}
\headline{\ifnum\pageno=118\hfil\else{\headfont\hfil\chptitle\hfil}\fi}
\def\chptitle{On the asymptotic cardinal function for the multiquadric}
\bigskip
\centerline{\bf 7 : On the asymptotic cardinal function for the multiquadric}


\def\pts{{(x_j)_{j=1}^n}}

\def\sjk{{\sum_{j,k=1}^n y_j y_k\,}}
\def\djk{{\| x_j - x_k \|}}
\def\xjk{{x_j-x_k}}
\def\yjk{{|\sum_{j=1}^n y_j e^{ix_j t}|^2}}
\def\trig{{|\sum_{j=1}^n y_j e^{ijt}|^2}}

\def\skzd{{\sum_{k\in\Zd}}}

\def\rint{{\int_{-\infty}^{\infty}}}
\def\hint{{\int_0^\infty}}

\def\Sdm1{{S^{d-1}}}

\def\ysqr{{\Vert y \Vert^2}}

\def\ainv{{A_n^{-1}}}
\def\aninv{{\Vert A_n^{-1} \Vert_2}}

\def\zsum{{\sum_{k=-\infty}^\infty}}

\def\hi{{(2\pi)^{-1}}}
\def\tint{{\hi\int_0^{2\pi}}}
\def\tintd{{(2\pi)^{-d}\int_{\Td}}}
\def\fhat{{\hat f}}
\def\ajxi{{\Bigl|\sum_{j\in\Zd} a_j \exp(ij\xi)\Bigl|^2}}
\def\Rdm0{{\Rd\setminus\{0\}}}
\def\ghat{{\hat g}}
\def\ytrigo{{\Bigl|\sum_{j\in\Zd} y_j \exp(ix_j \xi)\Bigl|^2}}
\def\ytrig2{{\Bigl|\sum_{j\in\Zd} y_j \exp(ij \xi)\Bigl|^2}}
\def\ynj{{y_j^{(n)}}}

\def\phichat{{\phihat_c}}
\def\chichat{{{\hat \chi}_c}}

\def\min{{\rm min}}

\def\fhat{{\hat f}}
\def\phic{{\phi_c}}
\def\phichat{{\phihat_c}}
\def\chic{{\chi_c}}
\def\chichat{{{\hat \chi}_c}}

\sect {\chapno1. Introduction}
The radial basis function approach to interpolating a function
$f\colon \Rd \to \R1$ on the integer lattice $\Zd$ is as follows. Given a
continuous univariate function $\phi \colon [0,\infty) \to \R1$, we
seek a {\it cardinal function}
$$ \chi(x) = \sum_{j\in\Zd} a_j \phi(\|x-j\|), \qquad x \in \Rd,
\eqno{(\chapno1.1)}$$
that satisfies
$$ \chi(k) = \delta_{0,k}, \qquad k \in \Zd. $$
Therefore 
$$ If(x) = \sum_{j\in\Zd} f(j) \chi(x-j), \qquad x \in \Rd,
\eqno{(\chapno1.2)}$$
is an interpolant to $f$ on the integer lattice whenever (\chapno1.2) is well
defined. Here $\|\cdot\|$ is the Euclidean norm on $\Rd$. This
approach provides a useful and flexible family of approximants for
many choices of $\phi$, but
here we concentrate on the Hardy multiquadric $\phi_c(r) =
(r^2+c^2)^{1/2}$.
For this function, Buhmann (1990) has shown that a cardinal function $\chi_c$
exists and its Fourier tranform is given by the equation
$$ \chichat(\xi) = {\phichat(\|\xi\|) \over {\sum_{k\in\Zd}
\phichat(\|\xi+2\pi k\|)}}, \qquad \xi \in \Rd,
\eqno{(\chapno1.3)}$$
where $\{\phichat(\|\xi\|) : \xi \in \Rd\}$ is the generalized Fourier
transform of $\{\phic(\|x\|): x \in \Rd\}$.
Further, $\chic$ possesses a classical Fourier transform (see Jones (1982)
or Schwartz (1966)). 
In this chapter, we prove that $\chichat$ enjoys the following property:
$$ \lim_{c\to\infty} \chichat(\xi) = 
   \cases{1, \qquad \xi \in (-\pi,\pi)^d,\cr 0 ,\qquad \xi \notin
[-\pi,\pi]^d,}
\eqno{(\chapno1.4)}$$
which sheds new light on the approximation properties of the
multiquadric as $c \to \infty$. For example, in the case $d=1$, (\chapno1.4)
implies that $\lim_{c\to\infty} \chic(x) = \hbox{sinc}(x)$, providing
a perhaps unexpected link with sampling theory and the classical
theory of the Whittaker cardinal spline. 
Further, our work has links with the error analysis of Buhmann and
Dyn (1991) and illuminates the explicit calculation of Section 4 of
Powell (1991). It may also be compared with the results of Madych and
Nelson (1990) and Madych (1990), because these papers present analogous
results for polyharmonic cardinal splines.

\sect {\chapno2. Some properties of the multiquadric}
The generalized Fourier transform of $\phi_c$ is given by 
$$ \phichat(\|\xi\|) = -\pi^{-1} (2\pi c/\|\xi\|)^{(d+1)/2}
K_{(d+1)/2}(c\|\xi\|),
\eqno{(\chapno2.1)} $$
for nonzero $\xi \in \Rd$ (see Jones (1982)). Here $\{ K_\nu(r): r > 0 \}$ are the modified
Bessel functions, which are  positive and smooth in $\R1^+$, have a pole at the
origin, and decay exponentially (see Abramowitz and Stegun (1970)). There is
an integral representation for these modified Bessel functions
(Abramowitz and Stegun (1970), equation 9.6.23) which transforms (\chapno2.1)
into a highly useful formula for $\phichat$:
$$ \phichat(\|\xi\|) 
= - \lambda_d c^{d+1} \int_1^\infty \exp(-cx\|\xi\|) (x^2-1)^{d/2}\,dx
,
\eqno{(\chapno2.2)}$$
where $\lambda_d = \pi^{d/2} /\Gamma(1+d/2)$.
A simple consequence of (\chapno2.2) is the following lemma, which 
bounds the exponential decay of $\phichat$.

\proclaim {Lemma \chapno2.1}. If $\|\xi\| > \|\eta\| > 0$, then
$$ |\phichat(\|\xi\|)| \le \exp[-c(\|\xi\|-\|\eta\|)]\,|\phichat(\|\eta\|)|. $$

\pf Applying (\chapno2.2), we obtain 
$$\eqalign{
 |\phichat(\|\xi\|)| 
  &=  \lambda_d c^{d+1} \int_1^\infty
    \exp[-cx(\|\xi\|-\|\eta\|)] \exp(-cx\|\eta\|)\,(x^2-1)^{d/2}\,dx \cr 
  &\le  \exp(-c(\|\xi\|-\|\eta\|))\,|\phichat(\|\eta\|)|, 
}$$ 
providing the desired bound.~\qed

We now prove our main result. We let $I \colon \Rd \to \R1$ be the
{\it characteristic function} of the cube $[-\pi,\pi]^d$, that is
$$ I(\xi) = \cases{1, \qquad \xi \in [-\pi,\pi]^d,\cr 0 ,\qquad \xi \notin
[-\pi,\pi]^d.}$$

\proclaim {Proposition \chapno2.2}. Let $\xi$ be any fixed point of $\Rd$. We have
$$ \lim_{c\to\infty} \chichat(\xi) = I(\xi), $$
if $\|\xi\|_\infty \ne \pi$, that is $\xi$ does not lie in the
boundary of $[-\pi,\pi]^d$.

\pf First, suppose that $\xi \notin [-\pi,\pi]^d$. Then there exists a nonzero integer
$k_0$ such that $\|\xi+2\pi k_0\| < \|\xi\|$, and Lemma \chapno2.1 provides the
bounds
$$\eqalign{
 |\phichat(\|\xi\|)| 
 &\le     \exp[-c(\|\xi\|-\|\xi+2\pi k_0\|)] |\phichat(\|\xi+2\pi k_0\|)| \cr 
 &\le  \exp[-c(\|\xi\|-\|\xi+2\pi k_0\|)] \sum_{k\in\Zd} |\phichat(\|\xi+2\pi k\|)| .
}$$
Thus, applying (\chapno1.3) and remembering that $\phichat$ does not change sign, we have 
$$ 0 \le \chichat(\xi) \le \exp[-c(\|\xi\|-\|\xi+2\pi k_0\|)], \qquad
\xi \notin [-\pi,\pi]^d. \eqno{(\chapno2.3)}$$
The upper bound of (\chapno2.3) converges to zero as $c\to\infty$,
which completes the proof for this range of $\xi$.

Suppose now that $\xi \in (-\pi,\pi)^d$. Further, we shall assume that
$\xi$ is nonzero, because we know that $\chichat(0) = 1$ for all
values of $c$.
Then $\|\xi+2\pi k\| > \|\xi\|$, for every
nonzero integer $k \in \Zd$. Now (\chapno1.3) provides the expression
$$
 \chichat(\xi) 
 =  \Bigl( 1 + \sum_{k\in\Zd\setminus\{0\}} 
    \Bigl| {\phichat(\|\xi+2\pi k\|) \over \phichat(\|\xi\|)}\Bigr| \Bigr)^{-1}.
\eqno{(\chapno2.4)}$$
We shall show that 
$$ \lim_{c\to\infty}
   \sum_{k\in\Zd\setminus\{0\}} \left|{{\phichat(\|\xi+2\pi k\|)}\over{\phichat(\|\xi\|)}}\right|
    = 0, \qquad \xi \in (-\pi,\pi)^d, \eqno{(\chapno2.5)}$$
which, together with (\chapno2.4), implies that 
$\lim_{c\to\infty} \chichat(\xi) = 1$.

Now Lemma \chapno2.1 implies that
$$ \sum_{k\in\Zd\setminus\{0\}} \left|{\phichat(\|\xi+2\pi k\|)\over\phichat(\|\xi\|)}\right|
 \le \sum_{k\in\Zd\setminus\{0\}} \exp[-c(\|\xi+2\pi k\| - \|\xi\|)], \eqno{(\chapno2.6)} $$
and each term of the series on the right converges to zero as $c\to \infty$,
since $\|\xi+2\pi k\| > \|\xi\|$ for every nonzero integer $k$. Therefore we need only
deal with the tail of the series. Specifically, we derive the equation
$$ \lim_{c\to\infty} 
   \sum_{\|k\| \ge 2\|e\|} \exp[-c(\|\xi+2\pi k\| - \|\xi\|)] = 0 ,\eqno{(\chapno2.7)}$$
where $e = [1, 1, \ldots, 1]^T$. Now, if $\|k\| \ge 2\|e\|$, then
$$ \|\xi+2\pi k\| - \|\xi\| \ge 2\pi(\|k\|-\|e\|) \ge \pi\|k\|,$$
remembering that we have $\|\xi\| \le \pi \| e \|$.
Hence
$$ \sum_{\|k\| \ge 2\|e\|} \exp[-c(\|\xi+2\pi k\| - \|\xi\|)]
   \le \sum_{\|k\| \ge 2\|e\|} \exp(-\pi c \|k\|). \eqno{(\chapno2.8)} $$
It is a simple exercise to prove that the series
$\sum_{\|k\| \ge 2\|e\|} \exp(-\pi \|k\|)$ is convergent. Therefore, given any
$\epsilon > 0$, there exists a positive number $R \ge 1$ such that
$$ \sum_{\|k\| \ge 2R\|e\|} \exp(-\pi \|k\|) \le \epsilon .$$
Consequently, when $c \ge \lceil R \rceil$ we have the inequality
$$ \sum_{\|k\| \ge 2\|e\|} \exp(-\pi c \|k\|)
  \le \sum_{\|k\| \ge 2R\|e\|} \exp(-\pi \|k\|) \le \epsilon,$$
which establishes (\chapno2.5). The proof is complete.~\qed

\def\fhat{{\hat f}}
\def\phichat{{\phihat_c}}
\def\chichat{{{\hat \chi}_c}}
\def\Td{{[-\pi,\pi]^d}}

\sect {\chapno3. Multiquadrics and entire functions of exponential
type $\bf \pi$}
\proclaimdefn {Definition \chapno3.1} Let $f \in L^2(\Rd)$. We shall say that $f$ is a {\it function
of exponential type $A$} if its Fourier transform $\fhat$ is supported by the
cube $[-A,A]^d$. We shall denote the set of all functions of exponential
type $A$ by $E_A(\Rd)$.

We remark that the Paley-Wiener theorem implies that $f$ may be
extended to an entire function on $\CC^d$ satisfying a certain growth
condition at infinity (see Stein and Weiss (1971), pages 108ff),
although we
do not need this result.
 
\proclaim {Lemma \chapno3.2}. Let $f \in E_\pi(\Rd) \cap L^2(\Rd)$ be a
continuous function. Then we have the equation
$$ \sum_{k\in\Zd} \fhat(\xi+2\pi k) = \sum_{k\in\Zd} f(k) \exp(-ik\xi)
, \eqno{(\chapno3.1)}$$
the second series being convergent in $L^2(\Rd)$.

\pf Let 
$$g(\xi) = \sum_{k\in\Zd} \fhat(\xi+2\pi k) ,  \qquad \xi \in \Rd.$$
At any point $\xi \in \Rd$, this series contains at most one nonzero
term, because of the condition on the support of $\fhat$. Hence $g$ is
well defined. Further, we have the relations
$$ \int_{[-\pi,\pi]^d} |g(\xi)|^2\,d\xi
  =  \int_\Rd |\fhat(\xi)|^2\,d\xi < \infty, $$
since the Parseval theorem implies that $\fhat$ is an element of
$L^2(\Rd)$. Thus $g \in L^2([-\pi,\pi]^d)$ and its Fourier series
$$g(\xi) = \sum_{k\in\Zd} g_k \exp(ik\xi), $$
is convergent in $L^2([-\pi,\pi]^d)$. The Fourier coefficients are
given by the expressions
$$
  g_k 
   =  (2\pi)^{-d} \int_{\Td} \fhat(\xi) \exp(-ik\xi)\,d\xi  
   =  (2\pi)^{-d} \int_\Rd \fhat(\xi) \exp(-ik\xi)\,d\xi  
   =  f(-k),
$$
where the final equation uses the Fourier inversion theorem for
$L^2(\Rd)$. The proof is complete.~\qed

We observe that an immediate consequence of the lemma is the convergence of the
series $\sum_{k\in\Zd} [f(k)]^2$, by the Parseval theorem. 

For the following results, we shall need the fact that $\chi_c \in
L^2(\Rd)$, which is a consequence of the analysis of Buhmann (1990).

\proclaim {Lemma \chapno3.3}. Let $f \in E_\pi(\Rd) \cap L^2(\Rd)$ be a
continuous function. For each
positive integer $n$, we define the function
$$ \widehat{S_c^n f}(\xi) = \left(\sum_{\|k\|_1 \le n} f(k) \exp(-ik\xi)
\right) \chichat(\xi), \qquad \xi \in \Rd. \eqno{(\chapno3.2)}  $$
Then $\{S_c^n f : n = 1, 2, \ldots \}$ forms a Cauchy sequence in $L^2(\Rd)$.

\pf Let
$Q_n\colon\Rd\to \R1$ be the trigonometric
polynomial
$$ Q_n(\xi) = \sum_{\|k\|_1\le n} f(k) \exp(-ik\xi), \eqno{(\chapno3.3)} $$
so that $\widehat{S_c^n f}(\xi) = Q_n(\xi)\chichat(\xi)$. 
It is a consequence of Lemma \chapno3.2 that this sequence of functions
forms a Cauchy sequence in $L^2(\Td)$. Indeed, we shall prove that
for $m \ge n$ we have
$$ \|\widehat{S_c^m f} - \widehat{S_c^n f} \|_{L^2(\Rd)}
      \le \|Q_m - Q_n\|_{L^2(\Td)}, \eqno{(\chapno3.4)} $$
so that the sequence of functions $\{\widehat{S_c^n f}: n = 1, 2, \ldots \}$ is
a Cauchy sequence in $L^2(\Rd)$. 

Now Fubini's theorem provides the relation
$$\eqalign{
 \|\widehat{S_c^m f} - \widehat{S_c^n f} \|_{L^2(\Rd)}^2
  &=  \int_\Rd \left|Q_m(\xi)-Q_n(\xi)\right|^2 \chichat^2(\xi)\,d\xi \cr 
  &=  \int_{\Td} 
         \left|Q_m(\xi)-Q_n(\xi)\right|^2
         \left( \sum_{l\in\Zd} \chichat^2(\xi+2\pi l)\right)\,d\xi .
}\eqno{(\chapno3.5)}$$ 
However, (\chapno1.3) gives the bound
$$\eqalign{
 \sum_{l\in\Zd} \chichat^2(\xi+2\pi l)  
   &=    {\sum_{l\in\Zd} \phichat^2(\|\xi+2\pi l\|)\Bigl/\!
       (\sum_{k\in\Zd} \phichat(\|\xi+2\pi k\|))^2} \cr 
   &\le  1,
}\eqno{(\chapno3.6)}$$
which, together with (\chapno3.5), yields inequality (\chapno3.4).~\qed

Thus we may define
$$ \widehat{S_c f} (\xi) = \chichat(\xi) \sum_{k\in\Zd} f(k) \exp(-ik\xi), \eqno{(\chapno3.7)}$$
and the series is convergent in $L^2(\Rd)$. Applying the inverse Fourier
transform term by term, we obtain the useful equation
$$ S_c f(x) = \sum_{k\in\Zd} f(k) \chi_c(x-k), \quad x \in \Rd .$$

\proclaim {Theorem \chapno3.4}. Let $f \in E_\pi(\Rd) \cap L^2(\Rd)$ be a
continuous function. We have
$$ \lim_{c \to \infty} S_c f(x) = f(x), $$
and the convergence is uniform on $\Rd$.

\pf We have the equation
$$ S_c f(x) - f(x)
 = (2\pi)^{-d} \int_\Rd \sum_{k\in\Zd}  \fhat(\xi+2\pi k) 
     \left(\chichat(\xi) - I(\xi) \right) \exp(ix\xi)\,d\xi .$$
Thus we deduce the bound
$$ |S_c f(x) - f(x) | \hbox{\hskip  8cm } $$
$$\eqalign{
   &\le   (2\pi)^{-d} \int_{\Td} |\fhat(\xi)|
            \sum_{k\in\Zd} \Bigl|\chichat(\xi+2\pi k) - I(\xi+2\pi k)\Bigr|\,d\xi 
             \hbox{\hskip  3cm }\cr 
   &=     (2\pi)^{-d} \int_{\Td} |\fhat(\xi)|  
            \left(1 - \chichat(\xi) + \sum_{k\in{\Zd\setminus\{0\}}}
                    \chichat(\xi+2\pi k) \right)\,d\xi ,
}\eqno{(\chapno3.8)}$$
using the fact that $\chichat$ is non-negative, and we observe that
this upper bound is independent of $x$. Therefore we prove that the
upper bound converges to zero as $c \to \infty$.

Applying (\chapno1.3), we obtain the relation
$$  \sum_{k\in{\Zd\setminus\{0\}}}  \chichat(\xi+2\pi k) = 1 -
\chichat(\xi) , \eqno{(\chapno3.9)}$$
whence
$$  \left|S_c f(x) - f(x) \right|
      \le 2  (2\pi)^{-d} \int_{\Td} |\fhat(\xi)|
(1-\chichat(\xi))\,d\xi . \eqno{(\chapno3.10)}$$
Now $\fhat \in L^2(\Td)$ implies $\fhat \in L^1([-\pi,\pi]^d)$,
by the Cauchy-Schwartz inequality. Further, Proposition \chapno2.2 gives the limit $\lim_{c\to\infty}
\chichat(\xi) = 1$, for $\xi \in (-\pi,\pi)^d$, and we have $0 \le 1 -
\chichat(\xi) \le 1$, by (\chapno1.3).
Therefore the dominated convergence theorem implies that
$$ \lim_ {c\to\infty}  
  (2\pi)^{-d} \int_{\Td} |\fhat(\xi)| (1-\chichat(\xi))\,d\xi               
  = 0. $$
The proof is complete.~\qed

\sect{\chapno4. Discussion}
Section 4 of Powell (1991) provides an explicit calculation that is
analogous to the proof of Theorem \chapno3.4 when $f(x)=x^2$. Of course, this
function does not satisfy the conditions of Theorem \chapno3.4. Therefore
extensions of this result are necessary, but the final form of the theorem is not
clear at present.

Theorem \chapno3.4 encourages the use of large $c$ for certain functions.
Indeed, it suggests that large $c$ will provide high accuracy
interpolants for univariate functions that are well approximated by integer
translates of the sinc function. Thus, in exact arithmetic, a large
value of $c$ should be useful whenever the function is well
approximated by the Whittaker cardinal series. However, we recall that
the linear systems arising when $c$ is large can be rather
ill-conditioned. Indeed, in Chapter 4 we proved that the smallest
eigenvalue of the interpolation matrix generated by a finite regular
grid converges to zero exponentially quickly as $c \to \infty$. We
refer the reader to Table 4.1 for further information.
Therefore special techniques are required for the effective use of large
$c$.

\vfill \eject
\def\chapno{8.}
\headline{\ifnum\pageno=126\hfil\else{\headfont\hfil\chptitle\hfil}\fi}
\def\chptitle{Conclusions}
\bigskip
\centerline{\bf 8 : Conclusions}
There seems to be no interest in using non-Euclidean norms for radial
basis functions at present, possibly because of the poor approximation
properties of the $\ell^1$-norm $\|\cdot\|_1$ reported by several
workers. Thus Chapter 2 does not seem to have any practical
applications yet. However, it may be useful to use $p$-norms
($1 < p < 2$), or functions of $p$-norms, when there is a known preferred
direction in the underlying function, because radial basis functions
based on the Euclidean norm can perform poorly in this context. On a
purely theoretical note, we observe that the construction of Section
2.4 can be applied to any norm enjoying the symmetries of the cube.

The greatest weakness --  and the greatest strength -- of the norm estimates
of Chapters 3--6 lies in their dependence on regular grids. However,
we note that the upper bounds on norms of inverses apply to sets of centres
which can be arbitrary subsets of a regular grid. In other words,
contiguous subsets of grids are not required. Furthermore, we
conjecture that a useful upper bound on the norm of the inverse
generated by an arbitrary set of centres with minimal separation
distance $\delta$ (that is $\|x_j - x_k\| \ge \delta > 0$ if $j \ne
k$) will be provided by the upper bound for the inverse generated by a
regular grid of spacing $\delta$.

Probably the most important practical finding of this dissertation is
that the number of steps required by the conjugate gradient algorithm
can be independent of the number of centres for suitable
preconditioners. We hope to discover preconditioners with this
property for arbitrary sets of centres.

The choice of constant in the multiquadric is still being investigated
(see, for instance, Kansa and Carlson (1992)). Because the approximation of
band-limited functions is of some practical importance, our findings
may be highly useful. In short, we suggest using as large a value of
the constant as the condition number allows. Hence there is some irony
in our earlier discovery that the condition number of the
interpolation matrix can increase exponentially quickly as the
constant increases.

Let us conclude with the remark that radial basis functions are
extremely rich mathematical objects, and there is much left to be
discovered. It is our hope that the strands of research 
initiated in this thesis will enable some of these future discoveries.

\vfill \eject
\headline{\ifnum\pageno=128\hfil\else{\headfont\hfil\chptitle\hfil}\fi}
\def\chptitle{References}
\bigskip
\centerline{\bf References}

\vskip 6pt\noindent
Abramowitz, M., and I. A. Stegun (1970), {\sl Handbook of Mathematical Functions},
Dover Publications (New York).

\vskip 6pt\noindent
Ball, K. M. (1989), Invertibility of Euclidean distance matrices and
radial basis interpolation, CAT report no. 201, Texas A \& M University,
College Station.

\vskip 6pt\noindent
Barrodale, I., M. Berkley and D. Skea (1992), ``Warping digital images
using thin plate splines'', presented at the Sixth Texas International
Symposium on Approximation Theory (Austin, January 1992).

\vskip 6pt\noindent
Baxter, B. J. C. (1991a), ``Conditionally positive functions and
$p$-norm distance matrices'', {\sl Constr. Approx.} {\bf 7}, pp.
427--440.

\vskip 6pt\noindent
Baxter, B. J. C. (1991b), ``Norm estimates for Toeplitz distance
matrices'', Report NA16, University of Cambridge.

\vskip 6pt\noindent
Baxter, B. J. C. (1992a), ``Norm estimates for inverses of distance
matrices'', in {\sl Mathematical Methods in Computer Aided Geometric
Design}, T. Lyche and L. L. Schumaker (eds.), Academic Press (New
York), pp. 9--18.

\vskip 6pt\noindent
Baxter, B. J. C. (1992b), ``On the asymptotic cardinal function of the
multiquadric $\phi(r) = (r^2 + c^2)^{1/2}$ as $c \to \infty$''. To
appear in a special issue of {\sl An International Journal of
Computers and Mathematics with Applications}.

\vskip 6pt\noindent
Baxter, B. J. C. (1992c), ``Norm estimates and preconditioning for
conjugate gradient solution of RBF linear systems''. In preparation. 

\vskip 6pt\noindent
Baxter, B. J. C. and C. A. Micchelli, ``Norm estimates for the $\ell^2$-inverses of
multivariate Toeplitz matrices''. In preparation.

\vskip 6pt\noindent
Beatson, R. K. and G. N. Newsam (1992), ``Fast evaluation of radial
basis functions: 1''. To
appear in a special issue of {\sl An International Journal of
Computers and Mathematics with Applications}.

\vskip 6pt \noindent
de Boor, C. (1987), ``Multivariate approximation'', in the {\sl The
State of the Art in Numerical Analysis}, A. Iserles and M. J. D.
Powell (eds.), Oxford University Press (Oxford), pp. 87--110.

\vskip 6pt\noindent
Buhmann, M. D. (1990), ``Multivariate cardinal interpolation with radial basis
functions'',
{\sl Constr. Approx.}~{\bf 6}, pp. 225--256.

\vskip 6pt\noindent
Buhmann, M. D. and M. J. D. Powell (1990), ``Radial basis function
interpolation on an infinite regular grid'', in {\sl Algorithms
for Approximation II}, J. C. Mason and M. G. Cox (eds.), Chapmann and
Hall (London), pp. 146--169.

\vskip 6pt\noindent
Buhmann, M. D. and N. Dyn (1991), ``Error estimates for multiquadric
interpolation'', in {\sl Curves and Surfaces}, P.-J. Laurent, A. Le
M\'ehaut\'e, and L. L. Schumaker (eds.), Academic Press (New York), pp. 51--58.

\vskip 6pt\noindent
Buhmann, M. D. and C. A. Micchelli (1991), ``Multiply monotone functions
for cardinal interpolation'', {\sl Advances in Applied Mathematics} {\bf 12} ,
pp. 358--386.

\vskip 6pt\noindent
Chui, C. K. (1988), {\sl Multivariate Splines}, SIAM Regional
Conference Series in Applied Mathematics {\bf 54} (Philadelphia).

\vskip 6pt\noindent
Chui, C. K. (1992), {\sl An Introduction to Wavelets}, Academic Press
(New York).

\vskip 6pt\noindent 
Davis, P. J. (1975), {\sl Interpolation and Approximation}, Dover
Publications (New York).

\vskip 6pt\noindent
Dahmen, W. and C. A. Micchelli (1983), ``Recent progress in
multivariate splines'', in {\sl Approximation Theory IV}, C. K. Chui,
L. L. Schumaker and J. D. Ward (eds), Academic Press (New York), pp. 27--121.

\vskip 6pt\noindent
Duchon, J. (1975), ``Fonctions-spline du type plaque mince en
dimension 2'', Technical Report 231, Universit\'e de Grenoble.

\vskip 6pt\noindent
Duchon, J. (1976), ``Fonctions-spline \`a energie invariante par
rotation'', Technical Report 27, Universit\'e de Grenoble.

\vskip 6pt\noindent
Dyn, N., W. A. Light and E. W. Cheney (1989), ``Interpolation by
piecewise linear radial basis functions'', {\sl J.~Approx.
Theory.}~{\bf 59}, pp. 202--223. 

\vskip 6pt\noindent
Dyn, N., D. Levin and S. Rippa (1986), ``Numerical procedures for
surface fitting of scattered data by radial functions'', {\sl SIAM J.
Sci. Stat. Comput.} {\bf 7}, pp. 639--659.

\vskip 6pt\noindent
Dyn, N., D. Levin and S. Rippa (1990), ``Data dependent triangulations
for piecewise linear interpolation'', {\sl IMA J. of Numer. Anal.}
{\bf 10}, pp. 137--154.

\vskip 6pt\noindent
Edrei, A. (1953), ``On the generating function of
doubly infinite, totally positive sequences'', {\sl Trans. Amer. Math. Soc.}
{\bf 74}, pp. 367--383.

\vskip 6pt\noindent
Fletcher, R. (1987), {\sl Practical Methods of Optimization}, Wiley \&
Sons (Chicester).

\vskip 6pt\noindent
Franke, R. (1982), ``Scattered data interpolation: tests of some
methods'', {\sl Math. of Comp.} {\bf 38}, pp. 181--200.

\vskip 6pt \noindent
Franke, R. (1987), ``Recent advances in the approximation of surfaces
from scattered data'', in {\sl Topics in Multivariate Approximation},
C. K. Chui, L. L. Schumaker and F. I. Utreras (eds), Academic Press
(New York), pp. 79--98.

\vskip 6pt\noindent
Golub, G. H. and C. F. Van Loan (1989), {\sl Matrix Computations}, The
John Hopkins University Press (Baltimore).

\vskip 6pt\noindent
Grenander, U. and G. Szeg\H{o} (1984),  {\sl Toeplitz Forms}, Chelsea (New York).

\vskip 6pt\noindent
Hardy, G. H., J. E. Littlewood and G. P\'olya (1952), {\sl Inequalities},
Cambridge University Press (Cambridge).

\vskip 6pt \noindent
Hayes, J. G. (1987), ``Fitting surfaces to data'', in {\sl The
Mathematics of Surfaces}, R. R. Martin (ed.), Oxford University Press
(Oxford), pp. 17--38.

\vskip 6pt\noindent
Hille, E. (1962), {\sl Analytic Function Theory, Volume II}, Ginn and
Co. (Waltham, Massachusetts).

\vskip 6pt \noindent
Jackson, I. R. H. (1988), {\sl Radial basis function methods for
multivariable approximation}, Ph.D. Dissertation, University of Cambridge.

\vskip 6pt\noindent
Jones, D. S. (1982), {\sl The Theory of Generalised Functions}, 
Cambridge University Press (Cambridge).

\vskip 6pt\noindent
Kansa, E. J. and R. E. Carlson (1992), ``Improved accuracy of
multiquadric interpolation using variable shape parameters''. To
appear in a special issue of {\sl An International Journal of
Computers and Mathematics with Applications}.

\vskip 6pt\noindent
Karlin, S. (1968), {\sl Total Positivity, Volume I}, Stanford University
Press (Stanford, California).

\vskip 6pt\noindent
Katznelson, Y. (1976), {\sl An Introduction to Harmonic Analysis},
Dover Publications (New York).

\vskip 6pt\noindent
Light, W. A. and E. W. Cheney (1986), {\sl Approximation Theory in
Tensor Product Spaces}, Lecture Notes in Mathematics {\bf 1169},
Springer Verlag (Berlin).

\vskip 6pt\noindent
Madych, W. R.  (1990), ``Polyharmonic splines, multivariate analysis and
entire functions'', {\sl International Series of Numerical Analysis} {\bf
94}, pp. 205--216.

\vskip 6pt\noindent
Madych, W. R. and S. A. Nelson (1990), ``Polyharmonic cardinal
splines'', {\sl J. Approx. Theory} {\bf 60}, pp. 141--156.

\vskip 6pt\noindent
Micchelli, C. A. (1986), ``Interpolation of scattered data:
distance matrices and conditionally positive functions'', {\sl Constr.
Approx.} {\bf 2}, pp. 11-22.

\vskip 6pt\noindent
Narcowich, F. J. and J. D. Ward (1990), ``Norm estimates for inverses
of scattered data interpolation matrices associated with completely
monotonic functions'', preprint.

\vskip 6pt\noindent
Narcowich, F. J. and J. D. Ward (1991), ``Norms of inverses and
condition numbers of matrices associated with scattered data'',
{\sl J. Approx. Theory} {\bf 64}, pp. 69--94.

\vskip 6pt\noindent
von Neumann, J. and I. J. Schoenberg (1941), ``Fourier integrals
and metric geometry'', {\sl Trans. Amer. Math. Soc.} {\bf 50}, pp. 226--251.

\vskip 6pt \noindent
Powell, M. J. D. (1981), {\sl Approximation Theory and Methods},
Cambridge University Press (Cambridge).

\vskip 6pt\noindent
Powell, M. J. D. (1991), ``Univariate multiquadric interpolation: some
recent results'',in {\sl Curves and Surfaces}, P.-J. Laurent, A. Le
M\'ehaut\'e, and L. L. Schumaker (eds.), Academic Press (New York), pp. 371--381. 

\vskip 6pt \noindent
Powell, M. J. D. Powell (1992), ``The theory of radial basis function
approximation in 1990'', in {\sl Advances in Numerical Analysis II:
Wavelets, Subdivision Algorithms and Radial Functions}, ed. W. A.
Light, Oxford University Press (Oxford), pp. 105--210.

\vskip 6pt\noindent
Rudin, W. (1973), {\sl Functional Analysis}, McGraw Hill (London).

\vskip 6pt\noindent
Schoenberg, I. J. (1935), ``Remarks to Maurice Fr\'echet's
article `Sur la definition d'une classe d'espace distanci\'es
vectoriellement applicable sur l'espace d'Hilbert.' '', {\sl Ann. of
Math.}~{\bf 36}, pp. 724-732.

\vskip 6pt\noindent
Schoenberg, I. J. (1937), ``On certain metric spaces arising
from Euclidean space by a change of metric and their embedding in
Hilbert space'', {\sl Ann. of Math.} {\bf 38}, pp. 787-793.

\vskip 6pt\noindent
Schoenberg, I. J. (1938), ``Metric spaces and completely
monotone functions'',  Ann. of Math. {\bf 39}, pp. 811-841.

\vskip 6pt\noindent
Schoenberg, I. J. (1951), ``On P\'olya frequency functions. I. The
totally positive functions and their Laplace transforms,'' {\sl J.
Analyse Math.} {\bf 1}, pp. 331--374.

\vskip 6pt\noindent
Schoenberg, I. J. (1973), {\sl Cardinal Spline Interpolation},
SIAM Regional
Conference Series in Applied Mathematics {\bf 12} (Philadelphia).

\vskip 6pt\noindent
Schwartz, L. (1966), {\sl Th\'eorie des Distributions}, Hermann (Paris).

\vskip 6pt\noindent
Stein, E. M. and G. Weiss (1971), {\sl Introduction to Fourier analysis on
Euclidean spaces}, Princeton Univ. Press (Princeton, New Jersey).

\vskip 6pt\noindent
Sun, X. (1990), ``Norm estimates for inverses of Euclidean distance
matrices'', preprint.

\vskip 6pt\noindent
Whittaker, E. T. and G. N. Watson (1927), {\sl A Course of Modern Analysis}, 
Cambridge University Press (Cambridge).

\vskip 6pt\noindent
Wilkinson, J. H. (1965), {\sl The Algebraic Eigenvalue Problem},
Oxford University Press (Oxford).

\vskip 6pt\noindent
Zygmund, A. (1979), {\sl Trigonometric Series},  Volumes I and II,  
Cambridge University Press (Cambridge).


\vfill \eject
\bye